\documentclass[12pt,fullpage]{article}

\usepackage{amsmath,amsfonts,amssymb,amsthm,bbm,graphics,psfrag}
\usepackage{xr,graphicx,epsfig,subfigure,color,soul, paralist}
\usepackage[active]{srcltx}
\def\R{\mathbb{R}}
\def\N{\mathbb{N}}

\def\S{\mathcal{S}}
\def\epsilon{\varepsilon}
\def\hat{\widehat}
\def\tilde{\widetilde}
\def\div{\mbox{div}}

\def\fin{$\Box$}

\newcommand{\SE}{\setcounter{equation}{0} \section}
\newcommand{\be}{\begin{equation}}
\newcommand{\ee}{\end{equation}}
\newcommand{\baa}{\begin{array}}
\newcommand{\eaa}{\end{array}}
\newcommand{\ba}{\begin{eqnarray}}
\newcommand{\ea}{\end{eqnarray}}

\newtheorem{theo}{\bf Theorem}[section]
\newtheorem{lem}[theo]{\bf Lemma}
\newtheorem{pro}[theo]{\bf Proposition}
\newtheorem{cor}[theo]{\bf Corollary}

\newtheorem{rem}[theo]{\bf Remark}

\begin{document}
\date{}
\title{\bf{Comparison results for semilinear elliptic equations using a new symmetrization method}}
\author{Fran\c cois Hamel$^{\hbox{\small{ a}}}$ and Emmanuel Russ$^{\hbox{\small{ b}}}$\thanks{This work has been carried out in the framework of the Labex Archim\`ede (ANR-11-LABX-0033) and of the A*MIDEX project (ANR-11-IDEX-0001-02), funded by the ``Investissements d'Avenir" French Government program managed by the French National Research Agency (ANR). The research leading to these results has received funding from the French ANR within the project PREFERED and from the European Research Council under the European Union's Seventh Framework Programme (FP/2007-2013) / ERC Grant Agreement n.321186 - ReaDi - Reaction-Diffusion Equations, Propagation and Modelling. Part of this work was also carried out during visits by the first author to the Departments of Mathematics of the University of California, Berkeley and of Stanford University, the hospitality of which is thankfully acknowledged.}\\
\\
\footnotesize{$^{\hbox{a }}$Aix Marseille Universit\'e, CNRS, Centrale Marseille}\\
\footnotesize{Institut de Math\'ematiques de Marseille, UMR 7373, 13453 Marseille, France}\\
\footnotesize{\& Institut Universitaire de France}\\
\footnotesize{$^{\hbox{b }}$ Universit\'e Joseph Fourier, Institut Fourier}\\
\footnotesize{100 rue des maths, BP 74, 38402 Saint-Martin d'H\`eres Cedex, France}}
\maketitle

\begin{abstract}
In this paper, we prove some pointwise comparison results between the solutions of some second-order semilinear elliptic equations in a domain $\Omega$ of $\R^n$ and the solutions of some radially symmetric equations in the equimeasurable ball $\Omega^*$. The coefficients of the symmetrized equations in~$\Omega^*$ satisfy similar constraints as the original ones in~$\Omega$. We consider both the case of equations with linear growth in the gradient and the case of equations with at most quadratic growth in the gradient. Lastly, we show some improved quantified comparisons when the original domain is not a ball. The method is based on a symmetrization of the second-order terms.
\end{abstract}

\tableofcontents


\SE{Introduction and main results}\label{intro}

Throughout all the paper, $n\ge 1$ is a given integer, $\Omega$ is a bounded domain of $\R^n$ of class~$C^2$ and $\Omega^{\ast}$ denotes the open Euclidean ball centered at $0$ such that
$$\left\vert\Omega^{\ast}\right\vert=\left\vert \Omega\right\vert,$$
where by a domain we mean a non-empty open connected subset of $\R^n$ and, for any measurable subset $E\subset\R^n$, $\left\vert E\right\vert$ stands for the Lebesgue measure of $E$.

We consider the following problem: given a bounded domain $\Omega\subset\R^n$ of class $C^2$ and given a positive solution $u$ of a quasilinear elliptic equation of the type
\be\label{equ}\left\{\baa{rcl}
-\div(A(x)\nabla u)+H(x,u,\nabla u) & = & 0\ \hbox{ in }\Omega,\vspace{3pt}\\
u & = & 0 \ \hbox{ on }\partial\Omega,\eaa\right.
\ee
in a sense to be detailed later, we want to compare $u$ to a solution $v$ of a similar problem in the ball $\Omega^*$, namely
\be\label{eqv}\left\{\baa{rcl}
-\div(\hat{A}(x)\nabla v)+\hat{H}(x,v,\nabla v) & = & 0\ \hbox{ in }\Omega^*,\vspace{3pt}\\
v & = & 0 \ \hbox{ on }\partial\Omega^*.\eaa\right.
\ee
Our goal is to show that, under some natural assumptions on $A$ and $H$, any solution~$u$ of~(\ref{equ}) in $\Omega$ will be controlled from above by a radially symmetric solution $v$ of a similar problem~(\ref{eqv}) in $\Omega^*$. More precisely, we will show that the distribution function of $u$ is not larger than that of $v$. Moreover, the coefficients of the problem~(\ref{eqv}) solved by $v$ in the ball $\Omega^*$ will actually be radially symmetric, and the solution $v$ itself will also be radially symmetric.

Let us list the precise notations and the assumptions attached to the problem~(\ref{equ}) and made throughout the paper. We denote by ${\mathcal S}_n(\R)$ the set of $n\times n$ symmetric matrices with real entries. We always assume that $A:\Omega\rightarrow {\mathcal S}_n(\R)$ is in $W^{1,\infty}(\Omega)$. This assumption will be denoted by $A=(A_{i,i'})_{1\le i,i'\le n}\in W^{1,\infty}(\Omega,{\mathcal S}_n(\R))$: all the components $A_{i,i'}$ are in $W^{1,\infty}(\Omega)$ and they can therefore be assumed to be continuous in $\overline{\Omega}$ up to a modification on a zero-measure set. We always assume that $A$ is uniformly elliptic in $\overline{\Omega}$: there exists $\lambda>0$ such that $A\ge\lambda\,\hbox{Id}$ in $\overline{\Omega}$, where $\hbox{Id}\in\S_n(\R)$ is the identity matrix, that is
$$\sum_{1\le i,i'\le n}A_{i,i'}(x)\xi_i\xi_{i'}=A(x)\xi\cdot\xi\ge\lambda\left\vert \xi\right\vert^2:=\lambda\sum_{1\le i\le n}(\xi_i)^2$$
for all $x\in\overline{\Omega}$ and $\xi\in\R^n$. Actually, in some statements we compare the matrix field $A$ with a matrix field of the type $x\mapsto\Lambda(x)\hbox{Id}$ in the sense that
\be\label{ALambda}
A(x)\ge\Lambda(x)\hbox{Id}\hbox{ a.e. in }\Omega,
\ee
where $\Lambda\in L^{\infty}_+(\Omega)$ and
$$L^{\infty}_+(\Omega)=\big\{\Lambda\in L^{\infty}(\Omega);\ \mathop{\hbox{ess inf}}_{\Omega}\Lambda>0\big\}=\big\{\Lambda\in L^{\infty}(\Omega);\ \exists\, \lambda>0,\ \Lambda(x)\geq \lambda\mbox{ a.e. in }\Omega\big\}.$$
The inequality~(\ref{ALambda}) means that $A(x)\xi\cdot\xi\ge\Lambda(x)|\xi|^2$ for a.e. $x\in\Omega$ and all $\xi\in\R^n$. For instance, if, for each $x\in\overline{\Omega}$, $\Lambda_A(x)$ denotes the smallest eigenvalue of the matrix $A(x)$, then $\Lambda_A\in L^{\infty}_+(\Omega)$ and there holds $A(x)\ge\Lambda_A(x)\hbox{Id}$ for all $x\in\overline{\Omega}$. The given function $H:\overline{\Omega}\times\R\times\R^n\to\R$ is assumed to be continuous and such that there exist a real number
$$1\le q\le 2,$$
and three bounded and continuous functions $a$, $b$ and $f:\overline{\Omega}\times\R\times\R^n\to\R$ such that
\be\label{hypH}\left\{\baa{l}
H(x,s,p)\ge-a(x,s,p)|p|^q+b(x,s,p)s-f(x,s,p),\vspace{3pt}\\
b(x,s,p)\ge 0\eaa\right.
\ee
for all $(x,s,p)\in\overline{\Omega}\times\R\times\R^n$. In particular, $H$ is bounded from below by an at most quadratic function in its last variable $p$. Notice that no bound from above is assumed a priori. The cases~$q=1$ and $1<q\le 2$ will actually be treated separately, since the existence and uniqueness results for problems~(\ref{equ}) or~(\ref{eqv}) are different whether $q$ be equal to or larger than~$1$, and since an additional condition will be used when $q>1$. 

We say that $u$ is a weak solution of~(\ref{equ}) if $u\in H^1_0(\Omega)$, $H(\cdot,u(\cdot),\nabla u(\cdot))\in L^1(\Omega)$ and
$$\int_{\Omega}A(x)\nabla u\cdot\nabla\varphi+\int_{\Omega}H(x,u,\nabla u)\varphi=0\ \hbox{ for all }\varphi\in H^1_0(\Omega)\cap L^{\infty}(\Omega).$$
The condition~$u=0$ on $\partial\Omega$ simply means that the trace of $u$ on $\partial\Omega$ is equal to $0$. When $H(\cdot,u(\cdot),\nabla u(\cdot))$ belongs to $L^2(\Omega)$, then this identity holds for all test functions $\varphi$ in $H^1_0(\Omega)$. Similarly, one defines the notion of weak solution $v$ of \eqref{eqv}. Throughout the paper, the solutions of~(\ref{equ}) and \eqref{eqv} are always understood as weak solutions, even if they may of course be stronger under some additional assumptions on the coefficients. Furthermore, we denote $W(\Omega)$ the space
$$W(\Omega)=\bigcap_{1\le p<+\infty}W^{2,p}(\Omega).$$
We recall from the Sobolev embeddings that any function $u$ in $W(\Omega)$ belongs to $C^{1,\alpha}(\overline{\Omega})$ for all~$\alpha\in[0,1)$, even if it means redefining $u$ in a negligible subset of $\overline{\Omega}$. Notice that if $u\in W(\Omega)$ is a weak solution of~(\ref{equ}), then $u$ is a strong solution, the equation~(\ref{equ}) is satisfied almost everywhere in~$\Omega$ and the boundary condition on $\partial\Omega$ holds in the pointwise sense.

In order to compare a solution of~(\ref{equ}) in $\Omega$ to another function defined in the equimeasurable ball~$\Omega^*$, the natural way is to use their distribution functions. Namely, for any function~$u\in L^1(\Omega)$, let $\mu_u$ be the distribution function of $u$ given by
$$\mu_u(t)=\big|\big\{x\in \Omega;\ u(x)>t\big\}\big|$$
for all $t\in\R$. Note that $\mu_u$ is right-continuous, non-increasing, $\mu_u(t)\to|\Omega|$ as $t\to-\infty$ and~$\mu_u(t)\to 0$ as $t\to+\infty$. For all  $x\in\Omega^{\ast}\backslash\{0\}$, define
$$u^{\ast}(x)=\min\big\{t\in\R;\ \mu_u(t)\leq\alpha_n\left\vert x\right\vert^n\big\},$$
where $\alpha_n=\pi^{n/2}/\Gamma(n/2+1)$ denotes the Lebesgue measure of the Euclidean unit ball in~$\R^n$. The function $u^{\ast}$, called the decreasing Schwarz rearrangement of $u$, is clearly radially symmetric, non-increasing in the variable $|x|$ and it satisfies
$$\big|\big\{x\in\Omega;\ u(x)>\zeta\big\}\big|=\big|\big\{x\in\Omega^{\ast};\ u^{\ast}(x)>\zeta\big\}\big|$$
for all $\zeta\in \R$. An essential property of the Schwarz symmetrization is the following one: if $u$ belongs to the space $H^1_0(\Omega)$, then $|u|^{\ast}\in H^1_0(\Omega^{\ast})$ and it is such that $\left\Vert\,|u|^{\ast}\right\Vert_{L^2(\Omega^{\ast})}=\left\Vert u\right\Vert_{L^2(\Omega)}$ and $\left\Vert \nabla |u|^{\ast}\right\Vert_{L^2(\Omega^{\ast})}\leq \left\Vert \nabla u\right\Vert_{L^2(\Omega)}$, see~\cite{bramanti,burchard,k,polyaszego}.

Our goal in this paper is to compare any given positive solution $u\in W(\Omega)$ of~(\ref{equ}) to a weak solution $v$ of a problem of the type~(\ref{eqv}), in the sense that $u^*\le v$ in $\Omega^*$, where the new coefficients $\hat{A}$ and $\hat{H}$ of~(\ref{eqv}) are radially symmetric with respect to $x\in\Omega^*$ and satisfy similar constraints or bounds as the given coefficients $A$ and $H$. As already mentioned, we consider two types of assumptions regarding the dependency of the lower bound~(\ref{hypH}) with respect to the variable~Ê$p$: namely, we treat separately the case where the lower bound is at most linear in $|p|$ (that is,~$q=1$) and the general case where $1<q\le 2$ and the lower bound is at most quadratic in~$|p|$, for which an additional assumption on the function~$b$ is made.


\subsection{Linear growth with respect to the gradient}\label{sec11}

Our first main result is concerned with the case where $H$ is linear in $|\nabla u|$ from below, in the sense that $q=1$ in~(\ref{hypH}). If~$g$ is a real number or a real-valued function, we set
$$g^+=\max(g,0).$$

\begin{theo}\label{th1}
Assume~$(\ref{ALambda})$ and~$(\ref{hypH})$ with $\Lambda\in L^{\infty}_+(\Omega)$ and $q=1$. Let $u\in W(\Omega)$ be a solution of~$(\ref{equ})$ such that $u>0$ in $\Omega$ and $|\nabla u|\neq 0$ everywhere on $\partial\Omega$. Then there are two radially symmetric functions $\hat{\Lambda}\in L^{\infty}_+(\Omega^*)$ and $\hat{a}\in L^{\infty}(\Omega^*)$ such that
\be\label{hats}\left\{\baa{l}
0<\displaystyle{\mathop{\rm{ess}\,\rm{inf}}_{\Omega}}\,\Lambda\le\displaystyle{\mathop{\rm{ess}\,\rm{inf}}_{\Omega^*}}\,\hat{\Lambda}\le\displaystyle{\mathop{\rm{ess}\,\rm{sup}}_{\Omega^*}}\,\hat{\Lambda}\le\displaystyle{\mathop{\rm{ess}\,\rm{sup}}_{\Omega}}\,\Lambda,\ \ \|\hat{\Lambda}^{-1}\|_{L^1(\Omega^*)}=\|\Lambda^{-1}\|_{L^1(\Omega)},\vspace{3pt}\\
0\le\displaystyle{\mathop{\inf}_{\Omega\times\R\times\R^n}}\,a^+\le\displaystyle{\mathop{\rm{ess}\,\rm{inf}}_{\Omega^*}}\,\hat{a}\le\displaystyle{\mathop{\rm{ess}\,\rm{sup}}_{\Omega^*}}\,\hat{a}\le\displaystyle{\mathop{\sup}_{\Omega\times\R\times\R^n}}\,a^+,\eaa\right.
\ee
and
\be\label{u*v}
u^*\le v\hbox{ a.e. in }\Omega^*,
\ee
where $v\in H^1_0(\Omega^*)\cap C(\overline{\Omega^*})$ is the unique weak solution of
\be\label{eqvbis}\left\{\baa{rcl}
-{\rm{div}}(\hat{\Lambda}(x)\nabla v)+\hat{H}(x,\nabla v) & = & 0\ \hbox{ in }\Omega^*,\vspace{3pt}\\
v & = & 0 \ \hbox{ on }\partial\Omega^*\eaa\right.
\ee
with
$$\hat{H}(x,p)=-\hat{a}(x)\,|p|-f_u^*(x)$$
and $f_u^*$ is the Schwarz rearrangement of the function $f_u$ defined in $\Omega$ by
\be\label{deffu}
f_u(y)=f(y,u(y),\nabla u(y))\hbox{ for all }y\in\Omega.
\ee
\end{theo}

Notice that, since $\hat{H}(\cdot,\nabla v(\cdot))\in L^2(\Omega^{\ast})$, the fact that $v\in H^1_0(\Omega^*)$ is a weak solution of~(\ref{eqvbis}) means that
$$\int_{\Omega^*}\,\hat{\Lambda}(x)\nabla v\cdot\nabla\varphi-\int_{\Omega^*}\hat{a}(x)\,|\nabla v|\,\varphi-\int_{\Omega^*}f^*_u(x)\,\varphi=0$$
for all $\varphi\in H^1_0(\Omega^*)$.

From Theorem~\ref{th1} and the maximum principle, the following corollary holds.

\begin{cor}\label{cor1}
Under the notations of Theorem~$\ref{th1}$, for any functions $\overline{a}$ and~$\overline{f}$ in~$L^{\infty}(\Omega^*)$ such that~$\hat{a}\le\overline{a}$ and~$f_u^*\le\overline{f}$ a.e. in $\Omega^*$, there holds
$$u^*\le\overline{v}\hbox{ a.e. in }\Omega^*,$$
where $\overline{v}\in H^1_0(\Omega^*)\cap C(\overline{\Omega^*})$ is the unique weak solution of
\be\label{eqvter}\left\{\baa{rcl}
-{\rm{div}}(\hat{\Lambda}(x)\nabla\overline{v})+\overline{H}(x,\nabla\overline{v}) & = & 0\ \hbox{ in }\Omega^*,\vspace{3pt}\\
\overline{v} & = & 0 \ \hbox{ on }\partial\Omega^*\eaa\right.
\ee
with
$$\overline{H}(x,p)=-\overline{a}(x)\,|p|-\overline{f}(x).$$
In particular, there holds $u^*\le\overline{V}$ a.e. in $\Omega^*$, where $\overline{V}\in H^1_0(\Omega^*)\cap C(\overline{\Omega^*})$ is the unique weak solution of~$(\ref{eqvter})$ with $\overline{a}=\sup_{\Omega\times\R\times\R^n}a^+$ and $\overline{f}=\sup_{\Omega\times\R\times\R^n}f$.
\end{cor}

\noindent{\bf{Proof.}} The proof is an immediate consequence of Corollary~2.1 of Porretta~\cite{Porretta}: indeed, with the notations of Theorem~\ref{th1} and Corollary~\ref{cor1}, there holds
$$-\hbox{div}(\hat{\Lambda}(x)\nabla v)+\overline{H}(x,\nabla v)=(\hat{a}(x)-\overline{a}(x))|\nabla v|+f^*_u(x)-\overline{f}(x)\le0$$
in the weak $H^1_0(\Omega^*)$ sense, that is
$$\int_{\Omega^*}\,\hat{\Lambda}(x)\nabla v\cdot\nabla\varphi-\int_{\Omega^*}\overline{a}(x)\,|\nabla v|\,\varphi-\int_{\Omega^*}\overline{f}(x)\,\varphi\le0$$
for all $\varphi\in H^1_0(\Omega^*)$ with $\varphi\ge 0$ a.e. in $\Omega^*$. In other words, $v\in H^1_0(\Omega^*)$ is a weak subsolution of~(\ref{eqvter}). The maximum principle (Corollary~2.1 of~\cite{Porretta}) yields $v\le\overline{v}$ a.e. in $\Omega^*$, whence $u^*\le\overline{v}$ a.e. in $\Omega^*$ from~(\ref{u*v}).\hfill\fin\break

For problems~(\ref{eqvbis}) and~(\ref{eqvter}) in~$\Omega^*$, the existence and uniqueness of the solutions $v$ and~$\overline{v}$ in $H^1_0(\Omega^*)$ is a direct consequence of Theorem 2.1 of Porretta~\cite{Porretta}. For~(\ref{eqvbis}), since all coefficients~$\hat{\Lambda}$,~$\hat{a}$ and~$f_u^*$ are radially symmetric, it follows from the uniqueness that $v$ is itself radially symmetric. Furthermore, $v$ is then H\"older continuous in $\overline{\Omega^*}$ from the radial symmetry and the local De Giorgi-Moser-Nash estimates (see Theorem~8.29 in~\cite{gt}). For problem~(\ref{eqvter}), the function~$\overline{v}$ is still locally H\"older continuous in $\Omega^{\ast}$. It may not be radially symmetric in general, since the functions $\overline{a}$ and $\overline{f}$ are not assumed to be radially symmetric. But it follows from the maximum principle and Corollary \ref{cor1} that $0\leq u^{\ast}\leq \overline{v}\leq \overline{V}$ in $\Omega^{\ast}$, and since $\overline{V}$ is continuous in $\overline{\Omega^{\ast}}$ and vanishes on $\partial\Omega^{\ast}$, the function $\overline{v}$ is continuous in~$\overline{\Omega^*}$ too.

As far as problem~(\ref{equ}) in~$\Omega$ is concerned, additional conditions guaranteeing the existence and uniqueness of a solution $u$ can be given. Namely, if one assumes that
\be\label{hypH2}\left\{\baa{l}
\exists\,\omega>0,\ \forall\,(x,s,s^{\prime},p)\in \Omega\times \R\times \R\times \R^n,\\
\qquad\qquad\qquad\qquad\qquad\qquad\omega^{-1}(s-s^{\prime})\leq H(x,s,p)-H(x,s^{\prime},p)\leq\omega(s-s^{\prime}),\vspace{3pt}\\
\exists\,\alpha\in L^{r}(\Omega),\ \forall\,(x,s,p,p')\in\Omega\times\R\times\R^n\times\R^n,\vspace{3pt}\\
\qquad\qquad\qquad\qquad\qquad\qquad|H(x,s,p)-H(x,s,p')|\le\alpha(x)(1+|s|^{2/n})|p-p'|,\eaa\right.
\ee
where $r=n^2/2$ if $n\geq 3$, $r\in (2,+\infty)$ if $n=2$ and $r=2$ if $n=1$, then there is at most one solution of~(\ref{equ}) in $H^1_0(\Omega)$ (see Theorem~\ref{thapp} below). If, in addition to~(\ref{hypH2}), one assumes that
\be\label{hypH3}
\exists\,\beta\in L^t(\Omega),\ \forall\,(x,s,p)\in\Omega\times\R\times\R^n,\ |H(x,s,p)|\le\beta(x)(1+|s|+|p|),
\ee
where $t=n$ if $n\ge 3$, $t$ is any real number in $(2,+\infty)$ if $n=2$ and $t=2$ if $n=1$, then there exists a (unique) solution $u$ of~(\ref{equ}) in $H^1_0(\Omega)$. Furthermore, if the function~$\beta$ in~(\ref{hypH3}) is such that
\be\label{hypH4}
\beta\in L^{\infty}(\Omega),
\ee
then the unique solution $u$ of~(\ref{equ}) belongs to the space $W(\Omega)$. Lastly, under the additional assumption
\be\label{hypH5}\left\{\baa{l}
\exists\,\varpi\ge0,\ \forall\,(x,p)\in\Omega\times\R^n,\ H(x,0,p)\le\varpi\,|p|,\vspace{3pt}\\
H(\cdot,0,0)\not\equiv 0\hbox{ in }\Omega,\eaa\right.
\ee
then $u>0$ in $\Omega$ and $\partial_{\nu}u:=\nu\cdot\nabla u<0$ on $\partial\Omega$, where $\nu$ denotes the outward unit normal to~$\partial\Omega$. These aforementioned existence and uniqueness results, which are inspired from~\cite{Porretta}, are summarized in Theorem~\ref{thapp} in the appendix (Section~\ref{secapp}). For further uniqueness results on semilinear problems of the type~(\ref{equ}) with linear growth in $|\nabla u|$, we refer to~\cite{bgm92,dg01}.

Let us now compare Theorem~\ref{th1} with some existing results in the literature. Theorem~\ref{th1} provides a comparison of the distribution functions of a given solution of~(\ref{equ}) in $\Omega$ and a solution $v$ of~(\ref{eqvbis}) in $\Omega^*$. The coefficients $\hat{\Lambda}$ and $\hat{H}$ in~(\ref{eqvbis}) satisfy the same type of pointwise bounds as the coefficients $A$ and $H$ in $\Omega$. Furthermore, the diffusion matrix $\hat{\Lambda}\,\hbox{Id}$ in the second-order term of~(\ref{eqvbis}) is proportional to the identity matrix at each point $x$, and the nonlinear term~$\hat{H}(x,\nabla v)$ is exactly affine in $|\nabla v|$. The first comparison result in the spirit of Theorem~\ref{th1} goes back to the seminal paper of Talenti~\cite{Talenti76}: if $A\in L^{\infty}(\Omega,\mathcal{S}_n(\R))$ and
$$A\ge\hbox{Id}\hbox{ a.e. in }\Omega,\ \ \Lambda=1\ \hbox{ and }\ H(x,s,p)=b(x)s-f(x)$$
with $b,\,f\in L^{\infty}(\Omega)$ and $b\ge 0$ a.e. in $\Omega$, then equation~(\ref{equ}) admits a unique solution $u\in H^1_0(\Omega)$ and $\left\vert u\right\vert^*\le v$ a.e. in $\Omega^*$, where $v\in W(\Omega^*)$ is the unique solution of $-\Delta v=\left\vert f\right\vert^*$ in~$\Omega^*$ with $v=0$ on $\partial\Omega^*$. If $A\in L^{\infty}(\Omega,\mathcal{S}_n(\R))$ and
\be\label{hyptalenti}\left\{\baa{l}
A\ge\hbox{Id}\hbox{ a.e. in }\Omega,\ \ \Lambda=1,\ H(x,s,p)=\alpha(x)\cdot p+b(x)s-f(x),\vspace{3pt}\\
\alpha\in L^{\infty}(\Omega,\R^n),\ b,\,f\in L^{\infty}(\Omega),\ b\ge 0\hbox{ a.e. in }\Omega,\eaa\right.
\ee
then it follows from Talenti~\cite{talenti} that the unique solution $u\in H^1_0(\Omega)$ of~(\ref{equ}) satisfies
\be\label{ineqtalenti}
\left\vert u\right\vert^*\le v\ \hbox{ a.e. in }\Omega^*,
\ee
where $v\in W(\Omega^*)\cap H^1_0(\Omega^*)$ is the unique solution of
\be\label{eqtalenti}
-\Delta v+\tilde{\alpha}\,e_r\cdot\nabla v=\left\vert f\right\vert^*\hbox{ in }\Omega^*
\ee
with $\tilde{\alpha}=\|\,|\alpha|\,\|_{L^{\infty}(\Omega)}$ and $e_r(x)=x/|x|$ for all $x\in\Omega^*\backslash\{0\}$. Actually, assuming that the matrix field $A$ is continuous in $\overline{\Omega}$, these seminal results of~\cite{Talenti76} and~\cite{talenti} can be recovered from Theorem~\ref{th1} of the present work, as will be explained in Section \ref{backtalenti} below.  
Furthermore, Trombetti and Vazquez~\cite{tv} proved that if, in particular,  $A\ge\hbox{Id}$, $\Lambda=1$, $H(x,s,p)=\alpha(x)\cdot p+bs-f(x)$ with~$\alpha\in W^{1,\infty}(\Omega,\R^n)$, $b,\,f\in L^{\infty}(\Omega)$, $\min(b,\hbox{div}(\alpha)+b)\ge0$ a.e. in $\Omega$, then~$\|u\|_{L^p(\Omega)}\le\|v\|_{L^p(\Omega^*)}$ for all~$1\le p\le+\infty$, where~$v\in W(\Omega^*)$ solves~$-\Delta v+\tilde{\alpha}\,e_r\cdot\nabla v+c_*v=\left\vert f\right\vert^*$ in $\Omega^*$ with $v=0$ on~$\partial\Omega^*$, $\tilde{\alpha}=\|\,|\alpha|\,\|_{L^{\infty}(\Omega)}$ and $c_*\in L^{\infty}(\Omega^*)$ is the nondecreasing symmetric rearrangement of any function~$c$ such that $0\le c\le\min(b,\hbox{div}(\alpha)+b)$ in~$\Omega$. Further results in this spirit can be found in~\cite{atl,atlm,inrz}.

In the references cited in the paragraph above, the function $\Lambda$ appearing in~(\ref{ALambda}) was assumed to be constant. When $\Lambda$ is given as a constant $\lambda>0$ in~(\ref{ALambda}) (this is always possible since the matrix field $A$ is assumed to be uniformly elliptic), then it follows from~(\ref{hats}) that the function~$\hat{\Lambda}$ appearing in Theorem~\ref{th1} is equal to~$\lambda$ and the principal part in~(\ref{eqvbis}) is proportional to the Laplacian. However, in the present paper, $\Lambda$ is not assumed to be constant and the function $\hat{\Lambda}$ is actually not constant in general (see Remark~5.5 of~\cite{hnrAM} for some examples). Notice in particular from~(\ref{hats}) that
$$\mathop{{\rm{ess}}\,{\rm{inf}}}_{\Omega^*}\hat{\Lambda}\ge\mathop{{\rm{ess}}\,{\rm{inf}}}_{\Omega}\Lambda\ \hbox{ and }\ \mathop{{\rm{ess}}\,{\rm{sup}}}_{\Omega^*}\hat{\Lambda}\le\mathop{{\rm{ess}}\,{\rm{sup}}}_{\Omega}\Lambda,$$
and that
$$\mathop{{\rm{ess}}\,{\rm{sup}}}_{\Omega^*}\hat{\Lambda}>\mathop{{\rm{ess}}\,{\rm{inf}}}_{\Omega}\Lambda\ \hbox{ and }\ \mathop{{\rm{ess}}\,{\rm{inf}}}_{\Omega^*}\hat{\Lambda}<\mathop{{\rm{ess}}\,{\rm{sup}}}_{\Omega}\Lambda$$
as soon as $\Lambda$ is not constant. To our best knowledge, the first occurrence of a non-constant function~$\Lambda$ in comparison results for elliptic problems of the type~(\ref{equ}) goes back to a paper by Alvino and Trombetti~\cite{at3}. In this paper, the authors consider the principal eigenvalue~$\lambda_1(\Omega,A)$ of the operator $-\hbox{div}(A(x)\nabla)$ in~$\Omega$ with Dirichlet boundary conditions and they prove that~$\lambda_1(\Omega,A)\ge\lambda_1(\Omega^*,\hat{\Lambda}\,\hbox{Id})$ for some function $\hat{\Lambda}$ which shares the same pro\-perties as that  given in Theorem~\ref{th1}. An analogous method is used in~\cite{at4} for the problem~$-\hbox{div}(A(x)\nabla u)=1$. We also mention~\cite{cms} for some results on the optimization of~$\lambda_1(\Omega,\Lambda\,{\rm{Id}})$ in the class of binary piecewise constant functions~$\Lambda$.

The problems described in the previous two paragraphs were concerned with functions~$H(x,u,\nabla u)$ which were linear with respect to $u$ and $\nabla u$. Analogous results for nonlinear problems have been established in~\cite{atl,bm,cianchi,fm,messano,talenti79}. In most of these works, the second-order coefficients for problems of the type~(\ref{equ}) in $\Omega$ are compared with a multiple of the Laplace operator for a problem of the type~(\ref{eqv}) in $\Omega^*$ (or with homogeneous second-order terms such as the $p$-Laplacian), and the comparisons between $u^*$ and $v$ are either pointwise or only integral and hold either in the whole ball $\Omega^*$ or only in a strict subdomain, depending on the assumptions on the lower-order coefficients. In Theorem~\ref{th1}, the original problem~(\ref{equ}) is nonlinear with respect to $u$ and $\nabla u$, the highest-order terms in~(\ref{equ}) are compared with equations having heterogeneous second-order terms, and the comparison between $u^*$ and $v$ are pointwise in the whole ball $\Omega^*$. One of the main novelties is also the method, which involves a symmetrization of the second-order terms with respect to the level sets of $u$. We refer to the proofs in the following sections for more details. 

As a matter of fact, the method used in the proof of Theorem~\ref{th1} leads to a quantified comparison result in the case where the original domain $\Omega$ is not a ball. We recall that $A=(A_{i,i'})_{1\le i,i'\le n}$ and we denote $\|A\|_{W^{1,\infty}(\Omega)}=\max_{1\le i,i'\le n}\|A_{i,i'}\|_{W^{1,\infty}(\Omega)}$.

\begin{theo}\label{th1bis}
Assume that $\Omega$ is not a ball, that~$(\ref{ALambda})$ and~$(\ref{hypH})$ hold with $\Lambda\in L^{\infty}_+(\Omega)$ and $q=1$. Then, under the notations of Theorem~$\ref{th1}$, there is a constant $\eta_u>0$, which depends on $\Omega$, $n$ and $u$, such that
\be\label{u*vetau}
(1+\eta_u)\,u^*\le v\hbox{ a.e. in }\Omega^*.
\ee
Furthermore, if there is $M>0$ such that
\be\label{hypAM}\left\{\baa{l}
\|A\|_{W^{1,\infty}(\Omega)}+\|\Lambda^{-1}\|_{L^{\infty}(\Omega)}+\|a\|_{L^{\infty}(\Omega\times\R\times\R^n)}+\|f\|_{L^{\infty}(\Omega\times\R\times\R^n)}\le M,\vspace{3pt}\\
|H(x,s,p)-H(x,0,0)|\le M\,(|s|+|p|)\hbox{ for all }(x,s,p)\in\Omega\times\R\times\R^n,\vspace{3pt}\\
-M\le H(x,0,0)\le 0\hbox{ for all }x\in\Omega,\ \ \displaystyle\int_{\Omega}H(x,0,0)\,dx\le-M^{-1}<0,\eaa\right.
\ee
then there is a constant $\eta>0$, which depends on $\Omega$, $n$ and $M$ but not on $u$, such that
\be\label{u*veta}
(1+\eta)\,u^*\le v\hbox{ a.e. in }\Omega^*.
\ee
\end{theo}

In addition to the aforementioned differences with respect to the existing results in the literature, the improved quantified compa\-risons stated in Theorem~\ref{th1bis} for problems of the type~(\ref{equ}) and~(\ref{eqvbis}) when $\Omega$ is not a ball have, to our knowledge, never been established in the literature, even in particular situations. 

\begin{rem}{\rm In our results, the solution $u$ of~(\ref{equ}) is assumed to be in $W(\Omega)$ (notice that this assumption is equivalent to $u\in C^1(\overline{\Omega})$, from the standard elliptic estimates applied to~\eqref{equ} and the fact that the function~$H$ is continuous). In the proof of our results, $u$ is approximated by some analytic functions $u_j$ in~$W^{2,p}(\Omega)$ weakly for all $1\le p<+\infty$ and in~$C^{1,\alpha}(\overline{\Omega})$ strongly for all $0\le\alpha<1$. Since $u>0$ in $\Omega$ and $|\nabla u|>0$ on $\partial\Omega$ by assumption, the functions $u_j$ satisfy these properties for large $j$ and one can then apply to them the symmetrization method described in Section~\ref{rearrangement}. In the approximation of $u$ by $u_j$, the $W^{2,p}$ theory is needed and the Lipschitz-continuity of the matrix field $A$ is used. It is actually beyond the scope of this paper to see under which minimal regularity assumptions on the coefficients of~\eqref{equ} the main results would still hold. Whereas the proof of Talenti's results are based on the Schwarz symmetrization, on P\'olya-Szeg\"o inequality and on the standard geometric isoperimetric inequality, our proofs are based on the symmetrization of the second-order terms and they require more regularity assumptions on the equation. However, our results also cover the case where the ellipticity functions $\Lambda$ in $\Omega$ and~$\hat{\Lambda}$ in $\Omega^*$ are not constant in general. Furthermore, they provide inequalities which can be quantitatively expressed in terms of some bounds on the coefficients and the domain and which can then be quantitatively improved when the domain is not a ball.}
\end{rem}


\subsection{At most quadratic growth with respect to the gradient}

Our second main result is concerned with the general case where $H$ is at most quadratic in~$|\nabla u|$ from below, under the additional assumption that $\inf_{\Omega\times\R\times\R^n}b>0$ in~(\ref{hypH}).

\begin{theo}\label{th2}
Assume that~$(\ref{ALambda})$ and~$(\ref{hypH})$ with $\Lambda\in L^{\infty}_+(\Omega)$, $1<q\le 2$ and
\be\label{hypb}
\inf_{\Omega\times\R\times\R^n}b>0.
\ee
Let $u\in W(\Omega)$ be a solution of~$(\ref{equ})$ such that $u>0$ in $\Omega$ and~$|\nabla u|\neq 0$ everywhere on $\partial\Omega$. Then there are three radially symmetric functions~$\hat{\Lambda}\in L^{\infty}_+(\Omega^*)$, $\hat{a}\in L^{\infty}(\Omega^*)$ and $\hat{f}\in L^{\infty}(\Omega^*)$, and a positive constant $\hat{\delta}$ such that
\be\label{hats4}\left\{\baa{l}
0<\displaystyle{\mathop{\rm{ess}\,\rm{inf}}_{\Omega}}\,\Lambda\le\displaystyle{\mathop{\rm{ess}\,\rm{inf}}_{\Omega^*}}\,\hat{\Lambda}\le\displaystyle{\mathop{\rm{ess}\,\rm{sup}}_{\Omega^*}}\,\hat{\Lambda}\le\displaystyle{\mathop{\rm{ess}\,\rm{sup}}_{\Omega}}\,\Lambda,\ \ \|\hat{\Lambda}^{-1}\|_{L^1(\Omega^*)}=\|\Lambda^{-1}\|_{L^1(\Omega)},\vspace{3pt}\\
0\le\Big(\displaystyle{\mathop{\inf}_{\Omega\times\R\times\R^n}}\,a^+\Big)\times\Big(\frac{\rm{ess}\,\rm{inf}_{\Omega}\Lambda}{\rm{ess}\,\rm{sup}_{\Omega}\Lambda}\Big)^{q-1}\vspace{3pt}\\
\qquad\qquad\qquad\qquad\le\displaystyle{\mathop{\rm{ess}\,\rm{inf}}_{\Omega^*}}\,\hat{a}\le\displaystyle{\mathop{\rm{ess}\,\rm{sup}}_{\Omega^*}}\,\hat{a}\le\Big(\displaystyle{\mathop{\sup}_{\Omega\times\R\times\R^n}}\,a^+\Big)\times\Big(\frac{\rm{ess}\,\rm{sup}_{\Omega}\Lambda}{\rm{ess}\,\rm{inf}_{\Omega}\Lambda}\Big)^{2q-2},\vspace{3pt}\\
\displaystyle{\mathop{\inf}_{\Omega\times\R\times\R^n}}\,f\le\displaystyle{\mathop{\rm{ess}\,\rm{inf}}_{\Omega^*}}\,\hat{f}\le\displaystyle{\mathop{\rm{ess}\,\rm{sup}}_{\Omega^*}}\,\hat{f}\le\displaystyle{\mathop{\sup}_{\Omega\times\R\times\R^n}}\,f,\ \ \displaystyle\int_{\Omega^*}\hat{f}=\int_{\Omega}f_u,\eaa\right.
\ee
and
\be\label{u*vbis}
u^*\le v\hbox{ a.e. in }\Omega^*,
\ee
where $f_u$ is defined as in~$(\ref{deffu})$ and $v\in H^1_0(\Omega^*)\cap L^{\infty}(\Omega^*)$ is the unique weak solution of
\be\label{eqv4}\left\{\baa{rcl}
-{\rm{div}}(\hat{\Lambda}(x)\nabla v)+\hat{H}(x,v,\nabla v) & = & 0\ \hbox{ in }\Omega^*,\vspace{3pt}\\
v & = & 0 \ \hbox{ on }\partial\Omega^*\eaa\right.
\ee
with
\be\label{defhatH}
\hat{H}(x,s,p)=-\hat{a}(x)\,|p|^q+\hat{\delta}\,s-\hat{f}(x).
\ee
Furthermore, for every $\epsilon>0$, there is a radially symmetric function $\hat{f}_{\epsilon}\in L^{\infty}(\Omega^*)$ such that~$\mu_{\hat{f}_{\epsilon}}=\mu_{f_u}$ and
\be\label{u*vepsilon}
\|(u^*-v_{\epsilon})^+\|_{L^{2^*}(\Omega^*)}\le\epsilon,
\ee
where $v_{\epsilon}\in H^1_0(\Omega^*)\cap L^{\infty}(\Omega^*)$ is the unique weak solution of
\be\label{eqv4eps}\left\{\baa{rcl}
-{\rm{div}}(\hat{\Lambda}(x)\nabla v_{\epsilon})+\hat{H}_{\epsilon}(x,v_{\epsilon},\nabla v_{\epsilon}) & = & 0\ \hbox{ in }\Omega^*,\vspace{3pt}\\
v_{\epsilon} & = & 0 \ \hbox{ on }\partial\Omega^*\eaa\right.
\ee
with
\be\label{defhatHeps}
\hat{H}_{\epsilon}(x,s,p)=-\hat{a}(x)\,|p|^q+\hat{\delta}\,s-\hat{f}_{\epsilon}(x).
\ee
In~$(\ref{u*vepsilon})$, $2^*=2n/(n-2)$ if $n\ge 3$, $2^*=\infty$ if $n=1$ and $2^*$ is any fixed real number in~$[1,+\infty)$ if $n=2$.
\end{theo}

Notice that, contrary to the conclusion of Theorem \ref{th1}, $\hat{f}$ and $f_u$ do not have the same distribution function in general, but $\hat{f}_{\varepsilon}$ and $f_u$ do.

Since $1<q\le 2$, the functions $\hat{a}\,|\nabla v|^q$ and $\hat{a}\,|\nabla v_{\epsilon}|^q$ are only in $L^1(\Omega^*)$ in general. Since~$\hat{H}(\cdot,v(\cdot),\nabla v(\cdot))\in L^1(\Omega^{\ast})$, the fact that $v\in H^1_0(\Omega^*)\cap L^{\infty}(\Omega^*)$ is a weak solution of~(\ref{eqv4}) means that
$$\int_{\Omega^*}\,\hat{\Lambda}(x)\nabla v\cdot\nabla\varphi-\int_{\Omega^*}\hat{a}(x)\,|\nabla v|^q\,\varphi+\int_{\Omega^*}\hat{\delta}\,v\,\varphi-\int_{\Omega^*}\hat{f}(x)\,\varphi=0$$
for all $\varphi\in H^1_0(\Omega^*)\cap L^{\infty}(\Omega^*)$, and similarly for $v_{\epsilon}$ with $\hat{f}_{\epsilon}$ instead of $\hat{f}$.

From Theorem~\ref{th2} and the maximum principle, the following corollary will be infered.

\begin{cor}\label{cor2}
Under the notations of Theorem~$\ref{th2}$, for any functions $\overline{a}$ and~$\overline{f}$ in~$L^{\infty}(\Omega^*)$ such that~$\hat{a}\le\overline{a}$ and~$\hat{f}\le\overline{f}$ a.e. in~$\Omega^*$, and for any constant~$\overline{\delta}$ such that $0<\overline{\delta}\le\hat{\delta}$, there holds
$$u^*\le\overline{v}\hbox{ a.e. in }\Omega^*,$$
where $\overline{v}\in H^1_0(\Omega^*)\cap L^{\infty}(\Omega^*)$ is the unique weak solution of
\be\label{eqv5}\left\{\baa{rcl}
-{\rm{div}}(\hat{\Lambda}(x)\nabla\overline{v})+\overline{H}(x,\overline{v},\nabla\overline{v}) & = & 0\ \hbox{ in }\Omega^*,\vspace{3pt}\\
\overline{v} & = & 0 \ \hbox{ on }\partial\Omega^*\eaa\right.
\ee
with
$$\overline{H}(x,s,p)=-\overline{a}(x)\,|p|^q+\overline{\delta}\,s-\overline{f}(x).$$
In particular, there holds $u^*\le\overline{V}$ a.e. in $\Omega^*$, where $\overline{V}\in H^1_0(\Omega^*)\cap L^{\infty}(\Omega^*)$ is the unique weak solution of~$(\ref{eqv5})$ with $\overline{\delta}=\hat{\delta}$, $\overline{a}=\big(\sup_{\Omega\times\R\times\R^n}a^+\big)\times\big({\rm{ess}}\,{\rm{sup}}_{\Omega}\Lambda/{\rm{ess}}\,{\rm{inf}}_{\Omega}\Lambda\big)^{2q-2}$ and $\overline{f}=\sup_{\Omega\times\R\times\R^n}f$.
\end{cor}

\noindent{\bf{Proof.}} With the notations of Theorem~\ref{th2} and Corollary~\ref{cor2}, together with the fact that~$v\ge u^*\ge0$ a.e. in $\Omega^*$ (from Theorem~\ref{th2}), there holds
$$-\hbox{div}(\hat{\Lambda}(x)\nabla v)+\overline{H}(x,v,\nabla v)=(\hat{a}(x)-\overline{a}(x))|\nabla v|^q+(\overline{\delta}-\hat{\delta})\,v+\hat{f}(x)-\overline{f}(x)\le0$$
in the weak $H^1_0(\Omega^*)\cap L^{\infty}(\Omega^*)$ sense, that is
$$\int_{\Omega^*}\,\hat{\Lambda}(x)\nabla v\cdot\nabla\varphi-\int_{\Omega^*}\overline{a}(x)\,|\nabla v|^q\,\varphi+\int_{\Omega^*}\overline{\delta}\,v\,\varphi-\int_{\Omega^*}\overline{f}(x)\,\varphi\le0$$
for all $\varphi\in H^1_0(\Omega^*)\cap L^{\infty}(\Omega^*)$ with $\varphi\ge0$ a.e. in $\Omega^*$. In other words, $v$ is a weak $H^1_0(\Omega^*)\cap L^{\infty}(\Omega^*)$ subsolution of~(\ref{eqv5}), with $\overline{\delta}>0$. The maximum principle (Theorem~2.1 of~\cite{bm95}) yields $v\le\overline{v}$ a.e. in $\Omega^*$, whence $u^*\le\overline{v}$ a.e. in $\Omega^*$ from~(\ref{u*vbis}).\hfill\fin\break

For problems~(\ref{eqv4}),~(\ref{eqv4eps}) and~(\ref{eqv5}), the existence of weak solutions $v$, $v_{\epsilon}$ and $\overline{v}$ in~$H^1_0(\Omega^*)\cap L^{\infty}(\Omega^*)$ follows from a paper by Boccardo, Murat and Puel (Th\'eor\`eme~2.1 and the following comments of~\cite{bmp84}, see also~\cite{bmp88}) and the uniqueness follows from Barles and Murat (Theorem~2.1 of~\cite{bm95}). Furthermore, even if it means redefining them in a negligible subset of~$\overline{\Omega^*}$,~$v$, $v_{\epsilon}$ and~$\overline{v}$ are locally H\"older continuous in $\Omega^*$ from Corollary~4.23 of~\cite{hl}, and~$v$ and $v_{\epsilon}$ are radially symmetric ($v$ and $v_{\epsilon}$ are then H\"older continuous in~$\overline{\Omega^*}$) by uniqueness since all coefficients of~(\ref{eqv4}) are radially symmetric. For problem~(\ref{eqv5}), the function~$\overline{v}$ may not be radially symmetric in general, since the functions $\overline{a}$ and $\overline{f}$ are not assumed to be radially symmetric. But it follows from the maximum principle and Corollary \ref{cor2} that~$0\leq u^{\ast}\leq \overline{v}\leq \overline{V}$ in $\Omega^{\ast}$, and since $\overline{V}$ is continous in $\overline{\Omega^{\ast}}$ and vanishes on $\partial\Omega^{\ast}$, the function $\overline{v}$ is continuous in~$\overline{\Omega^*}$ too.

More generally speaking, for problem~(\ref{equ}), it follows from~\cite{bmp88,bmp92} that there exists a solution $u\in H^1_0(\Omega)\cap L^{\infty}(\Omega)$ if
$$\left\{\baa{l}
H(x,s,p)=\beta(x,s,p)+\tilde{H}(x,s,p),\vspace{3pt}\\
\beta(x,s,p)s\ge\delta\,s^2,\ \ |\beta(x,s,p)|\le\kappa\,(\gamma(x)+|s|+|p|),\ \ |\tilde{H}(x,s,p)|\le\rho+\varrho(|s|)\,|p|^2\eaa\right.$$
with $\delta>0$, $\kappa>0$, $\rho>0$, $\gamma\in L^2(\Omega)$, $\gamma\ge 0$ a.e. in $\Omega$ and $\varrho:\R_+\to\R_+$ is increasing (see also~\cite{bbm88,dv87,Porretta} for further existence results). Furthermore, if $q<1+2/n$ and if there exist~$M\ge 0$ and $r\ge0$ such that $r(n-2)<n+2$ and
\be\label{hypHrq}
|H(x,s,p)|\le M(1+|s|^r+|p|^q)\hbox{ for all }(x,s,p)\in\Omega\times\R\times\R^n,
\ee
then any weak solution $u$ of~(\ref{equ}) belongs to~$W(\Omega)$, see Theorem~\ref{thapp2} below in the appendix. On the other hand, it follows from~\cite{bm95} that, if $\partial_sH$ and $\partial_pH$ exist with $\partial_sH(x,s,p)\ge\sigma>0$, $|\partial_pH(x,s,p)|\le\theta(|s|)\,(1+|p|)$, $|H(x,s,0)|\le\vartheta(|s|)$ for some continuous nonnegative functions~$\theta$ and~$\vartheta$, then there is at most one solution~$u$ of~(\ref{equ}) in $H^1_0(\Omega)\cap L^{\infty}(\Omega)$, and~$u$ is necessarily nonnegative if $H(\cdot,0,0)\le0$ in $\Omega$ (see also~\cite{bbgk99,Porretta} for other results in this direction). We refer to~\cite{bp06,bmmp05,gmp06} for further existence and uniqueness results for problems with strictly sub-quadratic dependence in $|\nabla u|$ (say, $q<2$ in~(\ref{eqv4})) and to~\cite{adp,fm00,fpr99,js,sirakov} for further existence results for problems of the type~(\ref{eqv4}),~(\ref{eqv4eps}),~(\ref{eqv5}), or more general ones, when~$\hat{\delta}\le0$ and $\hat{f}$, $\hat{f}_{\epsilon}$ or $\overline{f}$ are small in some spaces. However, it is worth pointing out that the existence and the uniqueness are not always guaranteed for general functions~$\hat{f}$, $\hat{f}_{\epsilon}$ or $\overline{f}$ if $\hat{\delta}\le0$, see in particular~\cite{adp,hbv,js,Porretta,sirakov} for problems of the type~(\ref{eqv4}),~(\ref{eqv4eps}),~(\ref{eqv5}) or for more general problems.

Comparisons between a solution $u$ of~(\ref{equ}) and a solution $v$ of a problem of the type~(\ref{eqv4}) when $H(x,s,p)$ is nonlinear in $p$ and grows at most quadratically in $|p|$ were first established by Alvino, Trombetti and Lions~\cite{atl} in the case $H=H(x,p)$ with
$$A\ge{\rm{Id}},\ \ \Lambda=1,\ \ |H(x,p)|\le f(x)+\kappa\,|p|^q$$
and $\kappa>0$, $f\in L^{\infty}(\Omega)$, $f\ge0$ a.e. in $\Omega$: in this case, there holds $u^*\le v$ a.e. in~$\Omega^*$, where~$v\in H^1_0(\Omega^*)\cap L^{\infty}(\Omega^*)$ is any weak solution of $-\Delta v=f^*+\kappa|\nabla v|^q$ in $\Omega^*$, provided such a solution~$v$ exists (it does if~$\|f\|_{L^{\infty}(\Omega)}$ is small enough). We refer to~\cite{mps,messano,pasic} for further results in this direction. On the other hand, Ferone and Posteraro~\cite{fp93} (see also~\cite{fpr99}) showed that if
$$A\ge{\rm{Id}},\ \ \Lambda=1,\ \ H(x,s,p)=\hbox{div}(F)+\tilde{H}(x,s,p),\ \ |\tilde{H}(x,s,p)|\le f(x)+|p|^2$$
with $F\in(L^r(\Omega))^n$, $f\in L^{r/2}(\Omega)$ and $r>n$, then $u^*\le v$ a.e. in $\Omega^*$ for any weak solution $v^*\in H^1_0(\Omega^*)\cap L^{\infty}(\Omega^*)$ of~$-\Delta v=f^*(x)+|\nabla v|^2+\hbox{div}(\hat{F}e_r)$ in $\Omega^*$, provided such a solution~$v$ exists (it does if the norms of~$f$ and~$F$ are small enough), where $\hat{F}$ shares some common properties with the function~$\hat{a}$ appearing in Theorems~\ref{th1} and~\ref{th2}. To our best knowledge, in the case of at most quadratic growth with respect to the gradient, the only comparison result involving non-constant functions $\Lambda$ is contained in a recent paper by Tian and Li~\cite{tl}: if
$$|H(x,s,p)|\le f(x)+\kappa\,\Lambda(x)^{2/q}|p|^q$$
with $\kappa>0$ and $f\ge 0$ a.e. in $\Omega$, then $u^*\le v$ a.e. in $\Omega^*$ for any weak solution $v^*\in H^1_0(\Omega^*)\cap L^{\infty}(\Omega^*)$ of~$-\hbox{div}(\tilde{\Lambda}(x)\nabla v)=f^*(x)+\kappa\,\tilde{\Lambda}(x)^{2/q}|\nabla v|^q$ in $\Omega^*$, provided such a  solution~$v$ exists (it does if~$f$ is small enough), where $\tilde{\Lambda}$ shares some common properties with the function~$\hat{\Lambda}$ appearing in Theorems~\ref{th1} and~\ref{th2} (in~\cite{tl}, the function $\Lambda$ can even be degenerate at some points, that is $\Lambda$ is nonnegative but is not necessarily in $L^{\infty}_+(\Omega)$, and $f/\Lambda$ belongs to $L^r(\Omega)$ for some suitable~$r$). Lastly, we refer to~\cite{ms87} for comparison results where $u$ is compared to a solution~$v$, if any, of an equation whose principal part is a homogeneous nonlinear term such as the~$p$-Laplacian.

In the references of the previous paragraph, a bound on the absolute value~$|H|$ of~$H$ in~$\Omega\times\R\times\R^n$ is used, the existence of solutions $v$ of some symmetrized problems in $\Omega^*$ is assumed and the functions $\hat{H}$ of these symmetrized problems only depend on $(x,p)\in\Omega^*\times\R^n$. As already mentioned, the existence of such solutions~$v$ is guaranteed only when some norms of the function $f$ (or $f_u$) are small, since the existence of $v$ does not hold for general functions~$f$ or $f_u$. However, when $\hat{H}=\hat{H}(x,p)$, roughly speaking, one can integrate the one-dimensional equation satisfied by the radially symmetric function~$v$ and it is then possible to compare $v$ with the solution of an equation involving the Schwarz rearrangement $f^*_u$ of the function $f_u$, as in the results of Section~\ref{sec11} (see also Lemma~\ref{lemzLz} below).

On the contrary, in Theorem~\ref{th2}, only the lower bound~(\ref{hypH}) is needed and the existence (and uniqueness) of the solutions~$v$ and $v_{\epsilon}$ of~(\ref{eqv4}) and~(\ref{eqv4eps}) is actually automatically guaranteed by the condition $\hat{\delta}>0$, which follows from the additional assumption $\inf_{\Omega\times\R\times\R^n}b>0$. The counterpart of the positivity of $\hat{\delta}$ in~(\ref{defhatH}) and~(\ref{defhatHeps}) is that one cannot integrate the one-dimensional equations satisfied by the functions $v$ and $v_{\epsilon}$. Thus, we do not know if it is possible to compare these functions $v$ and $v_{\epsilon}$ to the solution of the same type of equation with $f^*_u$ instead of~$\hat{f}$ or~$\hat{f}_{\epsilon}$ in~(\ref{defhatH}) or~(\ref{defhatHeps}). However, since $\inf_{\Omega\times\R\times\R^n}b>0$ and since $\hat{\delta}$ is shown to be positive, there is no need to assume that some norms of $f_u$, $\hat{f}$ or $\hat{f}_{\epsilon}$ are small. The function $f$ in~(\ref{hypH}) can be any bounded function.

We also point out an interesting particular case of Theorem~\ref{th2}: namely, if $\Lambda=\lambda>0$, $a=\alpha\ge0$ and $f=\gamma$ are assumed to be constant, then $u^*\le v$ almost everywhere in $\Omega^*$, where $v$ is the unique $H^1_0(\Omega^*)\cap L^{\infty}(\Omega^*)$ solution of
$$\left\{\baa{rcl}
-\lambda\Delta v-\alpha\,|\nabla v|^q+\hat{\delta}\,v & = & \gamma\ \hbox{ in }\Omega^*,\vspace{3pt}\\
v & = & 0 \ \hbox{ on }\partial\Omega^*.\eaa\right.$$
Even this particular case of Theorem~\ref{th2} is actually new.

As for Theorem~\ref{th1}, the method used in the proof of Theorem~\ref{th2} is based on a symmetrization of the second-order terms. However, a special attention has to be put on the constant $\hat{\delta}$ appearing in~(\ref{eqv4}) and~(\ref{eqv4eps}), for these problems to be well-posed.

Lastly, when $\Omega$ is not a ball, an improved quantified inequality can be established. To our knowledge, such an improved inequality for problems~(\ref{equ}) and~(\ref{eqv4}) had never been obtained before, even in particular situations.

\begin{theo}\label{th2bis}
Assume that $\Omega$ is not a ball, that~$(\ref{ALambda})$ and~$(\ref{hypH})$ hold with $\Lambda\in L^{\infty}_+(\Omega)$, $1<q\le 2$ and that~\eqref{hypb} is satisfied. Then, under the notations of Theorem~$\ref{th2}$, there is a constant $\eta_u>0$, which depends on $\Omega$, $n$ and $u$, such that
\be\label{u*vetaubis}
(1+\eta_u)\,u^*\le v\hbox{ a.e. in }\Omega^*
\ee
and
\be\label{u*vetaubiseps}
\|((1+\eta_u)\,u^*-v_{\epsilon})^+\|_{L^{2^*}(\Omega^*)}\le\epsilon.
\ee
Furthermore, if $q<1+2/n$ and if there are $M>0$ and $r\ge 0$ such that $r(n-2)<n+2$ and
\be\label{hypAMquadr}\left\{\baa{l}
\|A\|_{W^{1,\infty}(\Omega)}+\|\Lambda^{-1}\|_{L^{\infty}(\Omega)}+\|a\|_{L^{\infty}(\Omega\times\R\times\R^n)}+\|f\|_{L^{\infty}(\Omega\times\R\times\R^n)}\le M,\vspace{3pt}\\
b(x,s,p)\ge M^{-1}>0\hbox{ for all }(x,s,p)\in\Omega\times\R\times\R^n,\vspace{3pt}\\
\left\vert H(x,s,p)-H(x,0,0)\right\vert\le M\left(\left\vert s\right\vert^r+\left\vert p\right\vert^q\right)\hbox{ for all }(x,s,p)\in\Omega\times\R\times\R^n,\vspace{3pt}\\
|H(x,s,p)|\le M\,(1+|s|^r+|p|^q)\hbox{ for all }(x,s,p)\in\Omega\times\R\times\R^n,\vspace{3pt}\\
H(x,0,0)\le 0\hbox{ for all }x\in\Omega,\ \ \displaystyle\int_{\Omega}H(x,0,0)\,dx\le-M^{-1}<0,\eaa\right.
\end{equation}
then there is a constant $\eta>0$ depending only on $\Omega$, $n$, $q$, $M$ and $r$ such that
\be\label{u*vetabis}
(1+\eta)\,u^*\le v\hbox{ a.e. in }\Omega^*
\ee
and
\be\label{u*vetabiseps}
\|((1+\eta)\,u^*-v_{\epsilon})^+\|_{L^{2^*}(\Omega^*)}\le\epsilon.
\ee
\end{theo}

\noindent{\bf{Outline of the paper.}} In Section~\ref{rearrangement}, we set the basic ingredients for the proof of the main theorems. Namely, we prove some pointwise and differential rearrangement inequalities using a symmetrization of the second-order terms of~(\ref{equ}). Section~\ref{sec3} is devoted to the proofs of Theorems~\ref{th1} and~\ref{th1bis} for equation~(\ref{equ}) when the nonlinear term $H$ is bounded from below by an at most linear function of~$|\nabla u|$. In Section~\ref{sec4}, we consider the general case $q\in(1,2]$ in the lower bound~(\ref{hypH}) of~$H$ and we do the proofs of Theorems~\ref{th2} and~\ref{th2bis}. Lastly, Section~\ref{secapp} is concerned with the proof of some existence and uniqueness results for problems of the type~(\ref{equ}) in $\Omega$ under various assumptions on $H$.


\SE{Rearrangement inequalities} \label{rearrangement}

This section is devoted to the proofs of some rearrangement inequalities in the spirit of~\cite{hnrAM}, Section~3. These pointwise estimates and partial differential inequalities are of independent interest and are thus stated in a separate section. They will then be used in the proof of the main theorems of the paper in the next sections.


\subsection{Definitions of the symmetrizations}\label{sec21}

As in \cite{hnrAM}, let $\Omega$ be a $C^2$ bounded domain of $\R^n$, let $A_{\Omega}\in C^1(\overline{\Omega},{\mathcal S}_n(\R))$, $\Lambda_{\Omega}\in C^1(\overline{\Omega})\cap L^{\infty}_+(\Omega)$ and assume that
$$A_{\Omega}(x)\geq\Lambda_{\Omega}(x)\hbox{Id}\hbox{ for all }x\in\overline{\Omega}.$$
Let $\psi$ be a $C^1(\overline{\Omega})$ function, analytic and positive in $\Omega$, such that $\hbox{div}(A_{\Omega}\nabla\psi)\in L^1(\Omega)$, $\psi=0$ on $\partial\Omega$ and $|\nabla\psi(x)|\neq 0$ for all $x\in\partial\Omega$, so that $\nu\cdot\nabla\psi<0$ on $\partial\Omega$, where $\nu$ denotes the outward unit normal to $\partial\Omega$. For the sake of completeness, let us recall some definitions and notations already introduced in \cite{hnrAM}. Set
$$M=\max_{x\in \overline{\Omega}} \psi(x).$$
For all $a\in \left[0,M\right)$, define
\be\label{defOmegaa}
\Omega_a=\big\{x\in \Omega;\ \psi(x)>a\big\}
\ee
and, for all $a\in \left[0,M\right]$,
\[
\Sigma_a=\big\{x\in \overline{\Omega};\ \psi(x)=a\big\}.
\]
Notice that the $n$-dimensional Lebesgue measure $|\Sigma_a|$ of $\Sigma_a$ is equal to $0$ for every $a\in(0,M]$ by analyticity of $\psi$ in $\Omega$, and that the $n$-dimensional Lebesgue measure of $\Sigma_0=\partial\Omega$ is also equal to $0$. The set $\{x\in\overline{\Omega};\ \nabla\psi(x)=0\}$ is included in some compact set $K\subset\Omega$, which implies that the set
$$Z=\big\{a\in[0,M];\ \exists\,x\in\Sigma_{a},\ \nabla\psi(x)=0\big\}$$
of the critical values of $\psi$ is finite (\cite{soucek}) and can then be written as $Z=\{a_{1},\cdots,a_{m}\}$ for some positive integer $m\in\N^*=\N\backslash\{0\}$. Observe also that $M\in Z$ and that $0\not\in Z$. One can then assume without loss of generality that $0<a_{1}<\cdots<a_{m}=M$. The set $Y=[0,M]\backslash Z$ of the non-critical values of $\psi$ is open relatively to $[0,M]$ and can be written as
$$Y\ =\ [0,M]\backslash Z\ =\ [0,a_{1})\cup(a_{1},a_{2})\cup\cdots\cup(a_{m-1},M).$$
Denote by $R$ the radius of $\Omega^{\ast}$, that is $\Omega^*=B_R$, where, for $s>0$, $B_s$ denotes the open Euclidean ball centered at the origin with radius $s$. For all $a\in[0,M)$, let $\rho(a)\in(0,R]$ be defined so that
\begin{equation}\label{defrho}
|\Omega_{a}|=|B_{\rho(a)}|=\alpha_n\rho(a)^n,
\end{equation}
where $\alpha_n$ is the volume of the unit ball $B_1$. The function $\rho$ is extended at $M$ by $\rho(M)=0$. It follows then from~\cite{hnrAM} (Lemma~3.1 in~\cite{hnrAM}) and the fact that $|\Sigma_a|=0$ for every $a\in[0,M]$ that $\rho$ is a continuous decreasing map from $[0,M]$ onto $[0,R]$. Lastly, call
\be\label{defE}
E=\big\{x\in\overline{\Omega^*};\ |x|\in\rho(Y)\big\}.
\ee
The set $E$ is a finite union of spherical shells, it is open relatively to $\overline{\Omega^*}$ and can be written as
$$E=\big\{x\in\R^n;\ |x|\in(0,\rho(a_{m-1}))\cup\cdots\cup(\rho(a_{2}),\rho(a_{1}))\cup(\rho(a_{1}),R]\big\}$$
with $0=\rho(a_{m})=\rho(M)<\rho(a_{m-1})<\cdots<\rho(a_{1})<R$. Notice that $0\not\in E$.\par
Let us now recall the definition of the symmetrization $\hat{\psi}$ of $\psi$ introduced in \cite{hnrAM}. To do so, we first define a symmetrization of $\Lambda_{\Omega}$. Namely, for all $r\in \rho(Y)$, set
\be\label{defG}
G(r)=\frac{\displaystyle{\int_{\Sigma_{\rho^{-1}(r)}} \left\vert \nabla\psi(y)\right\vert^{-1}d\sigma_{\rho^{-1}(r)}}}{\displaystyle{\int_{\Sigma_{\rho^{-1}(r)}} \Lambda_{\Omega}(y)^{-1}\left\vert \nabla\psi(y)\right\vert^{-1}d\sigma_{\rho^{-1}(r)}}}>0,
\ee
where $\rho^{-1}\ :\ [0,R]\to[0,M]$ denotes the reciprocal of the function $\rho$ and $d\sigma_a$ denotes the surface measure on $\Sigma_a$ for $a\in Y$. For all $x\in E$, define
\be\label{defhatlambda}
\widehat{\Lambda}(x)=G(|x|)
\ee
and set, say, $\hat{\Lambda}(x)=0$ for all $x$ in the negligible set $\overline{\Omega^*}\backslash E=\Sigma_{a_1}\cup\cdots\cup\Sigma_{a_m}$. Notice that
\be\label{ineqlambda}
0<\min_{\overline{\Omega}}\Lambda_{\Omega}\le\mathop{{\rm{ess}}\,{\rm{inf}}}_{\Omega^*}\hat{\Lambda}\le\mathop{{\rm{ess}}\,{\rm{sup}}}_{\Omega^*}\hat{\Lambda}\le\max_{\overline{\Omega}}\Lambda_{\Omega}\ \hbox{ and }\ \int_{\Omega^*}\hat{\Lambda}^{-1}=\int_{\Omega}\Lambda_{\Omega}^{-1}
\ee
from the co-area formula and the fact that $|\Sigma_a|=0$ for all $a\in[0,M]$. Furthermore, the~$L^{\infty}_+(\Omega^*)$ function $\hat{\Lambda}$ is actually of class~$C^1$ in~$E\cap\Omega^*$. Define now $F(0)=0$ and, for all~$r\in\rho(Y)$, set
\be\label{defF}
F(r)=\frac{1}{n\alpha_nr^{n-1}G(r)}\int_{\Omega_{\rho^{-1}(r)}}\div(A_{\Omega}\nabla\psi)(x)\,dx.
\ee
This definition makes sense when $r\in\rho(Y)\backslash\{R\}$ since $A_{\Omega}\nabla\psi$ is of class $C^1$ in $\Omega$, and also when~$r=R$ since $\hbox{div}(A_{\Omega}\nabla\psi)$ is assumed to be in $L^1(\Omega)$. Let $\nu_a$ denote the outward unit normal to $\Omega_a$ for~$a\in Y$. From Green-Riemann formula, there holds
\be\label{eqF}\baa{rcl}
F(r) & = & \displaystyle\frac{1}{n\alpha_nr^{n-1}G(r)}\int_{\Sigma_{\rho^{-1}(r)}}A_{\Omega}(y)\nabla\psi(y)\cdot\nu_{\rho^{-1}(r)}(y)\,d\sigma_{\rho^{-1}(r)}\vspace{3pt}\\
& = & \displaystyle\frac{-1}{n\alpha_nr^{n-1}G(r)}\int_{\Sigma_{\rho^{-1}(r)}}|\nabla\psi(y)|\,A_{\Omega}(y)\nu_{\rho^{-1}(r)}(y)\cdot\nu_{\rho^{-1}(r)}(y)\,d\sigma_{\rho^{-1}(r)}<0\eaa
\ee
for all $r\in\rho(Y)\backslash\{R\}$, as well as for $r=R$ from Lebesgue's theorem and the smoothness of~$\partial\Omega$. From~(\ref{defF}) and~(\ref{eqF}), it follows that the function $F$ is actually continuous and bounded on the set
$$\rho(Y)\cup\{0\}=[0,\rho(a_{m-1}))\cup\cdots\cup(\rho(a_2),\rho(a_1))\cup(\rho(a_1),R]$$
(the continuity at $r=0$ follows from definition~(\ref{defF}) and the fact that $|\Omega_{\rho^{-1}(r)}|=|B_r|=\alpha_nr^n$). Finally, for all $x\in\overline{\Omega^*}$, set
\be\label{defhatpsi}
\hat{\psi}(x)=-\int_{|x|}^RF(r)dr.
\ee
We recall from~\cite{hnrAM} that $\hat{\psi}$ is positive in $\Omega^*$, equal to zero on $\partial\Omega^*=\partial B_R$, radially symmetric, decreasing with respect to $\left\vert x\right\vert$ in $\overline{\Omega^*}$, continuous in $\overline{\Omega^*}$, of class $C^1$ in $E\cup\left\{0\right\}$, of class~$C^2$ in~$E\cap\Omega^*$, and that $\hat{\psi}\in W^{1,\infty}(\Omega^*)\cap H^1_0(\Omega^{\ast})$.\par
Throughout this Section~\ref{rearrangement}, we are also given a real number $q$ such that
$$1\le q\le 2,$$
and two continuous functions $a_{\Omega}$ and $f_{\Omega}$ in $\overline{\Omega}$. We now define symmetrizations of the coefficients~$a_{\Omega}$ and $f_{\Omega}$. For all $x\in E$, define $\widehat{a}(x)$ by
\be\label{defhata}\widehat{a}(x)=\left\{\baa{ll}
\displaystyle\!\max_{y\in \Sigma_{\rho^{-1}(\left\vert x\right\vert)}} \left(a_{\Omega}^+(y)\Lambda_{\Omega}^{-1}(y)\right)\times\widehat{\Lambda}(x) & \!\!\hbox{if }q=2,\vspace{3pt}\\
\displaystyle\!\left(\frac{\displaystyle{\int_{\Sigma_{\rho^{-1}(|x|)}}\!\!\!a_{\Omega}^+(y)^{\frac 2{2-q}}\Lambda_{\Omega}(y)^{-\frac{q}{2-q}}\vert\nabla\psi(y)\vert^{-1}\,d\sigma_{\rho^{-1}(|x|)}}}{\displaystyle{\int_{\Sigma_{\rho^{-1}(|x|)}} \vert\nabla\psi(y)\vert^{-1}\,d\sigma_{\rho^{-1}(|x|)}}}\right)^{\frac{2-q}{2}}\!\!\!\!\!\!\times\widehat{\Lambda}(x)^{\frac{q}{2}} & \!\!\hbox{if }1\leq q<2,\eaa\right.
\ee
and $\hat{f}(x)$ by
\be\label{defhatf}
\widehat{f}(x)=\frac{\displaystyle{\int_{\Sigma_{\rho^{-1}(|x|)}} f_{\Omega}(y)\,\vert\nabla\psi(y)\vert^{-1}d\sigma_{\rho^{-1}(|x|)}}}{\displaystyle{\int_{\Sigma_{\rho^{-1}(|x|)}} \left\vert \nabla\psi(y)\right\vert^{-1}d\sigma_{\rho^{-1}(|x|)}}}.
\ee
Note that $\widehat{a}$ and $\widehat{f}$ are defined almost everywhere in $\overline{\Omega^{\ast}}$ (they can be extended by, say, $0$ on~$\overline{\Omega^*}\backslash E$).\par
Let us list here a few basic properties satisfied by the functions $\hat{a}$ and~$\hat{f}$, which will be used later in Section~\ref{sec3}. Firstly, the functions $\hat{a}$ and $\hat{f}$ are continuous in $E$. Secondly, from~(\ref{defhatlambda}), it follows immediately that, when $q=2$,
\be\label{hatamin}
\hat{a}(x)\ge\min_{\overline{\Omega}}a_{\Omega}^+\ \hbox{ for all }x\in E.
\ee
When $1\le q<2$, the H\"older inequality yields, for all $x\in E$,
$$\baa{rcl}
\displaystyle\int_{\Sigma_{\rho^{-1}(|x|)}}\Lambda_{\Omega}(y)^{-1}|\nabla\psi(y)|^{-1}d\sigma_{\rho^{-1}(|x|)} & \le & \displaystyle\left(\int_{\Sigma_{\rho^{-1}(|x|)}}\Lambda_{\Omega}(y)^{-\frac{q}{2-q}}|\nabla\psi(y)|^{-1}d\sigma_{\rho^{-1}(|x|)}\right)^{\frac{2-q}{q}}\vspace{3pt}\\
& & \displaystyle\times\left(\int_{\Sigma_{\rho^{-1}(|x|)}}|\nabla\psi(y)|^{-1}d\sigma_{\rho^{-1}(|x|)}\right)^{\frac{2q-2}{q}},\eaa$$
whence
$$\baa{rcl}
\hat{a}(x)\,\hat{\Lambda}(x)^{-\frac{q}{2}} & \ge & \displaystyle\Big(\min_{\overline{\Omega}}a_{\Omega}^+\Big)\times\left(\displaystyle\frac{\displaystyle\int_{\Sigma_{\rho^{-1}(|x|)}}\Lambda_{\Omega}(y)^{-\frac{q}{2-q}}|\nabla\psi(y)|^{-1}d\sigma_{\rho^{-1}(|x|)}}{\displaystyle\int_{\Sigma_{\rho^{-1}(|x|)}}|\nabla\psi(y)|^{-1}d\sigma_{\rho^{-1}(|x|)}}\right)^{\frac{2-q}{2}}\vspace{3pt}\\
& \ge & \displaystyle\Big(\min_{\overline{\Omega}}a_{\Omega}^+\Big)\times\left(\displaystyle\frac{\displaystyle\int_{\Sigma_{\rho^{-1}(|x|)}}\Lambda_{\Omega}(y)^{-1}|\nabla\psi(y)|^{-1}d\sigma_{\rho^{-1}(|x|)}}{\displaystyle\int_{\Sigma_{\rho^{-1}(|x|)}}|\nabla\psi(y)|^{-1}d\sigma_{\rho^{-1}(|x|)}}\right)^{\frac{q}{2}}\vspace{3pt}\\
& = & \displaystyle\Big(\min_{\overline{\Omega}}a_{\Omega}^+\Big)\times\hat{\Lambda}(x)^{-\frac{q}{2}}\eaa$$
from~(\ref{defhatlambda}). Therefore,~(\ref{hatamin}) holds for $1\le q<2$ as well. As for the upper bound of $\hat{a}$, it follows immediately from~(\ref{ineqlambda}) and~(\ref{defhata}) that, when $q=2$,
$$\hat{a}(x)\le\Big(\max_{\overline{\Omega}}a_{\Omega}^+\Big)\times\Big(\max_{\overline{\Omega}}\Lambda_{\Omega}^{-1}\Big)\times\Big(\max_{\overline{\Omega}}\Lambda_{\Omega}\Big)=\Big(\max_{\overline{\Omega}}a_{\Omega}^+\Big)\times\displaystyle\frac{\max_{\overline{\Omega}}\Lambda_{\Omega}}{\min_{\overline{\Omega}}\Lambda_{\Omega}}$$
for all $x\in E$. When $1\le q<2$, we get from~(\ref{defhatlambda}),~(\ref{ineqlambda}),~(\ref{defhata}) and by writing $\Lambda_{\Omega}(y)^{-q/(2-q)}=\Lambda_{\Omega}(y)^{-2(q-1)/(2-q)}\times\Lambda_{\Omega}(y)^{-1}$, that, for all $x\in E$,
$$\baa{rcl}
\hat{a}(x)\,\hat{\Lambda}(x)^{-\frac{q}{2}} & \!\!\le\!\! & \displaystyle\Big(\max_{\overline{\Omega}}a_{\Omega}^+\Big)\times\Big(\max_{\overline{\Omega}}\Lambda_{\Omega}^{-(q-1)}\Big)\times\left(\displaystyle\frac{\displaystyle\int_{\Sigma_{\rho^{-1}(|x|)}}\Lambda_{\Omega}(y)^{-1}|\nabla\psi(y)|^{-1}d\sigma_{\rho^{-1}(|x|)}}{\displaystyle\int_{\Sigma_{\rho^{-1}(|x|)}}|\nabla\psi(y)|^{-1}d\sigma_{\rho^{-1}(|x|)}}\right)^{\frac{2-q}{2}}\vspace{3pt}\\
& \!\!=\!\! & \displaystyle\Big(\max_{\overline{\Omega}}a_{\Omega}^+\Big)\times\Big(\max_{\overline{\Omega}}\Lambda_{\Omega}^{-(q-1)}\Big)\times\hat{\Lambda}(x)^{-\frac{2-q}{2}},\eaa$$
whence
$$\hat{a}(x)\le\Big(\max_{\overline{\Omega}}a_{\Omega}^+\Big)\times\Big(\max_{\overline{\Omega}}\Lambda_{\Omega}^{-(q-1)}\Big)\times\hat{\Lambda}(x)^{q-1}\le\Big(\max_{\overline{\Omega}}a_{\Omega}^+\Big)\times\left(\displaystyle\frac{\max_{\overline{\Omega}}\Lambda_{\Omega}}{\min_{\overline{\Omega}}\Lambda_{\Omega}}\right)^{q-1}.$$
To sum up, there holds
\be\label{ineqa}
\min_{\overline{\Omega}}a_{\Omega}^+\le\mathop{{\rm{ess}}\,{\rm{inf}}}_{\Omega^*}\hat{a}\le\mathop{{\rm{ess}}\,{\rm{sup}}}_{\Omega^*}\hat{a}\le\Big(\max_{\overline{\Omega}}a_{\Omega}^+\Big)\times\left(\displaystyle\frac{\max_{\overline{\Omega}}\Lambda_{\Omega}}{\min_{\overline{\Omega}}\Lambda_{\Omega}}\right)^{q-1}
\ee
in all cases $1\le q\le 2$. We also point out that
\be\label{ineqf}
\min_{\overline{\Omega}}f_{\Omega}\le\mathop{{\rm{ess}}\,{\rm{inf}}}_{\Omega^*}\hat{f}\le\mathop{{\rm{ess}}\,{\rm{sup}}}_{\Omega^*}\hat{f}\le\max_{\overline{\Omega}}f_{\Omega}\ \hbox{ and }\ \int_{\Omega^*}\hat{f}=\int_{\Omega}f_{\Omega}.
\ee
from~(\ref{defhatf}), the co-area formula and the fact that $|\Sigma_a|=0$ for all $a\in[0,M]$.\par
Lastly, we are given a nonnegative continuous function $b_{\Omega}$ in $\overline{\Omega}$.


\subsection{Inequalities for the symmetrized data}

Recall first that the function $\hat{\psi}$ satisfies the following key inequality (see \cite{hnrAM}, Corollary 3.6):

\begin{pro}\label{key1} For all $x\in\overline{\Omega^*}$ and all $y\in \Sigma_{\rho^{-1}(\left\vert x\right\vert)}$,
$$\hat{\psi}(x)\ge \psi(y)=\rho^{-1}(|x|)\ge0.$$
\end{pro}

We now establish a pointwise differential inequality involving all the symmetrizations defined in Section \ref{sec21}:

\begin{pro} \label{key2}
For all $x\in E\cap\Omega^*$, there exists $y\in \Sigma_{\rho^{-1}(\left\vert x\right\vert)}$ such that
$$\begin{array}{rl}
& \displaystyle -{\rm{div}}\big(\widehat{\Lambda}\nabla\hat{\psi}\big)(x)-\widehat{a}(x)\vert \nabla\hat{\psi}(x)\vert^q-\widehat{f}(x)\vspace{3pt}\\
\leq\!\! & \displaystyle-{\rm{div}}(A_{\Omega}\nabla\psi)(y)-a_{\Omega}(y)\left\vert \nabla\psi(y)\right\vert^q-f_{\Omega}(y)\vspace{3pt}\\
\leq\!\! & \displaystyle-{\rm{div}}(A_{\Omega}\nabla\psi)(y)-a_{\Omega}(y)\left\vert \nabla\psi(y)\right\vert^q+b_{\Omega}(y)\psi(y)-f_{\Omega}(y).\end{array}$$
\end{pro}

For the proof of Proposition \ref{key2}, we need the following lemma:

\begin{lem} \label{limits}
For all $x\in E$ with $\left\vert x\right\vert=r$, there holds
\be\label{ineqpsi1}
\lim_{t\rightarrow 0^+} \frac{\displaystyle \int_{\Omega_{\rho^{-1}(r)}\setminus \Omega_{\rho^{-1}(r-t)}} a_{\Omega}(y) \left\vert \nabla\psi(y)\right\vert^qdy}{\left\vert \Omega_{\rho^{-1}(r)}\setminus \Omega_{\rho^{-1}(r-t)}\right\vert} \leq \widehat{a}(x)\,\big\vert \nabla\hat{\psi}(x)\big\vert^q
\ee
and
\be\label{ineqpsi2}
\lim_{t\rightarrow 0^+} \frac{\displaystyle \int_{\Omega_{\rho^{-1}(r)}\setminus \Omega_{\rho^{-1}(r-t)}} f_{\Omega}(y)\,dy}{\left\vert \Omega_{\rho^{-1}(r)}\setminus \Omega_{\rho^{-1}(r-t)}\right\vert} = \widehat{f}(x).
\ee
\end{lem}

\noindent{\bf Proof.} Let $x\in E$ with $\left\vert x\right\vert=r$. Notice that $\left\{z\in \R^n;\ r-\eta\le\left\vert z\right\vert\le r\right\}\subset E$ for some~$\eta>0$, and that formula~(\ref{ineqpsi2}) is an immediate consequence of the co-area formula and the definition~(\ref{defhatf}) of $\widehat{f}$.\par
For the proof of~(\ref{ineqpsi1}), consider first the case where $q=2$. By the co-area formula,
$$
\lim_{t\rightarrow 0^+} \frac{\displaystyle \int_{\Omega_{\rho^{-1}(r)}\setminus \Omega_{\rho^{-1}(r-t)}} a_{\Omega}(y) \left\vert \nabla\psi(y)\right\vert^2dy}{\left\vert \Omega_{\rho^{-1}(r)}\setminus \Omega_{\rho^{-1}(r-t)}\right\vert} =\frac{\displaystyle \int_{\Sigma_{\rho^{-1}(r)}} a_{\Omega}(y)\left\vert \nabla\psi(y)\right\vert d\sigma_{\rho^{-1}(r)}}{\displaystyle \int_{\Sigma_{\rho^{-1}(r)}} \left\vert \nabla\psi(y)\right\vert^{-1} d\sigma_{\rho^{-1}(r)}}.
$$
As a consequence, by the definition of $\widehat{a}$ in~(\ref{defhata}),
$$
\lim_{t\rightarrow 0^+} \frac{\displaystyle \int_{\Omega_{\rho^{-1}(r)}\setminus \Omega_{\rho^{-1}(r-t)}} a_{\Omega}(y) \left\vert \nabla\psi(y)\right\vert^2dy}{\left\vert \Omega_{\rho^{-1}(r)}\setminus \Omega_{\rho^{-1}(r-t)}\right\vert} \leq \widehat{a}(x)\,\widehat{\Lambda}(x)^{-1} \frac{\displaystyle \int_{\Sigma_{\rho^{-1}(r)}} \Lambda_{\Omega}(y)\left\vert \nabla\psi(y)\right\vert d\sigma_{\rho^{-1}(r)}}{\displaystyle \int_{\Sigma_{\rho^{-1}(r)}} \left\vert \nabla\psi(y)\right\vert^{-1} d\sigma_{\rho^{-1}(r)}}.
$$
But inequality $(3.16)$ in \cite{hnrAM} yields
\be\label{ineq316}
\frac{\displaystyle \int_{\Sigma_{\rho^{-1}(r)}} \Lambda_{\Omega}(y)\left\vert \nabla\psi(y)\right\vert d\sigma_{\rho^{-1}(r)}}{\displaystyle \int_{\Sigma_{\rho^{-1}(r)}} \left\vert \nabla\psi(y)\right\vert^{-1} d\sigma_{\rho^{-1}(r)}}\leq \widehat{\Lambda}(x)\,\big\vert \nabla\hat{\psi}(x)\big\vert^2,
\ee
which ends the proof of~(\ref{ineqpsi1}) when $q=2$.\par
Consider now the case where $1\le q<2$. Then, using the co-area formula again, one has
\begin{equation}\label{coarea}
\lim_{t\rightarrow 0^+} \frac{\displaystyle \int_{\Omega_{\rho^{-1}(r)}\setminus \Omega_{\rho^{-1}(r-t)}} a_{\Omega}(y) \left\vert \nabla\psi(y)\right\vert^qdy}{\left\vert \Omega_{\rho^{-1}(r)}\setminus \Omega_{\rho^{-1}(r-t)}\right\vert} =\frac{\displaystyle \int_{\Sigma_{\rho^{-1}(r)}} a_{\Omega}(y)\left\vert \nabla\psi(y)\right\vert^{q-1} d\sigma_{\rho^{-1}(r)}}{\displaystyle \int_{\Sigma_{\rho^{-1}(r)}} \left\vert \nabla\psi(y)\right\vert^{-1} d\sigma_{\rho^{-1}(r)}}.
\end{equation}
The H\"older inequality yields
$$
\begin{array}{lll}
\displaystyle \int_{\Sigma_{\rho^{-1}(r)}}\!\!\!\!\!a_{\Omega}(y)\left\vert \nabla\psi(y)\right\vert^{q-1} d\sigma_{\rho^{-1}(r)} &  \!\!\le\!\! & \displaystyle \int_{\Sigma_{\rho^{-1}(r)}}\!\!\!\!\!a_{\Omega}^+(y) \Lambda_{\Omega}(y)^{-\frac q2} \left\vert \nabla\psi(y)\right\vert^{\frac q2-1}\!\Lambda_{\Omega}(y)^{\frac q2}\left\vert \nabla \psi(y)\right\vert^{\frac q2} d\sigma_{\rho^{-1}(r)} \\
& \!\!\leq\!\! & \displaystyle \left(\int_{\Sigma_{\rho^{-1}(r)}}\!\!\!a_{\Omega}^+(y)^{\frac 2{2-q}} \Lambda_{\Omega}(y)^{-\frac q{2-q}} \left\vert \nabla\psi(y)\right\vert^{-1}d\sigma_{\rho^{-1}(r)}\right)^{\frac{2-q}2} \\
& \!\!\!\! & \times \displaystyle \left(\int_{\Sigma_{\rho^{-1}(r)}}\!\!\Lambda_{\Omega}(y) \left\vert \nabla \psi(y)\right\vert d\sigma_{\rho^{-1}(r)}\right)^{\frac q2},
\end{array}
$$
so that, by the definition (\ref{defhata}) of $\widehat{a}$ and by \eqref{coarea},
$$\displaystyle \lim_{t\rightarrow 0^+} \frac{\displaystyle \int_{\Omega_{\rho^{-1}(r)}\setminus \Omega_{\rho^{-1}(r-t)}}\!\!\!\!\!a_{\Omega}(y) \left\vert \nabla\psi(y)\right\vert^qdy}{\left\vert \Omega_{\rho^{-1}(r)}\setminus \Omega_{\rho^{-1}(r-t)}\right\vert}\leq\displaystyle \widehat{a}(x)\,\widehat{\Lambda}(x)^{-\frac q2}\times\!\left(\frac{\displaystyle \int_{\Sigma_{\rho^{-1}(r)}}\!\!\!\Lambda_{\Omega}(y) \left\vert \nabla \psi(y)\right\vert d\sigma_{\rho^{-1}(r)}}{\displaystyle\int_{\Sigma_{\rho^{-1}(r)}}\!\!\left\vert \nabla \psi(y)\right\vert^{-1} d\sigma_{\rho^{-1}(r)}}\right)^{\!\!\frac q2}\!\!\!.$$
Using inequality~(\ref{ineq316}), one therefore concludes that
$$
\displaystyle \lim_{t\rightarrow 0^+} \frac{\displaystyle \int_{\Omega_{\rho^{-1}(r)}\setminus \Omega_{\rho^{-1}(r-t)}} a_{\Omega}(y) \left\vert \nabla\psi(y)\right\vert^qdy}{\left\vert \Omega_{\rho^{-1}(r)}\setminus \Omega_{\rho^{-1}(r-t)}\right\vert} \leq   \widehat{a}(x)\,\big\vert \nabla\hat{\psi}(x)\big\vert^q,
$$
as claimed. The proof of Lemma~\ref{limits} is thereby complete.\hfill\fin\break

\noindent{\bf Proof of Proposition \ref{key2}.} Let $x\in E\cap\Omega^*$ with $\left\vert x\right\vert=r$ and let $\eta>0$ be such that~$\left\{z\in \R^n;\ r-\eta\le\left\vert z\right\vert\le r\right\}\subset E\cap\Omega^*$. It follows from the definition of $\hat{\psi}$ and the Green-Riemann formula that, for all $t\in(0,\eta]$ and for all $z\in\R^n$ with $|z|=1$, one has
$$\baa{l}
\displaystyle\int_{\Omega_{\rho^{-1}(r)}\setminus \Omega_{\rho^{-1}(r-t)}}\!\!\!\!\hbox{div}(A_{\Omega}\nabla\psi)(y)\,dy\vspace{3pt}\\
\qquad\qquad=n\alpha_nr^{n-1}\hat{\Lambda}(rz)\nabla\hat{\psi}(rz)\cdot z-n\alpha_n(r-t)^{n-1}\hat{\Lambda}((r-t)z)\nabla\hat{\psi}((r-t)z)\cdot z\vspace{3pt}\\
\qquad\qquad=\displaystyle\int_{B_r\setminus B_{r-t}}\hbox{div}\big(\hat{\Lambda}\nabla\hat{\psi}\big)(y)\,dy,\eaa$$
since $\hat{\Lambda}$ is radially symmetric and of class $C^1$ in $E$ and $\hat{\psi}$ is radially symmetric and of class~$C^2$ in~$E\cap\Omega^*$. Hence,
\begin{equation} \label{divaomega}
\lim_{t\rightarrow 0^+}\!\frac{\displaystyle \int_{\Omega_{\rho^{-1}(r)}\setminus \Omega_{\rho^{-1}(r-t)}}\!\!\!\!\!\!\mbox{div}(A_{\Omega}\nabla\psi)(y)\,dy}{\left\vert \Omega_{\rho^{-1}(r)}\setminus \Omega_{\rho^{-1}(r-t)}\right\vert}=\lim_{t\rightarrow 0^+}\!\frac{\displaystyle \int_{B_r\setminus B_{r-t}}\!\!\!\!\hbox{div}\big(\hat{\Lambda}\nabla\hat{\psi}\big)(y)\,dy}{\left\vert B_r\setminus B_{r-t}\right\vert}=\mbox{div}\big(\widehat{\Lambda}\nabla\hat{\psi}\big)(x).
\end{equation}
Gathering \eqref{divaomega} and Lemma \ref{limits}, one obtains
$$
\begin{array}{l}
\displaystyle \lim_{t\rightarrow 0^+}Ê\frac{\displaystyle \int_{\Omega_{\rho^{-1}(r)}\setminus \Omega_{\rho^{-1}(r-t)}} \big(\mbox{div}(A_{\Omega}\nabla\psi)(y)+a_{\Omega}(y)\left\vert \nabla\psi(y)\right\vert^q+f_{\Omega}(y)\big)dy}{\left\vert \Omega_{\rho^{-1}(r)}\setminus \Omega_{\rho^{-1}(r-t)}\right\vert}\vspace{3pt}\\
\qquad\qquad\qquad\qquad\qquad\qquad\qquad\qquad\qquad\qquad\leq  \mbox{div}\big(\widehat{\Lambda}\nabla\hat{\psi}\big)(x)+\widehat{a}(x)\big\vert \nabla\hat{\psi}(x)\big\vert^q+\widehat{f}(x).
\end{array}
$$
Arguing as in the proof of Proposition 3.8 in \cite{hnrAM}, one therefore obtains the existence of a point $y\in \Sigma_{\rho^{-1}(r)}$ such that
$$
\displaystyle \mbox{div}(A_{\Omega}\nabla\psi)(y)+a_{\Omega}(y)\left\vert \nabla\psi(y)\right\vert^q+f_{\Omega}(y)
\leq\mbox{div}\big(\widehat{\Lambda}\nabla\hat{\psi}\big)(x)+\widehat{a}(x)\big\vert \nabla\hat{\psi}(x)\big\vert^q+\widehat{f}(x).
$$
Since both functions $b_{\Omega}$ and $\psi$ are nonnegative, the conclusion of Proposition \ref{key2} readily follows.\hfill\fin


\subsection{An improved inequality when $\min_{\overline{\Omega}} b_{\Omega}>0$}

Let us now state an improved version of Proposition \ref{key2}, assuming especially that $b_{\Omega}$ is bounded from below on $\overline{\Omega}$ by a positive constant. For all $N>0$ and all $\beta>0$, let $E_{N,\beta}(\Omega)$ be the set of functions $\phi\in C^1(\overline{\Omega})$ such that
$$\phi=0\hbox{ on }\partial\Omega,\ \left\Vert \phi\right\Vert_{C^1(\overline{\Omega})}\leq N\mbox{ and } \phi(x)\geq \beta\,d(x,\partial\Omega)\ge0\mbox{ for all }x\in \overline{\Omega},$$
where $d(\cdot,\partial\Omega)$ denotes the Euclidean distance to $\partial\Omega$ and
$$\|\phi\|_{C^1(\overline{\Omega})}=\|\phi\|_{L^{\infty}(\Omega)}+\|\,|\nabla\phi|\,\|_{L^{\infty}(\Omega)}.$$

\begin{pro} \label{key2bis}
In addition to the general assumptions of Section~$\ref{sec21}$, assume that
\be\label{defmbomega}
\min_{\overline{\Omega}} b_{\Omega}\ge m_b>0,
\ee
that $m_{\Lambda}>0$, $M_a\geq 0$, $M_f\geq 0$, $N>0$ and $\beta>0$ are such that
$$\min_{\overline{\Omega}}\Lambda_{\Omega}\ge m_{\Lambda}>0,\ \ \|a_{\Omega}^+\|_{L^{\infty}(\Omega)}\le M_a,\ \ \|f_{\Omega}^+\|_{L^{\infty}(\Omega)}\le M_f,$$
$\psi\in E_{N,\beta}(\Omega)$, and that there exists $\kappa\geq 0$ such that
\begin{equation} \label{subsol}
-{\rm{div}}(A_{\Omega}\nabla \psi)(y)-a_{\Omega}(y)\left\vert \nabla \psi(y)\right\vert^q+b_{\Omega}(y)\psi(y)-f_{\Omega}(y)\leq\kappa\mbox{ in }\Omega.
\ee
Then there exists a constant $\delta>0$ only depending on $\Omega$, $n$, $m_b$, $m_{\Lambda}$, $M_a$, $M_f$, $N$, $\beta$ and $\kappa$, with the following property: for all $x\in E\cap\Omega^*$, there exists $y\in \Sigma_{\rho^{-1}(\left\vert x\right\vert)}$ such that
$$
\begin{array}{rl}
& \displaystyle -{\rm{div}}\big(\widehat{\Lambda}\nabla\hat{\psi}\big)(x)-\widehat{a}(x)\big\vert \nabla\hat{\psi}(x)\big\vert^q+\delta\hat{\psi}(x)-\widehat{f}(x)\vspace{3pt}\\
\leq\!\! & -{\rm{div}}(A_{\Omega}\nabla\psi)(y)-a_{\Omega}(y)\left\vert \nabla\psi(y)\right\vert^q+b_{\Omega}(y)\psi(y)-f_{\Omega}(y).
\end{array}
$$
\end{pro}

The proof relies on the following observation:

\begin{lem} \label{ineqpsi}
Under the assumptions of Proposition~$\ref{key2bis}$, there exists a constant $\hat{\delta}>0$ only depending on $\Omega$, $n$, $m_{\Lambda}$, $M_a$, $M_f$, $N$, $\beta$ and $\kappa$, such that, for all $x\in \overline{\Omega^{\ast}}$ and all $y\in \Sigma_{\rho^{-1}(\left\vert x\right\vert)}$, there holds
$$
\hat{\delta}\,\hat{\psi}(x)\leq \psi(y).
$$
\end{lem}

\noindent{\bf Proof.} We first claim that there exists $\gamma>0$ only depending on $\Omega$ and $\beta$ such that, for all~$\phi\in E_{N,\beta}(\Omega)$ and all $a\in [0,\max_{\overline{\Omega}}\phi]$,
\begin{equation} \label{claimpsi}
0\leq \alpha_n(R^n-(\rho_{\phi}(a))^n)\leq \gamma a,
\end{equation}
where $\rho_{\phi}$ is the function $\rho$ associated with $\phi$, as defined in~\eqref{defOmegaa} and~\eqref{defrho} with $\phi$ instead of~$\psi$. Indeed, for $\phi\in E_{N,\beta}(\Omega)$, one has
$$
\left\vert \Omega_a\right\vert\geq \left\vert \left\{y\in \Omega;\ d(y,\partial\Omega)> \frac a{\beta}\right\}\right\vert\geq \left\vert \Omega\right\vert-\gamma a,
$$
using the fact that $\Omega$ is of class $C^1$. This yields the claim. \par
Let us now prove that there exists $\eta>0$ only depending on $\Omega$, $N$ and $\beta$ such that, for all~$\phi\in E_{N,\beta}(\Omega)$, all $x\in \overline{\Omega^{\ast}}$ and all $y\in \Sigma_{(\rho_{\phi})^{-1}(\left\vert x\right\vert)}$,
\begin{equation} \label{compardistance}
d(y,\partial\Omega)\geq \eta\,d(x,\partial\Omega^{\ast}).
\end{equation}
Let us assume by contradiction that this is not true. Then, there exist a sequence of functions~$(\phi_k)_{k\geq 1}\in E_{N,\beta}(\Omega)$ and two sequences of points $(x_k)_{k\geq 1}\in \overline{\Omega^{\ast}}$ and $(y_k)_{k\geq 1}\in \overline{\Omega}$ with
\begin{equation} \label{conditionyk}
y_k\in \Sigma_{(\rho_{\phi_k})^{-1}(\left\vert x_k\right\vert)}\mbox{ and } d(y_k,\partial\Omega)< \frac{d(x_k,\partial\Omega^{\ast})}{k}
\end{equation}
for all $k\ge 1$. This implies that $d(y_k,\partial\Omega)\rightarrow 0$, and since the $C^1(\overline{\Omega})$ norms of the functions~$\phi_k$ are uniformly bounded,  it follows that $\phi_k(y_k)\rightarrow 0$. Applying \eqref{claimpsi} with $a=\phi_k(y_k)$, one obtains that
$$
\rho_{\phi_k}(\phi_k(y_k))\rightarrow R,
$$
that is $\left\vert x_k\right\vert=\rho_{\phi_k}(\phi_k(y_k))\rightarrow R$ when $k\to+\infty$. Using again the uniform bound for the $C^1(\overline{\Omega})$ norms of $\phi_k$, one has, for all $k$ large enough,
$$
\phi_k(y_k)\leq Nd(y_k,\partial\Omega),
$$
whence
$$
\phi_k(y_k)\leq \frac{Nd(x_k,\partial\Omega^{\ast})}{k}=\frac{N(R-\left\vert x_k\right\vert)}{k}.
$$
Applying $\rho_{\phi_k}$ to both sides of this inequality and using the fact that $\rho_{\phi_k}$ is nonincreasing, it follows that
$$
\left\vert x_k\right\vert\geq \rho_{\phi_k}\left(\frac Nk(R-\left\vert x_k\right\vert)\right)
$$
for all $k$ large enough. Using \eqref{claimpsi} again, one easily deduces that, for all $k$ large enough,
$$
\alpha_n\left\vert x_k\right\vert^n\geq \alpha_nR^n-\frac{\gamma N}k (R-\left\vert x_k\right\vert),
$$
that is, for all $k$ large enough,
$$
\alpha_n \frac{R^n-\left\vert x_k\right\vert^n}{R-\left\vert x_k\right\vert} \leq \frac{\gamma N}k
$$
(note that, by \eqref{conditionyk}, $\left\vert x_k\right\vert<R$ for all $k\geq 1$), and this provides a contradiction when $k\rightarrow +\infty$ since $\left\vert x_k\right\vert\rightarrow R$. Thus, \eqref{compardistance} is proved.\par
Let us now end up the proof of Lemma \ref{ineqpsi}. By \eqref{compardistance}, one has
\begin{equation} \label{minorpsi}
\psi(y)\geq \beta\,d(y,\partial\Omega)\geq \beta\,\eta\,d(x,\partial\Omega^{\ast})=\beta\,\eta\,(R-\left\vert x\right\vert)
\end{equation}
for all $x\in\overline{\Omega^*}$ and $y\in\Sigma_{\rho^{-1}(|x|)}$. But, using~(\ref{defG}),~(\ref{defF}),~(\ref{subsol}) and the nonnegativity of $b_{\Omega}$ and $\psi$, one has, for all $r\in\rho(Y)$,
$$\baa{rcl}
-F(r) & \le & \displaystyle\frac{1}{n\alpha_nr^{n-1}G(r)}\int_{\Omega_{\rho^{-1}(r)}}\!\!\big(\|a_{\Omega}^+\|_{L^{\infty}(\Omega)}N^q+\kappa+\|f_{\Omega}^+\|_{L^{\infty}(\Omega)}\big)\,dy\vspace{3pt}\\
& \le & \displaystyle\frac{R\,\big(M_a\max(1,N)^2+\kappa+M_f\big)}{n\,m_{\Lambda}},\eaa$$
since $\min_{\overline{\Omega}}\Lambda_{\Omega}\ge m_{\Lambda}>0$ and $|\Omega_{\rho^{-1}(r)}|=|B_r|=\alpha_nr^n$. Together with~(\ref{defhatpsi}), it follows that there exists $\theta>0$ only depending on $\Omega$, $n$, $m_{\Lambda}$, $M_a$, $M_f$, $N$ and $\kappa$, such that
\begin{equation} \label{majorpsitilde}
\hat{\psi}(x)\leq \theta\,(R-\left\vert x\right\vert)
\end{equation}
for all $x\in \overline{\Omega^{\ast}}$. The conclusion of Lemma \ref{ineqpsi} then readily follows from \eqref{minorpsi} and \eqref{majorpsitilde}. \hfill\fin\break

\noindent{\bf Proof of Proposition \ref{key2bis}.} Let $x\in E\cap\Omega^*$. Proposition \ref{key2} provides the existence of a point $y\in \Sigma_{\rho^{-1}(\left\vert x\right\vert)}$ such that
\begin{equation} \label{provides}
\displaystyle -\mbox{div}\big(\widehat{\Lambda}\nabla\hat{\psi}\big)(x)-\widehat{a}(x)\big\vert \nabla\hat{\psi}(x)\big\vert^q-\widehat{f}(x)\leq-\mbox{div}(A_{\Omega}\nabla\psi)(y)-a_{\Omega}(y)\left\vert \nabla\psi(y)\right\vert^q-f_{\Omega}(y).
\end{equation}
Now, it follows from Lemma \ref{ineqpsi} and~(\ref{defmbomega}) that, with $\delta=\hat{\delta}\,m_b>0$, one has
\begin{equation} \label{psitildepsi}
\delta\,\hat{\psi}(x)\leq m_b\,\psi(y)\le b_{\Omega}(y)\,\psi(y),
\end{equation}
and it is therefore enough to sum up \eqref{provides} and \eqref{psitildepsi} to obtain the conclusion of Proposition~\ref{key2bis}.\hfill\fin


\subsection{The case where $\Omega$ is not a ball}\label{sec24}

Let us finally recall that Proposition \ref{key1} can be improved when $\Omega$ is not a ball. Following Section~4 of~\cite{hnrAM}, for all $\alpha\in (0,1)$, all $N>0$ and all $\beta>0$, let $E_{\alpha,N,\beta}(\Omega)$ be the set of functions $\phi\in E_{N,\beta}(\Omega)\cap C^{1,\alpha}(\overline{\Omega})$ such that $\left\Vert \phi\right\Vert_{C^{1,\alpha}(\overline{\Omega})}\leq N$, where
$$\|\phi\|_{C^{1,\alpha}(\overline{\Omega})}=\|\phi\|_{C^1(\overline{\Omega})}+\sup_{(x,y)\in\Omega\times\Omega,\,x\neq y}\frac{|\nabla\phi(x)-\nabla\phi(y)|}{|x-y|^{\alpha}}.$$
It was established in \cite{hnrAM} (Corollary~4.4 in~\cite{hnrAM}) that:

\begin{pro}\label{key1notball}
In addition to the general assumptions of Section~$\ref{sec21}$, assume that $\psi\in E_{\alpha,N,\beta}(\Omega)$ for some $\alpha\in(0,1)$, $N>0$, $\beta>0$, and that $\Omega$ is not a ball. Then, there exists~$\eta>0$ only depending on $\Omega$, $n$, $\alpha$, $N$ and $\beta$ such that, for all $x\in\overline{\Omega^*}$ and all $y\in \Sigma_{\rho^{-1}(\left\vert x\right\vert)}$,
$$\hat{\psi}(x)\ge (1+\eta)\,\psi(y).$$
\end{pro}


\SE{Linear growth with respect to the gradient}\label{sec3}

This section is devoted to the proofs of Theorems~\ref{th1} and~\ref{th1bis}. The proofs of some technical lemmas used in the proofs of these theorems are done in Section~\ref{sec32}. The proofs of Theorems~\ref{th1} and~\ref{th1bis}, that are done in Section~\ref{sec31}, follow the same general scheme. As a matter of fact, the only difference in the conclusions~(\ref{u*v}),~(\ref{u*vetau}) and~(\ref{u*veta}) is that the inequalities~(\ref{u*vetau}) and~(\ref{u*veta}) are quantified when $\Omega$ is not a ball, in that they involve a parameter~$\eta_u>0$ (resp.~$\eta>0$) which depends on $\Omega$, $n$ and $u$ (resp. $\Omega$, $n$ and~$M$ given in~(\ref{hypAM})). Most of the steps of the proofs of Theorems~\ref{th1} and~\ref{th1bis} will be identical, this is the reason why the proofs are done simultaneously. However, in some steps or in some arguments, we will consider specifically the case where $\Omega$ is not ball and where the assumption~(\ref{hypAM}) is made. Some more precise estimates will be proved in this case.\par
Throughout this section, we assume~(\ref{ALambda}) and~(\ref{hypH}) with $\Lambda\in L^{\infty}_+(\Omega)$ and
$$q=1,$$
that is the nonlinear $H$ is bounded from below by an at most linear function of~$|p|$. Furthermore, $u\in W(\Omega)$ denotes a solution of~(\ref{equ}) such that
\be\label{hypu}
u>0\hbox{ in }\Omega\hbox{ and }|\nabla u|\neq 0\hbox{ on }\partial\Omega.
\ee
We recall that, even if it means redefining $u$ on a negligible subset of~$\overline{\Omega}$, one can assume without loss of generality that $u\in C^{1,\alpha}(\overline{\Omega})$ for all $\alpha\in[0,1)$. 


\subsection{Proofs of Theorems~\ref{th1} and~\ref{th1bis}}\label{sec31}

The preliminary step (Step~1) of the proofs is concerned with some uniform bounds on $u$, depending only on $\Omega$, $n$ and~$M$, under the assumption~(\ref{hypAM}). These bounds, which are of independent interest, will be used later in the specific case where $\Omega$ is not a ball. Then, the general strategy consists firstly in approximating~$u$ in~$\Omega$ by smooth solutions~$u_j$ of some regularized equations (Step~2) and then in applying the rearrangement inequalities of Section~\ref{rearrangement} to the approximated solutions~$u_j$ and the coefficients appearing in the approximated equations (Step~3). In Step~3, these rearrangement inequalities are quantified when $\Omega$ is not a ball and~(\ref{hypAM}) is assumed. The ideas used in the next steps of the proofs are identical for Theorems~\ref{th1} and~\ref{th1bis}. More precisely, in Steps~4 and~5, we apply a maximum principle to the symmetrized functions in $\Omega^*$, called $\hat{\psi}_k=\hat{u_{j_k}}$, namely we compare them to the solutions~$v_k$ of some radially symmetric equations in $\Omega^*$. We then pass to the limit as~$k\to+\infty$ in~$\Omega^*$ (Steps~6 and~7). We also approximate the symmetrized coefficients~$\hat{f}_k$ in~$\Omega^*$ appearing in the proof by some functions in~$\Omega^*$ having the same distribution function as the function $f_u$ defined in~(\ref{deffu}) (Steps~8 and~9). Finally, in Steps~10 and~11, we pass to some limits and we use the Hardy-Littlewood inequality to compare some approximated solutions in $\Omega^*$ with the solution~$v$ of~(\ref{eqvbis}).

\subsubsection*{Step 1: uniform bounds on $u$ under assumption~(\ref{hypAM})}

In this step, some uniform pointwise and smoothness estimates  are established under assumption~(\ref{hypAM}). Actually, these quantified estimates will only be needed for the quantified inequality~(\ref{u*veta}) in Theorem~\ref{th1bis}. We recall that the sets $E_{\alpha,N,\beta}(\Omega)$ have been defined in Section~\ref{sec24}

\begin{lem}\label{lembounds}
Under assumption~$(\ref{hypAM})$, there are some real numbers $N>0$ and $\beta>0$, which depend only on $\Omega$, $n$ and $M$, such that $u\in E_{1/2,N,\beta}(\Omega)$.
\end{lem}

The proof of this lemma, which has its independent interest, is postponed in Section~\ref{sec32}. We prefer to directly go in the sequel on the main steps of the proofs of Theorems~\ref{th1} and~\ref{th1bis}.

\subsubsection*{Step 2: approximated coefficients and approximated solutions $u_j$ in $\Omega$}

Let $H_{\infty}:\overline{\Omega}\to\R$ be the continuous function defined by
$$H_{\infty}(x)=H(x,u(x),\nabla u(x))\ \hbox{ for all }x\in\overline{\Omega}$$
and let $(H_j)_{j\in\N}$ be a sequence of polynomial functions such that
\be\label{Hkinfty}
H_j(x)\mathop{\longrightarrow}_{j\to+\infty}H_{\infty}(x)=H(x,u(x),\nabla u(x))\ \hbox{ uniformly in }x\in\overline{\Omega}.
\ee\par
We recall that the given matrix field $A=(A_{i,i'})_{1\le i,i'\le n}$ is in $W^{1,\infty}(\Omega,{\mathcal{S}}_n(\R))$ and that all entries $A_{i,i'}$ can be assumed to be continuous in $\overline{\Omega}$ without loss of generality. Now, following Steps~1 and~2 of Section~5.2.1 of~\cite{hnrAM}, there is a sequence of~$C^{\infty}(\overline{\Omega},{\mathcal{S}}_n(\R))$ matrix fields~$(A_j)_{j\in\N}=((A_{j;i,i'})_{1\le i,i'\le n})_{j\in\N}$ with polynomial entries $A_{j;i,i'}$ and a sequence of~$C^{\infty}(\overline{\Omega})\cap L^{\infty}_+(\Omega)$ functions~$(\Lambda_j)_{j\in\N}$ such that
\be\label{Aklambdak}\left\{\baa{l}
\displaystyle A_{j;i,i'}\mathop{\longrightarrow}_{j\to+\infty}A_{i,i'}\hbox{ uniformly in }\overline{\Omega}\hbox{ for all }1\le i,i'\le n,\ \ \displaystyle\mathop{\sup}_{j\in\N}\|A_j\|_{W^{1,\infty}(\Omega)}<+\infty,\vspace{3pt}\\
A_j\ge\Lambda_j\hbox{Id}\hbox{ in }\overline{\Omega}\ \hbox{ and }\ \|\Lambda_j^{-1}\|_{L^1(\Omega)}=\|\Lambda^{-1}\|_{L^1(\Omega)}\ \hbox{ for all }j\in\N,\vspace{3pt}\\
0<\displaystyle\mathop{{\rm{ess}}\,{\rm{inf}}}_{\Omega}\Lambda\le\mathop{\liminf}_{j\to+\infty}\Big(\mathop{\min}_{\overline{\Omega}}\Lambda_j\Big)\le\mathop{\limsup}_{j\to+\infty}\Big(\mathop{\max}_{\overline{\Omega}}\Lambda_j\Big)\le\mathop{{\rm{ess}}\,{\rm{sup}}}_{\Omega}\Lambda\eaa\right.
\ee
(namely, one can take $A_j=A_{j,j}$ and $\Lambda_j=\alpha_{j,j}\underline{\Lambda}_{j,j}$ for all $j\in\N$, with the notations of Section~5.2.1 of~\cite{hnrAM}).\par
For each $j\in\N$, let $u_j$ be the solution of
\be\label{defuk}\left\{\baa{rcll}
-\hbox{div}(A_j\nabla u_j)(x) & = & -H_j(x) & \hbox{in }\Omega,\vspace{3pt}\\
u_j & = & 0 & \hbox{on }\partial\Omega.\eaa\right.
\ee
Each function $u_j$ belongs to $W(\Omega)\cap H^1_0(\Omega)$ and is analytic in $\Omega$. Furthermore, from the previous definitions and from standard elliptic estimates, the functions $u_j$ converge, up to extraction of a subsequence, in $W^{2,p}(\Omega)$ weakly  for all $1\le p<+\infty$ and in $C^{1,\alpha}(\overline{\Omega})$ strongly for all $0\le\alpha<1$ to the solution $u_{\infty}\in W(\Omega)\cap H^1_0(\Omega)$ of
$$\left\{\baa{rcll}
-\hbox{div}(A\nabla u_{\infty})(x) & = & -H_{\infty}(x)=-H(x,u(x),\nabla u(x)) & \hbox{in }\Omega,\vspace{3pt}\\
u_{\infty} & = & 0 & \hbox{on }\partial\Omega,\eaa\right.$$
which, by uniqueness of the solution of this linear problem with right-hand side $-H_{\infty}$, is necessarily equal to $u$. By uniqueness of the limit, one gets that the whole sequence $(u_j)_{j\in\N}$ converges to $u$ in $W^{2,p}(\Omega)$ weak for all $1\le p<+\infty$ and in $C^{1,\alpha}(\overline{\Omega})$ strong for all $0\le\alpha<1$.\par
On the other hand, there are some positive real numbers $N_u$ and $\beta_u$, which depend on~$u$, such that $u\in E_{1/2,N_u,\beta_u}(\Omega)$, because of~(\ref{hypu}) and the smoothness of $u$. Thus, it follows from the convergence of the sequence $(u_j)_{j\in\N}$ to $u$, in (at least) $C^{1,1/2}(\overline{\Omega})$ that, for all~$j$ large enough,
\be\label{ujNbeta}
|\nabla u_j|\neq 0\hbox{ on }\partial\Omega,\ \ u_j>0\hbox{ in }\Omega\ \hbox{ and }\ u_j\in E_{1/2,2N_u,\beta_u/2}(\Omega)
\ee
(notice that the properties $|\nabla u_j|\neq 0$ on $\partial\Omega$ and $u_j>0$ in $\Omega$ are actually automatically fulfilled when the third property $u_j\in E_{1/2,2N_u,\beta_u/2}(\Omega)$ is fulfilled, since $\beta_u>0$, but we prefer to write the three properties all together for the sake of clarity). We can assume that~(\ref{ujNbeta}) holds for all $j\in\N$ without loss of generality.\par
Furthermore, under assumption~(\ref{hypAM}) of Theorem~\ref{th1bis}, it follows from Lemma~\ref{lembounds} that $u\in E_{1/2,N,\beta}(\Omega)$ for some positive constants $N$ and $\beta$ which only depend on $\Omega$, $n$ and $M$ (and which do not depend on~$u$). As in the previous paragraph, one can then assume without loss of generality that, under assumption~(\ref{hypAM}),
\be\label{ujNbetaunif}
|\nabla u_j|\neq 0\hbox{ on }\partial\Omega,\ \ u_j>0\hbox{ in }\Omega\ \hbox{ and }\ u_j\in E_{1/2,2N,\beta/2}(\Omega)
\ee
for all $j\in\N$.

\subsubsection*{Step 3: symmetrized coefficients and the inequalities $u^*_{j_k}\le\hat{\psi}_k$, $(1+\eta_u)\,u^*_{j_k}\le\hat{\psi}_k$ and~$(1+\eta)\,u^*_{j_k}\le\hat{\psi}_k$ in~$\Omega^*$}

Let now $k\in\N$ be fixed in this step and in the next two ones. For all $j\in\N$ and $x\in\Omega$, denote
$$\baa{rcl}
B_j(x) & = & \displaystyle-\,\hbox{div}(A_j\nabla u_j)(x)-a(x,u(x),\nabla u(x))\,|\nabla u_j(x)|+b(x,u(x),\nabla u(x))\,u_j(x)\vspace{3pt}\\
& & \displaystyle-f(x,u(x),\nabla u(x))-2^{-k}\vspace{3pt}\\
& = & -H_j(x)-a(x,u(x),\nabla u(x))\,|\nabla u_j(x)|+b(x,u(x),\nabla u(x))\,u_j(x)\vspace{3pt}\\
& & \displaystyle-f(x,u(x),\nabla u(x))-2^{-k}.\eaa$$
Due to~(\ref{hypH}), (\ref{Hkinfty}) and the fact that $u_j\to u$ in (at least) $C^1(\overline{\Omega})$ as $j\to+\infty$, it follows that
$$\limsup_{j\to+\infty}\Big(\sup_{x\in\Omega}B_j(x)\Big)\le-2^{-k}<0.$$
Therefore, there is an integer $j_k\in\N$ such that $B_{j_k}(x)\le0$ for all $x\in\Omega$, that is
\be\label{Bkeps}\baa{l}
\displaystyle-\,\hbox{div}(A_{j_k}\nabla u_{j_k})(x)-a(x,u(x),\nabla u(x))\,|\nabla u_{j_k}(x)|+b(x,u(x),\nabla u(x))\,u_{j_k}(x)\vspace{3pt}\\
\displaystyle-f(x,u(x),\nabla u(x))-2^{-k}\ \le\ 0\ \hbox{ for all }x\in\Omega.\eaa
\ee
Without loss of generality, one can assume that
$$j_k\ge k.$$\par
One can then apply the general results of Section~\ref{rearrangement} to the coefficients
\be\label{defcoef}\left\{\baa{l}
A_{\Omega}(x)=A_{j_k}(x),\ \Lambda_{\Omega}(x)=\Lambda_{j_k}(x),\ \psi(x)=u_{j_k}(x),\vspace{3pt}\\
a_{\Omega}(x)=a(x,u(x),\nabla u(x)),\ b_{\Omega}(x)=b(x,u(x),\nabla u(x)),\vspace{3pt}\\
f_{\Omega}(x)=f(x,u(x),\nabla u(x))=f_u(x),\eaa\right.
\ee
and $q=1$. Call $\rho_k:[0,\max_{\overline{\Omega}}u_{j_k}]\to[0,R]$, $E_k$, $\hat{\Lambda}_k$, $\hat{\psi}_k$, $\hat{a}_k$ and $\hat{f}_k$ the symmetrized quantities defined as in~(\ref{defrho}), (\ref{defE}), (\ref{defhatlambda}), (\ref{defhatpsi}), (\ref{defhata}) and~(\ref{defhatf}). In particular, the set $E_k$ can be written as
\be\label{Ek}
E_k=\big\{x\in\R^n;\ |x|\in(0,\rho_k(a^k_{m_k-1}))\cup\cdots\cup(\rho_k(a^k_2),\rho_k(a^k_1))\cup(\rho_k(a^k_1),R]\big\}
\ee
where $0<a^k_1<\cdots<a^k_{m_k}=\max_{\overline{\Omega}}u_{j_k}$ denote the $m_k$ critical values of the function~$u_{j_k}$ in~$\overline{\Omega}$. The function $\hat{\psi}_k$ belongs to $W^{1,\infty}(\Omega^*)\cap H^1_0(\Omega^*)$ and is of class $C^1$ in $E_k\cup\{0\}$ and $C^2$ in $E_k\cap\Omega^*$, the function $\hat{\Lambda}_k\in L^{\infty}_+(\Omega^*)$ is of class $C^1$ in $E_k\cap\Omega^*$ and the functions $\hat{a}_k$ and~$\hat{f}_k\in L^{\infty}(\Omega^*)$ are continuous in $E_k$. All functions $\hat{\Lambda}_k$, $\hat{\psi}_k$, $\hat{a}_k$ and $\hat{f}_k$ are radially symmetric. It follows from~(\ref{ineqlambda}),~(\ref{ineqa}),~(\ref{ineqf}) and~(\ref{Aklambdak}) that
\be\label{coefepsilon}\left\{\baa{l}
\displaystyle0<\min_{\overline{\Omega}}\Lambda_{j_k}\le\mathop{{\rm{ess}}\,{\rm{inf}}}_{\Omega^*}\hat{\Lambda}_k\le\mathop{{\rm{ess}}\,{\rm{sup}}}_{\Omega^*}\hat{\Lambda}_k\le\max_{\overline{\Omega}}\Lambda_{j_k},\vspace{3pt}\\
\|\hat{\Lambda}_k^{-1}\|_{L^1(\Omega^*)}=\|\Lambda_{j_k}^{-1}\|_{L^1(\Omega)}=\|\Lambda^{-1}\|_{L^1(\Omega)},\vspace{3pt}\\
\displaystyle\inf_{\Omega\times\R\times\R^n}a^+\le\min_{\overline{\Omega}}a(\cdot,u(\cdot),\nabla u(\cdot))^+\le\mathop{{\rm{ess}}\,{\rm{inf}}}_{\Omega^*}\hat{a}_k\le\mathop{{\rm{ess}}\,{\rm{sup}}}_{\Omega^*}\hat{a}_k\le\cdots\\
\qquad\qquad\qquad\qquad\qquad\qquad\qquad\qquad\displaystyle\cdots\le\max_{\overline{\Omega}}a(\cdot,u(\cdot),\nabla u(\cdot))^+\le\sup_{\Omega\times\R\times\R^n}a^+,\vspace{3pt}\\
\displaystyle\inf_{\Omega\times\R\times\R^n}f\le\min_{\overline{\Omega}}f_u\le\mathop{{\rm{ess}}\,{\rm{inf}}}_{\Omega^*}\hat{f}_k\le\mathop{{\rm{ess}}\,{\rm{sup}}}_{\Omega^*}\hat{f}_k\le\max_{\overline{\Omega}}f_u\le\sup_{\Omega\times\R\times\R^n}f.\eaa\right.
\ee\par
For the proof of Theorem~\ref{th1}, it follows then from Proposition~\ref{key1} that
\be\label{ineqkey1}
\hat{\psi}_k(x)\ge u_{j_k}(y)\ge 0\ \hbox{ for all }x\in\overline{\Omega^*}\hbox{ and }y\in\Sigma_{\rho_k^{-1}(|x|)}.
\ee
That means that
\be\label{schwarzuk}
0\le u_{j_k}^*(x)\le\hat{\psi}_k(x)\ \hbox{ for all }x\in\overline{\Omega^*},
\ee
where $u_{j_k}^*$ denotes the Schwarz symmetrization of the function $u_{j_k}$.\par
On the other hand, if $\Omega$ is not a ball, it follows from~(\ref{ujNbeta}) and Proposition~\ref{key1notball} that
$$\hat{\psi}_k(x)\ge(1+\eta_u)\,u_{j_k}(y)\ge 0\ \hbox{ for all }x\in\overline{\Omega^*}\hbox{ and }y\in\Sigma_{\rho_k^{-1}(|x|)},$$
where $\eta_u>0$ only depends on $\Omega$, $n$, $N_u$ and $\beta_u$, that is on $\Omega$, $n$ and $u$. Therefore,
\be\label{schwarzukbisu}
0\le(1+\eta_u)\,u_{j_k}^*(x)\le\hat{\psi}_k(x)\ \hbox{ for all }x\in\overline{\Omega^*}.
\ee\par
Furthermore, if $\Omega$ is not a ball and the assumption~(\ref{hypAM}) of Theorem~\ref{th1bis} is made, it follows from~(\ref{ujNbetaunif}) and Proposition~\ref{key1notball} that
$$\hat{\psi}_k(x)\ge(1+\eta)\,u_{j_k}(y)\ge 0\ \hbox{ for all }x\in\overline{\Omega^*}\hbox{ and }y\in\Sigma_{\rho_k^{-1}(|x|)},$$
where $\eta>0$ only depends on $\Omega$, $n$, $N$ and $\beta$, that is on $\Omega$, $n$ and $M$. Therefore,
\be\label{schwarzukbis}
0\le(1+\eta)\,u_{j_k}^*(x)\le\hat{\psi}_k(x)\ \hbox{ for all }x\in\overline{\Omega^*}
\ee
in this case.\par
Lastly, for both Theorems~\ref{th1} and~\ref{th1bis}, Proposition~\ref{key2} implies that, for all $x\in E_k\cap\Omega^*$, there exists a point~$y\in\Sigma_{\rho_k^{-1}(|x|)}$ such that
$$\begin{array}{rl}
& \displaystyle -\,{\rm{div}}\big(\widehat{\Lambda}_k\nabla\hat{\psi}_k\big)(x)-\widehat{a}_k(x)\vert \nabla\hat{\psi}_k(x)\vert-\widehat{f}_k(x)-2^{-k}\vspace{3pt}\\
\leq\!\! & \displaystyle-\,{\rm{div}}(A_{j_k}\nabla u_{j_k})(y)-a(y,u(y),\nabla u(y))\left\vert \nabla u_{j_k}(y)\right\vert+b(y,u(y),\nabla u(y))\,u_{j_k}(y)\vspace{3pt}\\
& -f(y,u(y),\nabla u(y))-2^{-k}\vspace{3pt}\\
\leq\!\! & 0,\end{array}$$
where the last inequality follows from~(\ref{Bkeps}). In other words,
\be\label{ineqpsieps}
-\,{\rm{div}}\big(\widehat{\Lambda}_k\nabla\hat{\psi}_k\big)(x)+\hat{a}_k(x)e_r(x)\cdot\nabla\hat{\psi}_k(x)\le g_k(x)\ \hbox{ for all }x\in E_k\cap\Omega^*,
\ee
where $e_r(x)=x/|x|$ for all $x\in\R^n\backslash\{0\}$ and
\be\label{defalphak}
g_k(x)=\widehat{f}_k(x)+2^{-k}
\ee
for all $x\in E_k$ (remember indeed that $\nabla\hat{\psi}_k(x)$ points in the direction of $-e_r(x)$ for all $x\in E_k$, from~(\ref{eqF}) and~(\ref{defhatpsi})).

\subsubsection*{Step 4: the functions $\hat{\psi}_k$ are $H^1_0(\Omega^*)$ weak subsolutions of~$(\ref{ineqpsieps})$}

The inequality~(\ref{ineqpsieps}) holds in $E_k\cap\Omega^*$, whence almost everywhere in $\Omega^*$. But the quantities appearing in~(\ref{ineqpsieps}) might be discontinuous across the critical spheres $\partial E_k$ in general. The goal of this step is to show that~(\ref{ineqpsieps}) holds nevertheless in the~$H^1_0(\Omega^*)$ weak sense as well, as stated in the following lemma.

\begin{lem}\label{lemsubpsik}
There holds
\be\label{subpsieps}
\int_{\Omega^*}\hat{\Lambda}_k\nabla\hat{\psi}_k\cdot\nabla\varphi+\int_{\Omega^*}\big(\hat{a}_ke_r\cdot\nabla\hat{\psi}_k\big)\varphi-\int_{\Omega^*}g_k\varphi\le0
\ee
for all $k\in\N$ and for all $\varphi\in H^1_0(\Omega^*)$ with $\varphi\ge0$ a.e. in $\Omega^*$.
\end{lem}

In order not to lengthen the main scheme of the proofs of Theorems~\ref{th1} and~\ref{th1bis}, the proof of Lemma~\ref{lemsubpsik} is postponed in Section~\ref{sec32} below.

\subsubsection*{Step 5: the inequalities $u^*_{j_k}\le v_k$, $(1+\eta_u)\,u^*_{j_k}\le v_k$ and $(1+\eta)\,u^*_{j_k}\le v_k$ in $\Omega^*$}

We first point out that $\hat{a}_k$ and~$g_k$ are in $L^{\infty}(\Omega^*)$. Let then $v_k$ be the unique $H^1_0(\Omega^*)$ solution of
\be\label{defvk}\left\{\baa{rcll}
-\,{\rm{div}}\big(\widehat{\Lambda}_k\nabla v_k\big)+\hat{a}_k\,e_r\cdot\nabla v_k & = & g_k & \hbox{in }\Omega^*,\vspace{3pt}\\
v_k & = & 0 & \hbox{on }\partial\Omega^*,\eaa\right.
\ee
where the above equation is understood in the weak sense, that is
$$\int_{\Omega^*}\hat{\Lambda}_k\,\nabla v_k\cdot\nabla\varphi+\int_{\Omega^*}\hat{a}_k\,(e_r\cdot\nabla v_k)\,\varphi-\int_{\Omega^*}g_k\,\varphi=0$$
for every $\varphi\in H^1_0(\Omega^*)$. The existence and uniqueness of~$v_k$ is guaranteed by Theorem~8.3 of~\cite{gt}. Since $\hat{\psi}_k$ is an $H^1_0(\Omega^*)$ subsolution of this problem, in the sense of Lemma~\ref{lemsubpsik} of Step~4, it then follows from the weak maximum principle (see Theorem~8.1 of~\cite{gt}) that
$$\hat{\psi}_k\le v_k\ \hbox{ a.e. in }\Omega^*.$$
Hence,~(\ref{schwarzuk}) yields
\be\label{inequjkvk}
0\le u_{j_k}^*\le v_k\ \hbox{ a.e. in }\Omega^*
\ee
under the assumptions of Theorem~\ref{th1}, whereas~(\ref{schwarzukbisu}) implies that
\be\label{inequjkvkbisu}
0\le(1+\eta_u)\,u_{j_k}^*\le v_k\ \hbox{ a.e. in }\Omega^*
\ee
if $\Omega$ is not a ball, and (\ref{schwarzukbis}) yields
\be\label{inequjkvkbis}
0\le(1+\eta)\,u_{j_k}^*\le v_k\ \hbox{ a.e. in }\Omega^*,
\ee
if $\Omega$ is not a ball and~(\ref{hypAM}) is assumed, where $\eta_u>0$ and $\eta>0$ are as in Step~3.

\subsubsection*{Step 6: the limiting inequalities $u^*\le v_{\infty}$, $(1+\eta_u)\,u^*\le v_{\infty}$ and $(1+\eta)\,u^*\le v_{\infty}$ in $\Omega^*$}

First of all, since $u_j\to u$ as $j\to+\infty$ in (at least)~$C^1(\overline{\Omega})$ and since $j_k\ge k$ for all $k\in\N$, it follows from~\cite{ct} that $u_{j_k}^*\to u^*$ in $L^1(\Omega^*)$ as~$k\to+\infty$. Up to extraction of a subsequence, one can then assume that
\be\label{ujk*}
u_{j_k}^*(x)\to u^*(x)\hbox{ a.e. in }\Omega^*\hbox{ as }k\to+\infty.
\ee\par
Let us now pass to the limit in the $H^1_0(\Omega^*)$ solutions $v_k$ of~(\ref{defvk}). Notice first, from~(\ref{Aklambdak}),~(\ref{coefepsilon}) and~(\ref{defalphak}), that the sequences $(\hat{\Lambda}_k)_{k\in\N}$, $(\hat{\Lambda}_k^{-1})_{k\in\N}$, $(\hat{a}_k)_{k\in\N}$ and $(g_k)_{k\in\N}$ are bounded in $L^{\infty}(\Omega^*)$. It follows then from Corollary~8.7 of~\cite{gt} that the sequence $(v_k)_{k\in\N}$ is bounded in $H^1_0(\Omega^*)$. Therefore, up to extraction of a subsequence, there exists a radially symmetric function $v_{\infty}\in H^1_0(\Omega^*)$ such that
\be\label{vkv}
v_k\rightharpoonup v_{\infty}\hbox{ in }H^1_0(\Omega^*)\hbox{ weak},\ v_k\to v_{\infty}\hbox{ in }L^2(\Omega^*)\hbox{ strong and a.e. in }\Omega^*\hbox{ as }k\to+\infty.
\ee\par
Together with~(\ref{inequjkvk}),~(\ref{inequjkvkbisu}),~(\ref{inequjkvkbis}) and~(\ref{ujk*}), one gets that
\be\label{u*vinfty}
0\le u^*\le v_{\infty}\ \hbox{ a.e. in }\Omega^*
\ee
under the assumptions of Theorem~\ref{th1}, that
\be\label{u*vinftybisu}
0\le(1+\eta_u)\,u^*\le v_{\infty}\ \hbox{ a.e. in }\Omega^*
\ee
if $\Omega$ is not a ball, and that
\be\label{u*vinftybis}
0\le(1+\eta)\,u^*\le v_{\infty}\ \hbox{ a.e. in }\Omega^*
\ee
if $\Omega$ is not a ball and assumption~(\ref{hypAM}) is made, where $\eta_u>0$ and $\eta>0$ are as in Step~3.

\subsubsection*{Step 7: a limiting equation satisfied by~$v_{\infty}$ in $\Omega^*$}

Let us now pass to the limit in the coefficients $\hat{\Lambda}_k$, $\hat{a}_k$ and $g_k$ of~(\ref{defvk}). From~(\ref{Aklambdak}),~(\ref{coefepsilon}) and~(\ref{defalphak}), there exist three radially symmetric functions $\hat{\Lambda}\in L^{\infty}_+(\Omega^*)$, $\hat{a}\in L^{\infty}(\Omega^*)$ and $\hat{f}\in L^{\infty}(\Omega^*)$ such that, up to extraction of some subsequence,
\be\label{chapeaux}
\hat{\Lambda}_k^{-1}\rightharpoonup\hat{\Lambda}^{-1},\ \ \hat{\Lambda}_k^{-1}\hat{a}_k\rightharpoonup\hat{\Lambda}^{-1}\hat{a}\ \hbox{ and }\ g_k\rightharpoonup\hat{f}\ \hbox{ in }L^{\infty}(\Omega^*)\hbox{ weak-* }\hbox{ as }k\to+\infty,
\ee
whence
\be\label{hats2}
\displaystyle0<\mathop{{\rm{ess}}\,{\rm{inf}}}_{\Omega}\Lambda\le\mathop{{\rm{ess}}\,{\rm{inf}}}_{\Omega^*}\hat{\Lambda}\le\mathop{{\rm{ess}}\,{\rm{sup}}}_{\Omega^*}\hat{\Lambda}\le\mathop{{\rm{ess}}\,{\rm{sup}}}_{\Omega}\Lambda\ \hbox{ and }\ \|\hat{\Lambda}^{-1}\|_{L^1(\Omega^*)}=\|\Lambda^{-1}\|_{L^1(\Omega)}.
\ee
Namely, the function $\hat{a}=\hat{\Lambda}\,\hat{\Lambda}^{-1}\,\hat{a}$ is defined as $\hat{\Lambda}$ times the $L^{\infty}(\Omega^*)$ weak-* limit of the sequence~$(\hat{\Lambda}_k^{-1}\hat{a}_k)_{k\in\N}$. Furthermore,~$\hat{a}$ is thus the $L^{\infty}(\Omega^*)$ weak-* limit of the functions $\hat{\Lambda}\,\hat{\Lambda}_k^{-1}\,\hat{a}_k$. Since
$$\min_{\overline{\Omega}}a(\cdot,u(\cdot),\nabla u(\cdot))^+\le\mathop{{\rm{ess}}\,{\rm{inf}}}_{\Omega^*}\hat{a}_k\le\mathop{{\rm{ess}}\,{\rm{sup}}}_{\Omega^*}\hat{a}_k\le\max_{\overline{\Omega}}a(\cdot,u(\cdot),\nabla u(\cdot))^+$$
from~(\ref{coefepsilon}), while $(0<)\,\hat{\Lambda}\,\hat{\Lambda}_k^{-1}\rightharpoonup 1$ in the $L^{\infty}(\Omega^*)$ weak-* sense as $k\to+\infty$, it follows that
\be\label{hats3}\baa{l}
\displaystyle0\le\inf_{\Omega\times\R\times\R^n}\!a^+\!\le\min_{\overline{\Omega}}a(\cdot,u(\cdot),\nabla u(\cdot))^+\!\le\mathop{{\rm{ess}}\,{\rm{inf}}}_{\Omega^*}\hat{a}\le\mathop{{\rm{ess}}\,{\rm{sup}}}_{\Omega^*}\hat{a}\le\cdots\vspace{3pt}\\
\qquad\qquad\qquad\qquad\qquad\qquad\qquad\cdots\le\max_{\overline{\Omega}}a(\cdot,u(\cdot),\nabla u(\cdot))^+\!\le\sup_{\Omega\times\R\times\R^n}\!a^+\!.\eaa
\ee\par
The main goal of this step is to show that $v_{\infty}$ is a weak $H^1_0(\Omega^*)$ solution of the limiting equation obtained by passing formally to the limit as $k\to+\infty$ in~(\ref{defvk}).

\begin{lem}\label{vinfty}
The function $v_{\infty}$ is a weak $H^1_0(\Omega^*)$ solution of
\be\label{eqvlimit}\left\{\baa{rcll}
-\,{\rm{div}}\big(\widehat{\Lambda}\nabla v_{\infty}\big)+\hat{a}\,e_r\cdot\nabla v_{\infty} & = & \hat{f} & \hbox{in }\Omega^*,\vspace{3pt}\\
v_{\infty} & = & 0 & \hbox{on }\partial\Omega^*,\eaa\right.
\ee
in the sense that
$$\int_{\Omega^*}\hat{\Lambda}\,\nabla v_{\infty}\cdot\nabla\varphi+\int_{\Omega^*}(\hat{a}\,e_r\cdot\nabla v_{\infty})\,\varphi-\int_{\Omega^*}\hat{f}\,\varphi=0$$
for every $\varphi\in H^1_0(\Omega^*)$.
\end{lem}

In order to go on the last steps of the proofs of Theorems~\ref{th1} and~\ref{th1bis}, the proof of Lemma~\ref{vinfty} is postponed in Section~\ref{sec32}.\par
The radially symmetric functions $\hat{\Lambda}\in L^{\infty}_+(\Omega^*)$ and $\hat{a}\in L^{\infty}(\Omega^*)$ will be those of the statements of Theorems~\ref{th1} and~\ref{th1bis}. Notice in particular that the properties~(\ref{hats}) follow from~(\ref{hats2}) and~(\ref{hats3}). Furthermore, we already know from~(\ref{u*vinfty}),~(\ref{u*vinftybisu}) and~(\ref{u*vinftybis}) that $0\le u^* \le v_{\infty}$ a.e. in $\Omega^*$ (respectively~$0\le(1+\eta_u)\,u^*\le v_{\infty}$ and $0\le(1+\eta)\,u^*\le v_{\infty}$) under the assumptions of Theorem~\ref{th1} (respectively Theorem~\ref{th1bis} when $\Omega$ is not a ball, without or with assumption~(\ref{hypAM})), where $v_{\infty}$ solves~(\ref{eqvlimit}). However, the right-hand side of~(\ref{eqvlimit}) involves a function $\hat{f}$ which may not be the Schwarz rearrangement $f_u^*$ of the function $f_u$ defined in~(\ref{deffu}), or may even not have the same distribution function as $f_u$. In the remaining four steps of the proofs of Theorems~\ref{th1} and~\ref{th1bis}, one shall then approximate the function~$v_{\infty}$ by some functions~$w_k$ and~$z_l$ (different from the~$v_k$'s in general) and the functions $\hat{f}_k$ by some functions having the same distribution function as~$f_u=f(\cdot,u(\cdot),\nabla u(\cdot))$, before finally comparing a function $z_L$ for $L$ large enough to the solution~$v$ of~(\ref{eqvbis}) with $H$ and $f_u^*$.

\subsubsection*{Step 8: approximation of $v_{\infty}$ by some functions $w_k$ in $\Omega^*$}

Let $(\hat{f}_k)_{k\in\N}$ be the sequence of radially symmetric functions defined in Step~3, and remember that the sequence $(\hat{f}_k)_{k\in\N}$ is bounded in $L^{\infty}(\Omega^*)$ from~(\ref{coefepsilon}) (one could replace $\hat{f}_k$ by $g_k$ without any change in the conclusions). Since $\hat{\Lambda}\in L^{\infty}_+(\Omega^*)$ and $\hat{a}\in L^{\infty}(\Omega^*)$, it follows from Theorem~8.3 of~\cite{gt} that, for each $k\in\N$, there is a unique weak $H^1_0(\Omega^*)$ solution~$w_k$ of
\be\label{eqwk}\left\{\baa{rcll}
-\,{\rm{div}}\big(\widehat{\Lambda}\nabla w_k\big)+\hat{a}\,e_r\cdot\nabla w_k & = & \hat{f}_k & \hbox{in }\Omega^*,\vspace{3pt}\\
w_k & = & 0 & \hbox{on }\partial\Omega^*,\eaa\right.
\ee
in the sense that
\be\label{weakwk}
\int_{\Omega^*}\hat{\Lambda}\,\nabla w_k\cdot\nabla\varphi+\int_{\Omega^*}(\hat{a}\,e_r\cdot\nabla w_k)\,\varphi-\int_{\Omega^*}\hat{f}_k\,\varphi=0\ \hbox{ for all }\varphi\in H^1_0(\Omega^*).
\ee
Furthermore, the functions $w_k$ are all radially symmetric by uniqueness and since all the coefficients $\hat{\Lambda}$, $\hat{a}$ and $\hat{f}_k$ are radially symmetric. Lastly, as in Step~6 above, the sequence $(w_k)_{k\in\N}$ is bounded in $H^1_0(\Omega^*)$ from Corollary~8.7 of~\cite{gt}. There exists then a function $w_{\infty}\in H^1_0(\Omega^*)$ such that, up to extraction of a subsequence, one has
$$w_k\rightharpoonup w_{\infty}\hbox{ in }H^1_0(\Omega^*)\hbox{ weak and }w_k\to w_{\infty}\hbox{ in }L^2(\Omega^*)\hbox{ strong as }k\to+\infty.$$
Since $\hat{f}_k=g_k-2^{-k}\rightharpoonup\hat{f}$ in $L^{\infty}(\Omega^*)$ weak-* as $k\to+\infty$ from~(\ref{defalphak}) and~(\ref{chapeaux}), it follows by passing to the limit as $k\to+\infty$ in~(\ref{weakwk}) that
$$\int_{\Omega^*}\hat{\Lambda}\,\nabla w_{\infty}\cdot\nabla\varphi+\int_{\Omega^*}(\hat{a}\,e_r\cdot\nabla w_{\infty})\,\varphi-\int_{\Omega^*}\hat{f}\,\varphi=0\ \hbox{ for all }\varphi\in H^1_0(\Omega^*).$$
In other words, $w_{\infty}$ is a weak $H^1_0(\Omega^*)$ solution of~(\ref{eqvlimit}). Referring again to Theorem~8.3 of~\cite{gt} for the uniqueness of the solution of~(\ref{eqvlimit}), one concludes that $w_{\infty}=v_{\infty}$ and that, by uniqueness of the limit, the whole sequence $(w_k)_{k\in\N}$ converges to $v_{\infty}$ in the following sense:
\be\label{wkvinfty}
w_k\rightharpoonup v_{\infty}\hbox{ in }H^1_0(\Omega^*)\hbox{ weak and }w_k\to v_{\infty}\hbox{ in }L^2(\Omega^*)\hbox{ strong as }k\to+\infty.
\ee

\subsubsection*{Step 9: approximation of the function $w_K$ for $K$ large by some functions $z_l$ in $\Omega^*$}

Let $\epsilon>0$ be arbitrary and given until the end of the proofs of Theorems~\ref{th1} and~\ref{th1bis}. From~(\ref{wkvinfty}), there is an integer $K\in\N$ large enough such that
\be\label{wKv}
\|w_K-v_{\infty}\|_{L^2(\Omega^*)}\le\frac{\epsilon}{2}.
\ee
The function $w_K$ is the weak $H^1_0(\Omega^*)$ solution of~(\ref{eqwk}) with $k=K$. The radially symmetric function $\hat{f}_K\in L^{\infty}(\Omega^*)$ in the right-hand side of~(\ref{eqwk}) is given by~(\ref{defhatf}) with
$$f_{\Omega}(y)=f(y,u(y),\nabla u(y))=f_u(y)$$
and $\rho=\rho_K$ is given by the function $\psi=u_{j_K}$: namely, with the general notations of Section~\ref{rearrangement}, there holds
$$\widehat{f}_K(x)=\frac{\displaystyle{\int_{\Sigma_{\rho_K^{-1}(|x|)}} f_u(y)\,\vert\nabla u_{j_K}(y)\vert^{-1}d\sigma_{\rho_K^{-1}(|x|)}}}{\displaystyle{\int_{\Sigma_{\rho_K^{-1}(|x|)}} \left\vert \nabla u_{j_K}(y)\right\vert^{-1}d\sigma_{\rho_K^{-1}(|x|)}}}$$
for all $x\in E_K$, where the set $E_K\subset\overline{\Omega^*}$ is given by~(\ref{Ek}) with $k=K$. As already noticed in the general properties of Section~\ref{rearrangement}, one knows that $\big|\big\{y\in\Omega;\ u_{j_K}(y)=a\big\}\big|=0$ for every~$a\in[0,\max_{\overline{\Omega}}u_{j_K}]$. It follows then from the co-area formula that
$$\int_{S_{\rho_K(b),\rho_K(a)}}\hat{f}_K=\int_{\Omega_{a,b}}f_u$$
for all $0\le a<b\le\max_{\overline{\Omega}}u_{j_K}$, where $\Omega_{a,b}=\big\{y\in\Omega;\ a<u_{j_K}(y)<b\big\}$. Therefore, one infers from Lemma~5.1 of~\cite{hnrAM} (see also Lemma 1.1 of~\cite{at3} and Lemma 2.2 of~\cite{at0}) that there is a sequence $(h_l)_{l\in\N}$ of $L^{\infty}(\Omega^*)$ radially symmetric functions such that
$$h_l\rightharpoonup\hat{f}_K\hbox{ as }l\to+\infty\hbox{ in }L^{\infty}(\Omega^*)\hbox{ weak-* and }\mu_{h_l}=\mu_{f_u}\hbox{ for all }l\in\N,$$
that is the functions $h_l$ have all the same distribution function as the function $f_u$.\par
On the other hand, as in Step~8, there is for every $l\in\N$ a unique weak $H^1_0(\Omega^*)$, and radially symmetric, solution $z_l$ of
\be\label{eqzl}\left\{\baa{rcll}
-\,{\rm{div}}\big(\widehat{\Lambda}\nabla z_l\big)+\hat{a}\,e_r\cdot\nabla z_l & = & h_l & \hbox{in }\Omega^*,\vspace{3pt}\\
z_l & = & 0 & \hbox{on }\partial\Omega^*\eaa\right.
\ee
and the functions $z_l$ converge to $w_K$ in $H^1_0(\Omega^*)$ weak and in $L^2(\Omega^*)$ strong as $l\to+\infty$. In particular, there is $L\in\N$ large enough such that $\|z_L-w_K\|_{L^2(\Omega^*)}\le\epsilon/2$, whence
\be\label{zLvinfty}
\|z_L-v_{\infty}\|_{L^2(\Omega^*)}\le\epsilon
\ee
from~(\ref{wKv}).

\subsubsection*{Step 10: the inequality $z_L\le z$ in $\Omega^*$}

Remember that the distributions functions $\mu_{h_L}$ and $\mu_{f_u}$ of the functions $h_L\in L^{\infty}(\Omega^*)$ and $f_u\in L^{\infty}(\Omega)$ are identical, and let $z$ be the unique weak $H^1_0(\Omega^*)$, and radially symmetric, solution of
\be\label{eqz}\left\{\baa{rcll}
-\,{\rm{div}}\big(\widehat{\Lambda}\nabla z\big)+\hat{a}\,e_r\cdot\nabla z & = & f_u^* & \hbox{in }\Omega^*,\vspace{3pt}\\
z & = & 0 & \hbox{on }\partial\Omega^*,\eaa\right.
\ee
where $f^*_u\in L^{\infty}(\Omega^*)$ is the radially symmetric Schwarz rearrangement of the function $f_u$. The key-point of this step is the following lemma.

\begin{lem}\label{lemzLz}
There holds
\be\label{zLz}
z_L\le z\hbox{ a.e. in }\Omega^*.
\ee
\end{lem}

Actually, the same inequality holds with $z_l$ for every $l\in\N$, but we will only use it with the function~$z_L$. Notice also that, from the radial symmetry and the De Giorgi-Moser-Nash regularity theory (see Theorem 8.29 in~\cite{gt}), both functions $z_L$ and $z$ can be assumed to be continuous in $\overline{\Omega^*}$, even if it means redefining them on the negligible subset of~$\overline{\Omega^*}$. Therefore, the inequality $z_L(x)\le z(x)$ holds for all $x\in\overline{\Omega^*}$ without loss of generality.\par
The proof of this lemma is postponed in Section~\ref{sec32}. Let us now conclude the proofs of Theorems~\ref{th1} and~\ref{th1bis}.

\subsubsection*{Step 11: conclusion of the proofs of Theorems~\ref{th1} and~\ref{th1bis}}

Let $v$ be the unique weak $H^1_0(\Omega^*)$ solution of the equation~(\ref{eqvbis}) of Theorems~\ref{th1} and~\ref{th1bis}, that is
$$\left\{\baa{rcll}
-\,{\rm{div}}\big(\widehat{\Lambda}\nabla v\big)-\hat{a}\,|\nabla v| & = & f_u^* & \hbox{in }\Omega^*,\vspace{3pt}\\
v & = & 0 & \hbox{on }\partial\Omega^*.\eaa\right.$$
Remember that the existence and uniqueness of $v$ is guaranteed by Theorem~2.1 of~\cite{Porretta}. We also recall that $\hat{a}\ge0$ a.e. in $\Omega^*$ from~(\ref{hats3}). In particular,
$$-\hbox{div}(\hat{\Lambda}\nabla v)+\hat{a}\,e_r\cdot\nabla v\ge-\hbox{div}(\hat{\Lambda}\nabla v)-\hat{a}\,|\nabla v|=f^*_u$$
in the weak $H^1_0(\Omega^*)$ sense, that is
$$\int_{\Omega^*}\hat{\Lambda}\,\nabla v\cdot\nabla\varphi+\int_{\Omega^*}(\hat{a}\,e_r\cdot\nabla v)\,\varphi-\int_{\Omega^*}f^*_u\,\varphi\ge0$$
for every $\varphi\in H^1_0(\Omega^*)$ with $\varphi\ge0$ a.e. in $\Omega^*$. In other words, the function $v$ is a weak $H^1_0(\Omega^*)$ supersolution of the equation~(\ref{eqz}) satisfied by~$z$. The maximum principle (Theorem~8.1 of~\cite{gt}) implies that
$$z\le v\hbox{ a.e. in }\Omega^*.$$
Together with~(\ref{zLz}), one gets that
\be\label{zLv}
z_L\le v\hbox{ a.e. in }\Omega^*.
\ee\par
As a conclusion, it follows from~(\ref{u*vinfty}),~(\ref{u*vinftybisu}),~(\ref{u*vinftybis}),~(\ref{zLvinfty}) and~(\ref{zLv}) that, under the assumptions of Theorem~\ref{th1},
\be\label{uastv}\baa{rcl}
\|(u^*-v)^+\|_{L^2(\Omega^*)} & \le & \|(u^*-v_{\infty})^++(v_{\infty}-z_L)^++(z_L-v)^+\|_{L^2(\Omega^*)}\vspace{3pt}\\
& \le & \|(u^*\!-\!v_{\infty})^+\|_{L^2(\Omega^*)}+\|(v_{\infty}\!-\!z_L)^+\|_{L^2(\Omega^*)}+\|(z_L\!-\!v)^+\|_{L^2(\Omega^*)}\vspace{3pt}\\
& \le & 0+\epsilon+0=\epsilon,\eaa
\ee
whereas, under the assumptions of Theorem~\ref{th1bis}, there holds similarly
$$\|((1+\eta_u)u^*-v)^+\|_{L^2(\Omega^*)}\le\epsilon$$
if $\Omega$ is not a ball and
$$\|((1+\eta)u^*-v)^+\|_{L^2(\Omega^*)}\le\epsilon$$
if $\Omega$ is not a ball and assumption~(\ref{hypAM}) is made, where $\eta_u>0$ (resp. $\eta>0$) only depends on~$\Omega$,~$n$ and~$u$ (resp. on $\Omega$, $n$ and the constant $M$ in~(\ref{hypAM})).\par
Since $\epsilon>0$ can be arbitrarily small and $u^*$ and $v$ do not depend on $\epsilon$, one concludes that~$\|(u^*-v)^+\|_{L^2(\Omega^*)}=0$ under the assumptions of Theorem~\ref{th1}, that is
\be\label{u*v2}
u^*\le v\hbox{ a.e. in }\Omega^*,
\ee
whereas $\|((1+\eta_u)u^*-v)^+\|_{L^2(\Omega^*)}=0$ (resp. $\|((1+\eta)u^*-v)^+\|_{L^2(\Omega^*)}=0$), that is
$$(1+\eta_u)\,u^*\le v\hbox{ a.e. in }\Omega^*$$
(resp. $(1+\eta)\,u^*\le v$ a.e. in $\Omega^*$), if $\Omega$ is not a ball (resp. if $\Omega$ is not a ball and assumption~(\ref{hypAM}) is made). They are the desired conclusions. The proofs of Theorems~\ref{th1} and~\ref{th1bis} are thereby complete.\hfill\fin

\begin{rem}\label{remu*z}{\rm By replacing $v$ by $z$ in~\eqref{uastv} and by using directly~\eqref{zLz} instead of~\eqref{zLv}, it also follows that
\be\label{u*z}
u^*\le z\hbox{ a.e. in }\Omega^*
\ee
under the assumptions of Theorem~$\ref{th1}$, whereas $(1+\eta_u)\,u^*\le z$ a.e. in $\Omega^*$ (resp. $(1+\eta)\,u^*\le z$ a.e. in $\Omega^*$) if $\Omega$ is not a ball (resp. if $\Omega$ is not a ball and assumption~(\ref{hypAM}) is made).}
\end{rem}


\subsection{Proofs of Lemmas~\ref{lembounds},~\ref{lemsubpsik},~\ref{vinfty} and~\ref{lemzLz}}\label{sec32}

This section is devoted to the proof of four technical lemmas which have been used for the proofs of Theorems~\ref{th1} and~\ref{th1bis} in the previous section.\hfill\break

\noindent{\bf{Proof of Lemma~\ref{lembounds}.}} We first prove uniform $L^{\infty}(\Omega)$ bounds for the function $u$, from which uniform $C^{1,1/2}(\overline{\Omega})$ estimates will follow. Next, we prove a uniform lower bound on $u$. These estimates will imply that $u\in E_{1/2,N,\beta}(\Omega)$ for some constants positive $N$ and $\beta$ which do not depend on $u$. In the proof, we denote $C_i$ some constants which may depend on $\Omega$, $n$ and~$M>0$ given in~(\ref{hypAM}), but which do not depend on the given solution~$u$ of~(\ref{equ}). We recall that~(\ref{ALambda}),~(\ref{hypH}) with $q=1$ and~(\ref{hypu}) are assumed throughout the proof.\par
First of all, it follows from~(\ref{hypH}) and~(\ref{hypu}) that
$$-\hbox{div}(A(x)\nabla u)+q_u(x)\cdot\nabla u\le f_u(x)$$
a.e. in $\Omega$ (we recall that $A\in W^{1,\infty}(\Omega,{\mathcal{S}}_n(\R))$ and $u\in W(\Omega)$), where $f_u$ is given in~(\ref{deffu}) and
$$q_u(x)=\left\{\baa{ll}
\displaystyle -a(x,u(x),\nabla u(x))\frac{\nabla u(x)}{|\nabla u(x)|} & \hbox{if }|\nabla u(x)|\neq 0,\vspace{3pt}\\
0 & \hbox{if }|\nabla u(x)|=0.\eaa\right.$$
It follows then from the maximum principle (Theorem~8.1 in~\cite{gt}) that
$$u\le U\ \hbox{ a.e. in }\Omega,$$
where $U\in H^1_0(\Omega)\cap W(\Omega)$ denotes the solution of
$$\left\{\baa{rcll}
-\div(A(x)\nabla U)+q_u(x)\cdot\nabla U & = & f_u(x) & \hbox{in }\Omega,\vspace{3pt}\\
U & = & 0 & \hbox{on }\partial\Omega,\eaa\right.$$
Since $\|q_u\|_{L^{\infty}(\Omega,\R^n)}\le\|a\|_{L^{\infty}(\Omega\times\R\times\R^n)}\le M$, $\|f_u\|_{L^{\infty}(\Omega)}\le\|f\|_{L^{\infty}(\Omega\times\R\times\R^n)}\le M$ and $\|A\|_{W^{1,\infty}(\Omega)}+\|\Lambda^{-1}\|_{L^{\infty}(\Omega)}\le M$, it follows from standard elliptic estimates that $\|U\|_{L^{\infty}(\Omega)}\le C_1$ for some positive constant $C_1$ (which depends on $\Omega$, $n$ and $M$ but not on $u$). Since~$0\le u\le U$, one concludes that
\be\label{C1}
\|u\|_{L^{\infty}(\Omega)}\le C_1.
\ee\par
By testing~(\ref{equ}) against $u$ itself and using~(\ref{ALambda}) and~(\ref{hypAM}), one gets that
$$\baa{rcl}
\displaystyle M^{-1}\int_{\Omega}|\nabla u(x)|^2dx & \le & \displaystyle\int_{\Omega}\Lambda(x)\,|\nabla u(x)|^2dx\vspace{3pt}\\
& \le & \displaystyle\int_{\Omega}A(x)\nabla u(x)\cdot\nabla u(x)=\int_{\Omega}-H(x,u(x),\nabla u(x))\,u(x)dx.\eaa$$
But $u\ge 0$ in $\Omega$ and
$$-H(x,u,\nabla u)\le a(x,u,\nabla u)\,|\nabla u|-b(x,u,\nabla u)\,u+f(x,u,\nabla u)\le M(1+|\nabla u|)$$
from~(\ref{hypH}) and~(\ref{hypAM}). Hence,~(\ref{C1}) yields
$$M^{-1}\int_{\Omega}|\nabla u|^2\le M\int_{\Omega}(1+|\nabla u|)\,u\le M\,C_1\int_{\Omega}(1+|\nabla u|),$$
from which one infers that
$$\|u\|_{H^1_0(\Omega)}\le C_2$$
for some positive constant $C_2$. It follows now from~(\ref{hypAM}) and~(\ref{C1}) that
\be\label{boundsH}
|H(x,u,\nabla u)|\le|H(x,u,\nabla u)-H(x,0,0)|+|H(x,0,0)|\le M(C_1+|\nabla u|)+M
\ee
in $\Omega$, whence $\|H(\cdot,u(\cdot),\nabla u(\cdot))\|_{L^2(\Omega)}\le C_3$ for some positive constant $C_3$. Standard elliptic estimates, together with~(\ref{hypAM}), then yield the existence of a positive constant $C_4$ such that~$\|u\|_{H^2(\Omega)}\le C_4$, whence $\|u\|_{W^{1,2^*}(\Omega)}\le C_5$ for some positive constant $C_5$ with $2^*=2n/(n-2)$ if $n\ge 3$ ($2^*=\infty$ if $n=1$, and $2^*$ denotes an arbitrarily fixed real number larger than $2$ if $n=2$). Using again~(\ref{boundsH}) and a standard bootstrap argument, it follows that $\|u\|_{W^{2,p}(\Omega)}\le \tilde{C}_p$ for every $1\le p<+\infty$ where $\tilde{C}_p$ depends only on $\Omega$, $n$, $M$ and $p$. In particular, one gets that
$$\|u\|_{C^{1,1/2}(\overline{\Omega})}\le C_6$$
for some positive constant $C_6$.\par
To complete the proof of Lemma~\ref{lembounds}, one just needs to show the existence of a positive constant~$\beta$, depending on $\Omega$, $n$ and $M$ but not on $u$, such that
$$u(x)\ge\beta\,d(x,\partial\Omega)\hbox{ for all }x\in\Omega.$$
Assume by contradiction that there is no such constant $\beta>0$. Then there are a sequence~$(A_m)_{m\in\N}$ of $W^{1,\infty}(\Omega,{\mathcal{S}}_n(\R))$ matrix fields, a sequence~$(\Lambda_m)_{m\in\N}$ of $L^{\infty}_+(\Omega)$ functions and four sequences~$(H_m)_{m\in\N}$, $(a_m)_{m\in\N}$, $(b_m)_{m\in\N}$, $(f_m)_{m\in\N}$ of continuous functions in~$\overline{\Omega}\times\R\times\R^n$, satisfying~(\ref{ALambda}),~(\ref{hypH}) and~(\ref{hypAM}) with the same parameter $M>0$, as well as a sequence~$(u_m)_{m\in\N}$ of~$W(\Omega)$ solutions of~(\ref{equ}) satisfying~(\ref{hypu}) and a sequence~$(x_m)_{m\in\N}$ of points in $\Omega$ such that
\be\label{xm}
u_m(x_m)<2^{-m}d(x_m,\partial\Omega).
\ee
From the previous paragraph, the sequence $(u_m)_{m\in\N}$ is actually bounded in $W^{2,p}(\Omega)$ for every $1\le p<+\infty$. Thus, the sequence $(H_m(\cdot,u_m(\cdot),\nabla u_m(\cdot)))_{m\in\N}$ is bounded in $L^{\infty}(\Omega)$ from~(\ref{hypAM}). As a consequence, there are a symmetric matrix field $A_{\infty}\in W^{1,\infty}(\Omega,{\mathcal{S}}_n(\R))$, a function $u_{\infty}\in W(\Omega)\cap H^1_0(\Omega)$, a point $x_{\infty}\in\overline{\Omega}$ and two functions $H_{\infty},\ H^0\in L^{\infty}(\Omega)$ such that, up to extraction of a subsequence and as $m\to+\infty$,
\be\label{Ainfty}\left\{\baa{l}\!
A_m\to A_{\infty}\hbox{ in }L^{\infty}(\Omega,{\mathcal{S}}_n(\R))\hbox{ with }A_{\infty}\ge M^{-1}\hbox{Id}\hbox{ in }\Omega,\vspace{3pt}\\
u_m\to u_{\infty}\hbox{ in }W^{2,p}(\Omega)\hbox{ weak for all }1\!\le\!p\!<\!+\infty\hbox{ and in }C^{1,\alpha}(\overline{\Omega})\hbox{ for all }0\!\le\!\alpha\!<\!1,\vspace{3pt}\\
x_m\to x_{\infty},\vspace{3pt}\\
H_m(\cdot,u_m(\cdot),\nabla u_m(\cdot))\rightharpoonup H_{\infty}\hbox{ and }H_m(\cdot,0,0)\rightharpoonup H^0\hbox{ in }L^{\infty}(\Omega)\hbox{ weak-*},\eaa\right.
\ee
and $u_{\infty}$ is a weak $H^1_0(\Omega)$ solution of
\be\label{uinfini}\left\{\baa{rcl}
-\div(A_{\infty}(x)\nabla u_{\infty})+H_{\infty}(x) & = & 0\ \hbox{ in }\Omega,\vspace{3pt}\\
u_{\infty} & = & 0 \ \hbox{ on }\partial\Omega,\eaa\right.
\ee
The function $u_{\infty}$ is then a strong $W(\Omega)$ solution of the above equation by elliptic regularity. Furthermore, $u_{\infty}\ge 0$ in $\Omega$ and it follows from~(\ref{xm}) that $u_{\infty}(x_{\infty})=0$ and that
\be\label{hopf}
\hbox{either }x_{\infty}\in\Omega,\hbox{ or }x_{\infty}\in\partial\Omega\hbox{ and }|\nabla u_{\infty}(x_{\infty})|=0.
\ee
For every $m\in\N$ and $x\in\Omega$, one has
$$H_m(x,u_m(x),\nabla u_m(x))\le M\big(u_m(x)+|\nabla u_m(x)|\big)+H_m(x,0,0),$$
from~(\ref{hypAM}) and the nonnegativity of $u_m$, whence
$$H_{\infty}(x)\le M\big(u_{\infty}(x)+|\nabla u_{\infty}(x)|\big)+H^0(x)\ \hbox{ a.e. in }\Omega$$
by passing to the $L^{\infty}(\Omega)$ weak-* limit as $m\to+\infty$. Therefore,
\be\label{uinfty}
-\div(A_{\infty}(x)\nabla u_{\infty})+M\big(u_{\infty}+|\nabla u_{\infty}|\big)\ge-H^0(x)\ge 0\ \hbox{ a.e. in }\Omega,
\ee
where the last inequality follows from the nonpositivity of $H_m(\cdot,0,0)$ by~(\ref{hypAM}). Since $u_{\infty}\ge 0$ in $\Omega$, $u_{\infty}(x_{\infty})=0$ and~(\ref{hopf}) holds, it follows from the strong maximum principle and Hopf lemma that $u_{\infty}$ is identically equal to~$0$ (see in particular Theorem~9.6 of~\cite{gt} -- and the discussion there around the Hopf lemma and the strong maximum principle for strong solutions -- even if it means changing the function $u_{\infty}$ into $U(x)=-e^{\lambda x}u_{\infty}(x)$ for some suitable $\lambda\in\R$). Therefore,~(\ref{uinfty}) implies that $H^0(x)=0$ for a.e. $x\in\Omega$. But
$$\int_{\Omega}H^0\le-M^{-1}<0$$
since $\int_{\Omega}H_m(\cdot,0,0)\le-M^{-1}$ from~(\ref{hypAM}) and $H_m(\cdot,0,0)\rightharpoonup H^0$ as $m\to+\infty$ in $L^{\infty}(\Omega)$ weak-*. One has then reached a contradiction.\par
Finally, there is a positive constant $\beta$ which depends on $\Omega$, $n$ and $M$ but not on $u$ such that $u\ge\beta\,d(\cdot,\partial\Omega)$ in $\Omega$. Finally, $u\in E_{1/2,N,\beta}(\Omega)$ for some positive constants $N$ and $\beta$ only depending on $\Omega$, $n$ and $M$. The proof of Lemma~\ref{lembounds} is thereby complete.\hfill\fin\break

\noindent{\bf{Proof of Lemma~\ref{lemsubpsik}.}} Let $k\in\N$ be fixed. The proof of the inequality~(\ref{subpsieps}) consists in integrating the inequality~(\ref{ineqpsieps}) satisfied by $\hat{\psi}_k$ in the shells of $E_k$ against a nonnegative test function $\varphi$ and then in controlling the boundary terms coming from the critical spheres of~$\partial E_k$. First of all, by elementary arguments, it is enough to prove that~(\ref{subpsieps}) holds for every~$\varphi\in C^1_c(\Omega^*)$ with $\varphi\ge0$ in~$\Omega^*$, where~$C^1_c(\Omega^*)$ denotes the set of~$C^1(\Omega^*)$ functions with compact support included in~$\Omega^*$. Let~$\varphi$ be any such nonnegative~$C^1_c(\Omega^*)$ function and call
$$I=\int_{\Omega^*}\hat{\Lambda}_k\nabla\hat{\psi}_k\cdot\nabla\varphi+\int_{\Omega^*}\big(\hat{a}_ke_r\cdot\nabla\hat{\psi}_k\big)\varphi-\int_{\Omega^*}g_k\varphi.$$
From~(\ref{Ek}), one can write
$$E_k\cap\Omega^*=S_{\rho_k(a^k_{m_k}),\rho_k(a^k_{m_k-1})}\cup\cdots\cup S_{\rho_k(a^k_2),\rho_k(a^k_1)}\cup S_{\rho_k(a^k_1),R},$$
where $\rho_k(a^k_{m_k})\!=\!0$, the integer $m_k$ is the number of the critical values $0<a^k_1\!<\!\cdots\!<\!a^k_{m_k}\!=\!\max_{\overline{\Omega}}u_{j_k}$ of the function $u_{j_k}$ in $\overline{\Omega}$ and $S_{\sigma,\sigma'}$ denotes the spherical shell
$$S_{\sigma,\sigma'}=\big\{x\in\R^n;\ \sigma<|x|<\sigma'\big\}$$
for any $0\le\sigma<\sigma'$. For convenience, define $a^k_0=0$, so that $\rho_k(a^k_0)=R$. It follows from Lebesgue's dominated convergence theorem that
\be\label{Ilimits}
I=\lim_{\gamma\to0^+}\sum_{i=0}^{m_k-1}\underbrace{\int_{S_{\rho_k(a^k_{i+1})+\gamma,\rho_k(a^k_i)-\gamma}}\Big(\hat{\Lambda}_k\nabla\hat{\psi}_k\cdot\nabla\varphi+\big(\hat{a}_ke_r\cdot\nabla\hat{\psi}_k\big)\varphi-g_k\varphi\Big)}_{=:I_{i,\gamma}}.
\ee
For every $0\le i\le m_k-1$ and every $\gamma$ such that
$$0<2\gamma<\gamma_0:=\min_{0\le\iota\le m_k-1}(\rho_k(a^k_{\iota})-\rho_k(a^k_{\iota+1})),$$
the functions $\hat{\Lambda}_k$ and $\hat{\psi}_k$ are of class $C^1(\overline{S_{\rho_k(a^k_{i+1})+\gamma,\rho_k(a^k_i)-\gamma}})$ and $C^2(\overline{S_{\rho_k(a^k_{i+1})+\gamma,\rho_k(a^k_i)-\gamma}})$ respectively, whence Green-Riemann formula implies that
$$\baa{rcl}
I_{i,\gamma} & = & \displaystyle\int_{S_{\rho_k(a^k_{i+1})+\gamma,\rho_k(a^k_i)-\gamma}}\underbrace{\big(-\hbox{div}\big(\hat{\Lambda}_k\nabla\hat{\psi}_k\big)+\hat{a}_ke_r\cdot\nabla\hat{\psi}_k-g_k\big)}_{\le 0}\,\underbrace{\varphi}_{\ge0}\vspace{3pt}\\
& & \displaystyle+\int_{\partial B_{\rho_k(a^k_i)-\gamma}}\hat{\Lambda}_k\nabla\hat{\psi}_k\cdot e_r\,\varphi\,d\theta_{\rho_k(a^k_i)-\gamma}-\int_{\partial B_{\rho_k(a^k_{i+1})+\gamma}}\hat{\Lambda}_k\nabla\hat{\psi}_k\cdot e_r\,\varphi\,d\theta_{\rho_k(a^k_{i+1})+\gamma},\eaa$$
where $d\theta_s$ denotes the surface measure on the sphere $\partial B_s$ for $s>0$. The first integral in the right-hand side of the above equality is nonpositive because of~(\ref{ineqpsieps}) and since $\varphi$ is nonnegative. On the other hand, for all $s\in[0,R]\backslash\big\{\rho_k(a^k_{\iota});\,1\le\iota\le m_k\big\}$ and for all $x\in\partial B_s$, it follows from~(\ref{defhatlambda}),~(\ref{defF}),~(\ref{defhatpsi}) and~(\ref{defuk}) that
$$\hat{\Lambda}_k(x)\nabla\hat{\psi}_k(x)\!\cdot\!e_r(x)=\frac{1}{n\alpha_ns^{n-1}}\int_{\Omega_{\rho_k^{-1}(s)}}\!\!\!\!\!\!\!\hbox{div}(A_{j_k}\nabla u_{j_k})=\frac{1}{n\alpha_ns^{n-1}}\int_{\Omega_{\rho_k^{-1}(s)}}\!\!\!\!\!\!\!H_{j_k}=:J(s).$$
Remember that $H_{j_k}$ is a continuous function in $\overline{\Omega}$ (it is a polynomial function) and that $|\Omega_{\rho_k^{-1}(s)}|=|B_s|=\alpha_ns^n$ and $|\Sigma_{\rho_k^{-1}(s)}|=0$ for each $s\in[0,R]$. Therefore, the function $J$, which can be defined for all $s\in(0,R]$ by the right-hand side of the above displayed equality, is continuous and bounded on $(0,R]$ (the boundedness of $J$ follows from the boundedness of $H_{j_k}$ and the fact that $|\Omega_{\rho_k^{-1}(s)}|=\alpha_ns^n$). Therefore, one has
$$I_{i,\gamma}\le J(\rho_k(a^k_i)-\gamma)\int_{\partial B_{\rho_k(a^k_i)-\gamma}}\varphi\,d\theta_{\rho_k(a^k_i)-\gamma}-J(\rho_k(a^k_{i+1})+\gamma)\int_{\partial B_{\rho_k(a^k_{i+1})+\gamma}}\varphi\,d\theta_{\rho_k(a^k_{i+1})+\gamma}$$
for every $0\le i\le m_k-1$ and every $\gamma\in(0,\gamma_0/2)$. Finally, using the continuity of $J$ on $[0,R]$ and of $\varphi$ on $\overline{\Omega^*}$, it follows that, for every $0\le i\le m_k-2$,
$$\limsup_{\gamma\to0^+}I_{i,\gamma}\le J(\rho_k(a^k_i))\int_{\partial B_{\rho_k(a^k_i)}}\varphi\,d\theta_{\rho_k(a^k_i)}-J(\rho_k(a^k_{i+1}))\int_{\partial B_{\rho_k(a^k_{i+1})}}\varphi\,d\theta_{\rho_k(a^k_{i+1})},$$
while, for $i=m_k-1$,
$$
\limsup_{\gamma\to0^+}I_{i,\gamma}\le J(\rho_k(a^k_i))\int_{\partial B_{\rho_k(a^k_i)}}\varphi\,d\theta_{\rho_k(a^k_i)},
$$
since $J$ is bounded on $(0,R]$ and $\lim_{\rho\rightarrow 0^+}\int_{\partial B_{\rho}}\varphi\,d\theta_{\rho}=0$. As a conclusion,~(\ref{Ilimits}) yields
$$I\le J(\rho_k(a^k_0))\int_{\partial B_{\rho_k(a^k_0)}}\varphi\,d\theta_{\rho_k(a^k_0)}=J(R)\int_{\partial\Omega^*}\varphi\,d\theta_R=0$$
since $a^k_0=0$ (by convention), $\rho_k(a^k_0)=R$ and $\varphi$ is compactly supported in $\Omega^*$. As already emphasized,~(\ref{subpsieps}) then holds for every $\varphi\in H^1_0(\Omega^*)$ such that $\varphi\ge 0$ a.e. in $\Omega^*$. The proof of Lemma~\ref{lemsubpsik} is thereby complete.\hfill\fin\break

\noindent{\bf{Proof of Lemma~\ref{vinfty}.}} Equation~(\ref{eqvlimit}) is obtained formally by passing to the limit as $k\to+\infty$ in~(\ref{defvk}). However, since the convergence of the first-order derivatives of the functions~$v_k$ is only weak, as is that of the coefficients of~(\ref{defvk}) or some functions of them, one cannot pass directly to the limit in the first two terms of~(\ref{defvk}) and one needs more regularity. This regularity will be guaranteed by the radial symmetry, as shown in the next paragraphs.\par
Recall first that the sets~$E_k\subset\overline{\Omega^*}$ are given in~(\ref{Ek}). For every $k\in\N$, since $\hat{\Lambda}_k$ is of class~$C^1$ in~$E_k\cap\Omega^*$ and the functions $\hat{a}_k$ and $g_k$ are continuous in $E_k$, it follows from standard elliptic estimates that the function $v_k$ is in $W^{2,p}_{loc}(E_k\cap\Omega^*)$ for all $1\le p<+\infty$ whence in~$C^{1,\alpha}_{loc}(E_k\cap\Omega^*)$ for all $0\le\alpha<1$. Define
$${\mathcal{I}}_k=(0,R)\backslash\big\{\rho_k(a^k_i);\ 1\le i\le m_k\big\}=(0,\rho_k(a^k_{m_k-1}))\cup\cdots\cup(\rho_k(a^k_1),R),$$
where $0<a^k_1<\cdots<a^k_{m_k}=\max_{\overline{\Omega}}u_{j_k}$ are the critical values of the function $u_{j_k}$ in $\overline{\Omega}$. Let $\varsigma$ be any point in the unit sphere $\mathbb{S}^{n-1}$ and set
\be\label{deftilde}\left\{\baa{l}
\tilde{v}_k(r)=v_k(r\varsigma),\ \ \tilde{\Lambda}_k(r)=\hat{\Lambda}_k(r\varsigma),\vspace{3pt}\\
\tilde{w}_k(r)=\hat{\Lambda}_k(r\varsigma)\,\nabla v_k(r\varsigma)\cdot\varsigma=\tilde{\Lambda}_k(r)\,\tilde{v}_k'(r),\vspace{3pt}\\
\tilde{a}_k(r)=\hat{a}_k(r\varsigma),\ \tilde{g}_k(r)=g_k(r\varsigma),\eaa\right.
\ee
for all $r\in{\mathcal{I}}_k$. Denote also
\be\label{deftildebis}\left\{\baa{ll}
\displaystyle\tilde{\Lambda}(r)=\frac{1}{n\alpha_nr^{n-1}}\int_{\partial B_r}\hat{\Lambda}\,d\theta_r, & \displaystyle\tilde{v}_{\infty}(r)=\frac{1}{n\alpha_nr^{n-1}}\int_{\partial B_r}v_{\infty}\,d\theta_r,\vspace{3pt}\\
\displaystyle\tilde{a}(r)=\frac{1}{n\alpha_nr^{n-1}}\int_{\partial B_r}\hat{a}\,d\theta_r, & \displaystyle\tilde{f}(r)=\frac{1}{n\alpha_nr^{n-1}}\int_{\partial B_r}\hat{f}\,d\theta_r.\eaa\right.
\ee
From Fubini's theorem, the above quantities can be defined for almost every $r\in(0,R)$. Furthermore, the function~$\tilde{\Lambda}$ is in $L^{\infty}_+(0,R)$, the functions~$\tilde{a}$ and~$\tilde{f}$ are in $L^{\infty}(0,R)$, the function $\tilde{v}_{\infty}$ is in $H^1_{loc}((0,R])$, and it follows from~(\ref{vkv}),~(\ref{chapeaux}) and~(\ref{deftilde}) that
\be\label{tildes}\baa{l}
\tilde{\Lambda}_k^{-1}\rightharpoonup\tilde{\Lambda}^{-1},\ \ \tilde{\Lambda}_k^{-1}\tilde{a}_k\rightharpoonup\tilde{\Lambda}^{-1}\tilde{a}\ \hbox{ and }\ \tilde{g}_k\rightharpoonup\tilde{f},\ \hbox{ in }L^{\infty}(r_0,R)\hbox{ weak-*, }\hbox{ as }k\to+\infty,\vspace{3pt}\\
\qquad\qquad\qquad\qquad\qquad\qquad\qquad\qquad\qquad\qquad\qquad\qquad\qquad\hbox{ for every }r_0\in(0,R)\eaa
\ee
and
\be\label{vkvbis}\baa{l}
\tilde{v}_k\rightharpoonup\tilde{v}_{\infty}\hbox{ in }H^1(r_0,R)\hbox{ weak }\hbox{ and }\ \tilde{v}_k\to\tilde{v}_{\infty}\hbox{ in }L^2(r_0,R)\hbox{ strong, }\hbox{ as }k\to+\infty,\vspace{3pt}\\
\qquad\qquad\qquad\qquad\qquad\qquad\qquad\qquad\qquad\qquad\qquad\qquad\qquad\hbox{ for every }r_0\in(0,R).\eaa
\ee\par
Let us now pass to the limit as $k\to+\infty$ in the elliptic partial differential equation~(\ref{defvk}) and its associated one-dimensional ordinary differential equation. Namely, from the observations of the previous paragraph, there holds
\be\label{eqtildevk}
-\tilde{w}_k'=-\big(\tilde{\Lambda}_k\tilde{v}_k'\big)'=\frac{n-1}{r}\tilde{\Lambda}_k\tilde{v}_k'-\tilde{a}_k\tilde{v}_k'+\tilde{g}_k\ \hbox{ a.e. in }{\mathcal{I}}_k
\ee
for every $k\in\N$. Since the right-hand side of~(\ref{eqtildevk}) is continuous in ${\mathcal{I}}_k$, the continuous function~$\tilde{w}_k=\tilde{\Lambda}_k\tilde{v}_k'$ is actually of class $C^1$ in ${\mathcal{I}}_k$, whence $\tilde{v}_k$ is of class $C^2$ in ${\mathcal{I}}_k$ and the above equation~(\ref{eqtildevk}) is satisfied in the classical pointwise sense in ${\mathcal{I}}_k$. Furthermore, the sequence~$(\tilde{v}_k')_{k\in\N}$ is bounded in $L^2(r_0,R)$ for every $r_0\in(0,R)$ since $(v_k)_{k\in\N}$ is bounded in $H^1_0(\Omega^*)$ and the functions $v_k$ are radially symmetric. On the other hand, the sequences~$(\tilde{\Lambda}_k)_{k\in\N}$,~$(\tilde{a}_k)_{k\in\N}$ and~$(\tilde{g}_k)_{k\in\N}$ are bounded in~$L^{\infty}(0,R)$ from~(\ref{Aklambdak}),~(\ref{coefepsilon}),~(\ref{defalphak}) and~(\ref{deftilde}). It follows then from~(\ref{deftilde}) and~(\ref{eqtildevk}) that the sequence $(\tilde{w}_k)_{k\in\N}$ is bounded in $H^1(r_0,R)$ for every $r_0\in(0,R)$. Therefore, there is a function $\tilde{w}\in H^1_{loc}((0,R])$ such that, up to extraction of a subsequence,
\be\label{tildew}
\tilde{w}_k\mathop{\rightharpoonup}_{k\to+\infty}\tilde{w}\hbox{ in }H^1(r_0,R)\hbox{ weak},\ \tilde{w}_k\mathop{\rightarrow}_{k\to+\infty}\tilde{w}\hbox{ in }L^2(r_0,R)\hbox{ strong, for every }r_0\in(0,R).
\ee
By \eqref{tildes}, one therefore has
$$\tilde{v}_k'=\tilde{\Lambda}_k^{-1}\tilde{w}_k\mathop{\rightharpoonup}_{k\to+\infty}\tilde{\Lambda}^{-1}\tilde{w}\hbox{ in }L^2(r_0,R)\hbox{ weak, for every }r_0\in(0,R).$$
By uniqueness of the weak limit, it follows then from~(\ref{vkvbis}) that
\be\label{tildewbis}
\tilde{w}=\tilde{\Lambda}\,\tilde{v}'_{\infty}\hbox{ a.e. in }(0,R).
\ee\par
Finally, we are ready to show that $v_{\infty}$ is a weak $H^1_0(\Omega^*)$ solution of the limiting equation~(\ref{eqvlimit}), that is
$$I:=\int_{\Omega^*}\hat{\Lambda}\,\nabla v_{\infty}\cdot\nabla\varphi+\int_{\Omega^*}\big(\hat{a}\,e_r\cdot\nabla v_{\infty}\big)\varphi-\int_{\Omega^*}\hat{f}\,\varphi=0$$
for all $\varphi\in H^1_0(\Omega^*)$. Let~$\varphi\in H^1_0(\Omega^*)$ and $\epsilon>0$ be arbitrary. Since the functions $\hat{\Lambda}$,~$\hat{a}$,~$\hat{f}$ are in $L^{\infty}(\Omega^*)$, since the sequences~$(\hat{\Lambda}_k)_{k\in\N}$,~$(\hat{a}_k)_{k\in\N}$,~$(g_k)_{k\in\N}$ are bounded in $L^{\infty}(\Omega^*)$, since the function $v_{\infty}$ is in~$H^1_0(\Omega^*)$ and since the sequence $(v_k)_{k\in\N}$ is bounded in~$H^1_0(\Omega^*)$, it follows from Cauchy-Schwarz inequality and the fact that $\|\varphi\|_{H^1(B_r)}\to0$ as $r\to 0^+$ (from Lebesgue's dominated convergence theorem) that there exists $r_0\in(0,R)$ small enough so that
\be\label{kepsilon}
\Big|\int_{B_{r_0}}\Big(\hat{\Lambda}_k\nabla v_k\cdot\nabla\varphi+(\hat{a}_ke_r\cdot\nabla v_k)\varphi-g_k\varphi\Big)\Big|\le\epsilon\ \hbox{ for all }k\in\N
\ee
and
$$\Big|\int_{B_{r_0}}\Big(\hat{\Lambda}\,\nabla v_{\infty}\cdot\nabla\varphi+(\hat{a}\,e_r\cdot\nabla v_{\infty})\varphi-\hat{f}\,\varphi\Big)\Big|\le\epsilon.$$
Hence,
\be\label{IJ}
|I|\le\epsilon+\Big|\underbrace{\int_{\Omega^*\backslash B_{r_0}=S_{r_0,R}}\Big(\hat{\Lambda}\,\nabla v_{\infty}\cdot\nabla\varphi+(\hat{a}\,e_r\cdot\nabla v_{\infty})\varphi-\hat{f}\,\varphi\Big)}_{=:J}\Big|.
\ee
Call
\be\label{defphiPhi}
\phi(r)=\int_{\partial B_r}\varphi\,d\theta_r\ \hbox{ and }\ \Phi(r)=\int_{\partial B_r}e_r\cdot\nabla\varphi\,d\theta_r\ \hbox{ for }r\in(r_0,R).
\ee
From Fubini's theorem and Cauchy-Schwarz inequality, the functions $\phi$ and $\Phi$ are defined almost everywhere in $(r_0,R)$ and they belong to $L^2(r_0,R)$. Fubini's theorem also implies that
$$J=\int_{r_0}^R\Big(\tilde{\Lambda}(r)\,\tilde{v}_{\infty}'(r)\,\Phi(r)+\tilde{a}(r)\,\tilde{v}_{\infty}'(r)\,\phi(r)-\tilde{f}(r)\,\phi(r)\Big)\,dr.$$
Observe that
\be\label{tildelambdav}
\tilde{\Lambda}_k\tilde{v}_k'=\tilde{w}_k\to\tilde{w}=\tilde{\Lambda}\,\tilde{v}_{\infty}'\hbox{ in }L^2(r_0,R)\hbox{ strong as }k\to+\infty
\ee
from~(\ref{deftilde}),~(\ref{tildew}) and~(\ref{tildewbis}). Furthermore,
$$\tilde{a}_k\tilde{v}_k'=(\tilde{\Lambda}_k^{-1}\tilde{a}_k)\times(\tilde{\Lambda}_k\tilde{v}_k')\rightharpoonup(\tilde{\Lambda}^{-1}\tilde{a})\times(\tilde{\Lambda}\,\tilde{v}_{\infty}')=\tilde{a}\,\tilde{v}_{\infty}'\hbox{ in }L^2(r_0,R)\hbox{ weak as }k\to+\infty$$
from~(\ref{tildes}) and~(\ref{tildelambdav}). Lastly, $\tilde{g}_k\rightharpoonup\tilde{f}$ in $L^{\infty}(r_0,R)$ weak-* as $k\to+\infty$ from~(\ref{tildes}). Putting together all this limits leads to
$$J=\lim_{k\to+\infty}\Big[\underbrace{\int_{r_0}^R\Big(\tilde{\Lambda}_k(r)\,\tilde{v}_k'(r)\,\Phi(r)+\tilde{a}_k(r)\,\tilde{v}_k'(r)\,\phi(r)-\tilde{g}_k(r)\,\phi(r)\Big)\,dr}_{=:J_k}\Big].$$
Therefore, there is $K\in\N$ large enough such that
\be\label{JJK}
|J|\le|J_K|+\epsilon.
\ee
But
$$\baa{rcl}
J_K & = & \displaystyle\int_{\Omega^*\backslash B_{r_0}}\Big(\hat{\Lambda}_K\,\nabla v_K\cdot\nabla\varphi+(\hat{a}_Ke_r\cdot\nabla v_K)\varphi-g_K\,\varphi\Big)\vspace{3pt}\\
& = & \displaystyle-\int_{B_{r_0}}\Big(\hat{\Lambda}_K\,\nabla v_K\cdot\nabla\varphi+(\hat{a}_Ke_r\cdot\nabla v_K)\varphi-g_K\,\varphi\Big)\eaa$$
because the integral of $\hat{\Lambda}_K\,\nabla v_K\cdot\nabla\varphi+(\hat{a}_Ke_r\cdot\nabla v_K)\varphi-g_K\,\varphi$ over $\Omega^*$ is equal to $0$: indeed,~$\varphi\in H^1_0(\Omega^*)$ and~$v_K$ is the weak $H^1_0(\Omega^*)$ solution of~(\ref{defvk}) with $k=K$. Finally, it follows from~(\ref{kepsilon}) that $|J_K|\le\epsilon$, whence
$$|I|\le3\epsilon$$
from~(\ref{IJ}) and~(\ref{JJK}). Since $\epsilon>0$ was arbitrary, one concludes that $I=0$. Since $\varphi\in H^1_0(\Omega^*)$ was arbitrary, this means that $v_{\infty}$ is the weak $H^1_0(\Omega^*)$ solution of~(\ref{eqvlimit}). The proof of Lemma~\ref{vinfty} is thereby complete.\hfill\fin\break

\noindent{\bf{Proof of Lemma~\ref{lemzLz}.}} Remember that $z_L$ and $z$ are the unique weak $H^1_0(\Omega^*)$ solutions of~(\ref{eqzl}) with $l=L$ and~(\ref{eqz}), respectively. Furthermore, $z_L$ and $z$ are radially symmetric. In order to get the inequality~(\ref{zLz}), the general strategy is to integrate twice the one-dimensional equations associated to~(\ref{eqzl}) and~(\ref{eqz}), with a special care since these equations are only satisfied in the weak $H^1_0(\Omega^*)$ sense. We shall finally use the Hardy-Littlewood inequality to compare some integral terms.\par
First of all, from the De Giorgi-Moser-Nash regularity theory and the radial symmetry, it follows that, without loss of generality, $z_L$ and $z$ can be assumed to be continuous in $\overline{\Omega^*}$, even if it means redefining them on a negligible subset of $\overline{\Omega^*}$. Let $\varsigma$ be any point in the unit sphere $\mathbb{S}^{n-1}$ and define
\be\label{deftildezL}
\tilde{z}_L(r)=z_L(r\varsigma)\ \hbox{ and }\ \tilde{z}(r)=z(r\varsigma)\ \hbox{ for all }r\in[0,R].
\ee
The functions $\tilde{z}_L$ and $\tilde{z}$ are continuous on $[0,R]$, they belong to $H^1_{loc}((0,R])$ and, from Fubini's theorem, the integrals
\be\label{integrals}
\int_0^Rr^{n-1}\,\tilde{z}_L(r)^2\,dr,\ \int_0^Rr^{n-1}\,\tilde{z}'_L(r)^2\,dr,\ \int_0^Rr^{n-1}\,\tilde{z}(r)^2\,dr,\ \int_0^Rr^{n-1}\,\tilde{z}'(r)^2\,dr\hbox{ converge}.
\ee
Let $\tilde{\Lambda}\in L^{\infty}_+(0,R)$ and $\tilde{a}\in L^{\infty}(0,R)$ be defined as in~(\ref{deftildebis}) and set
\be\label{deftildeter}
\displaystyle\tilde{h}_L(r)=\frac{1}{n\alpha_nr^{n-1}}\int_{\partial B_r}h_L\,d\theta_r\hbox{ and }\tilde{f}_u(r)=\frac{1}{n\alpha_nr^{n-1}}\int_{\partial B_r}f_u^*\,d\theta_r,
\ee
where we recall that $h_L$ (resp. $f^*_u$) is the right-hand side of~(\ref{eqzl}) (resp.~(\ref{eqz})). These quantities can be defined for almost every $r\in(0,R)$ ($\tilde{f}_u$ is actually defined for all $0<r<R$) from Fubini's theorem, and $\tilde{h}_L$ and $\tilde{f}_u$ are in $L^{\infty}(0,R)$.\par
Consider now the equation~(\ref{eqzl}) with $l=L$. It follows from the definitions of $z_L$,~$\tilde{\Lambda}$,~$\tilde{a}$ and~$\tilde{h}_L$ that
$$\int_0^R\tilde{\Lambda}(r)\,\tilde{z}_L'(r)\,\varphi'(r)\,r^{n-1}\,dr+\int_0^R\tilde{a}(r)\,\tilde{z}_L'(r)\,\varphi(r)\,r^{n-1}\,dr=\int_0^R\tilde{h}_L(r)\,\varphi(r)\,r^{n-1}\,dr$$
for all $\varphi\in C^1_c(0,R)$. Define
\be\label{defzetaL}
\zeta_L(r)=r^{n-1}\,\tilde{\Lambda}(r)\,\tilde{z}_L'(r)
\ee
for almost every $r\in(0,R)$. Since $r\mapsto r^{(n-1)/2}\tilde{\Lambda}(r)$ is in $L^{\infty}(0,R)$, it follows from~(\ref{integrals}) that~$\zeta_L$ is in $L^2(0,R)$. There holds
$$\int_0^R\zeta_L(r)\,\varphi'(r)\,dr+\int_0^R\tilde{\Lambda}^{-1}(r)\,\tilde{a}(r)\,\zeta_L(r)\,\varphi(r)\,dr=\int_0^R\tilde{h}_L(r)\,\varphi(r)\,r^{n-1}\,dr$$
for all $\varphi\in C^1_c(0,R)$. Furthermore, $\tilde{\Lambda}^{-1}\,\tilde{a}\,\zeta_L\in L^2(0,R)$ and the function $r\mapsto\tilde{h}_L(r)\,r^{n-1}$ is in~$L^{\infty}(0,R)\subset L^2(0,R)$. Therefore, the function $\zeta_L$ is actually in $H^1(0,R)$ and
\be\label{zeta'L}
\zeta'_L(r)=\tilde{\Lambda}^{-1}(r)\,\tilde{a}(r)\,\zeta_L(r)-\tilde{h}_L(r)\,r^{n-1}\ \hbox{ a.e. in }(0,R).
\ee
Even if it means redefining $\zeta_L$ on a negligible subset of $[0,R]$, one can assume without loss of generality that $\zeta_L$ is then continuous on $[0,R]$. Define
$$\tilde{\Theta}(r)=\displaystyle e^{-\int_0^r\tilde{\Lambda}^{-1}(s)\,\tilde{a}(s)\,ds}\ \hbox{ for all }r\in[0,R].$$
The function $\tilde{\Theta}$, which does not depend on $L$, is continuous on $[0,R]$ and it belongs to~$W^{1,\infty}(0,R)$. Define also
\be\label{defomegaL}
\omega_L=\tilde{\Theta}\,\zeta_L\ \hbox{ on }[0,R].
\ee
The function $\omega_L$ is continuous on $[0,R]$ and it belongs to $H^1(0,R)$. It follows from~(\ref{zeta'L}) that
$$\omega'_L(r)=-\tilde{\Theta}(r)\,\tilde{h}_L(r)\,r^{n-1}\ \hbox{ a.e. in }(0,R).$$
As a consequence,
\be\label{omegaL}
\omega_L(r)=-\int_0^r\tilde{\Theta}(s)\,\tilde{h}_L(s)\,s^{n-1}\,ds+\omega_L(0)\ \hbox{ for all }r\in[0,R].
\ee\par
Let us now prove in this paragraph that $\omega_L(0)=0$ (remember that $\omega_L$ is continuous on~$[0,R]$). One has $\omega_L(0)=\zeta_L(0)$ since $\tilde{\Theta}(0)=1$. Consider first the case where the dimension~$n$ is such that $n\ge 2$. If $\zeta_L(0)\neq0$, then there would exist $r_0\in(0,R)$ and $\gamma>0$ such that~$|\zeta_L(r)|\ge\gamma>0$ for all $r\in[0,r_0]$, whence
$$\big|r^{n-1}\,\tilde{\Lambda}(r)\,\tilde{z}_L'(r)\big|\ge\gamma>0\ \hbox{ for a.e. }r\in(0,r_0).$$
Since $\tilde{\Lambda}\in L^{\infty}_+(0,R)$ and $\tilde{\Lambda}\le M_{\Lambda}:={\rm{ess}}\,{\rm{sup}}_{\Omega}\Lambda$ from~(\ref{hats2}) and~(\ref{deftildebis}), it would then follow that
$$r^{n-1}\,\tilde{z}_L'(r)^2\ge\frac{\gamma^2}{M_{\Lambda}^2r^{n-1}}\ \hbox{ for a.e. }r\in(0,r_0),$$
which contradicts the integrability of the integral~$\int_0^Rr^{n-1}\,\tilde{z}_L'(r)^2\,dr$ given in~(\ref{integrals}). Hence, $\zeta_L(0)=0$ and $\omega_L(0)=0$. \par
Consider now the case $n=1$. We work directly in the interval $\Omega^*=(-R,R)$. With similar arguments as in the previous paragraph, the $L^2(-R,R)$ function $\zeta^1_L:=\hat{\Lambda}\,z_L'$ is actually in~$H^1(-R,R)$, whence continuous on $[-R,R]$ without loss of generality, and $\zeta^1_L(x)=\zeta_L(x)$ for all $x\in[0,R]$. But since the $H^1(-R,R)$ function $z_L$ is even, the $L^2(-R,R)$ function $z_L'$ is odd, in the sense that
$$z_L'(-x)=-z_L'(x)\hbox{ for a.e. }x\in(-R,R).$$
Since $\hat{\Lambda}$ is even, the (continuous) function $\zeta^1_L=\hat{\Lambda}\,z_L'$ is odd on $[0,R]$. In particular, it vanishes at~Ê$0$. Hence, $\omega_L(0)=\zeta_L(0)=\zeta^1_L(0)=0$. To sum up, there holds
$$\omega_L(0)=0$$
in all dimensions $n\ge 1$.\par
From~(\ref{omegaL}), one then infers that
$$\omega_L(r)=-\int_0^r\tilde{\Theta}(s)\,\tilde{h}_L(s)\,s^{n-1}\,ds\ \hbox{ for all }r\in[0,R],$$
whence
\be\label{tildezL}
\tilde{\Theta}(r)\,r^{n-1}\,\tilde{\Lambda}(r)\,\tilde{z}_L'(r)=-\int_0^r\tilde{\Theta}(s)\,\tilde{h}_L(s)\,s^{n-1}\,ds\ \hbox{ for a.e. }r\in(0,R)
\ee
from~(\ref{defzetaL}) and~(\ref{defomegaL}). Similarly, by working with the equation~(\ref{eqz}) satisfied by $z$ with right-hand side $f^*_u$ (instead of $h_L$ in~(\ref{eqzl}) with $l=L$) and by using the notations~(\ref{deftildezL}) and~(\ref{deftildeter}), one gets that
\be\label{tildez}
\tilde{\Theta}(r)\,r^{n-1}\,\tilde{\Lambda}(r)\,\tilde{z}'(r)=-\int_0^r\tilde{\Theta}(s)\,\tilde{f}_u(s)\,s^{n-1}\,ds\ \hbox{ for a.e. }r\in(0,R).
\ee\par
Let us finally compare the right-hand sides of~(\ref{tildezL}) and~(\ref{tildez}), which are actually continuous with respect to $r\in[0,R]$. Define
$$\Theta(x)=\tilde{\Theta}(|x|)=\displaystyle e^{-\int_0^{|x|}\tilde{\Lambda}^{-1}(s)\,\tilde{a}(s)\,ds}\ \hbox{ for every }x\in\overline{\Omega^*}.$$
The function $\Theta$ is continuous in $\overline{\Omega^*}$, positive, radially symmetric and nonincreasing with respect to $|x|$ since $\tilde{\Lambda}>0$ and $\tilde{a}\ge 0$ a.e. in $(0,R)$, by~(\ref{hats2}),~(\ref{hats3}) and~(\ref{deftildebis}). In particular, the function $\Theta$ is equal to its Schwarz rearrangement $\Theta^*$ in the ball $\Omega^*$. For every $r\in(0,R]$, it follows then from Fubini's theorem and Hardy-Littlewood inequality that
\be\label{hardy}\baa{rcl}
\displaystyle\int_0^r\tilde{\Theta}(s)\,\tilde{h}_L(s)\,s^{n-1}\,ds & = & \displaystyle\frac{1}{n\alpha_n}\int_{B_r}\Theta(x)\,h_L(x)\,dx\vspace{3pt}\\
& \le & \displaystyle\frac{1}{n\alpha_n}\int_{B_r}\Theta^*(x)\,(h^r_L)^*(x)\,dx=\frac{1}{n\alpha_n}\int_{B_r}\Theta(x)\,(h^r_L)^*(x)\,dx,\eaa
\ee
where $(h^r_L)^*$ denotes the Schwarz rearrangement, in the ball $B_r$, of the restriction $h^r_L:=(h_L)_{|B_r}$ of the function $h_L$ in the ball $B_r$ (notice that the Hardy-Littlewood inequality is usually stated for nonnegative functions; here, the function $\Theta$ is nonnegative, but the function $h_L$ may not be nonnegative in general; however, due to the definition of the Schwarz rearrangement given in Section~\ref{intro} for general $L^1$ functions with no sign, the inequality $\int_{B_r}\Theta\,h_L\le\int_{B_r}\Theta^*(h^r_L)^*$ holds immediately from the standard  Hardy-Littlewood inequality, since $(h+\lambda)^*=h^*+\lambda$ for any $L^1$ function $h$ and any constant $\lambda\in\R$). On the other hand, elementary arguments imply that $(h^r_L)^*\le(h_L^*)_{|B_r}$ in~$B_r$, where~$(h_L^*)_{|B_r}$ denotes the restriction in $B_r$ of the Schwarz rearrangement $h_L^*$ of the function~$h_L$ in~$\Omega^*$. Hence
\be\label{schwarz}
(h^r_L)^*\le(f^*_u)_{|B_r}\ \hbox{ in }B_r
\ee
since $h_L^*=f_u^*$ in $\Omega^*$ (the functions $h_L\in L^{\infty}(\Omega^*)$ and $f_u\in L^{\infty}(\Omega)$ have indeed the same distribution functions $\mu_{h_L}=\mu_{f_u}$ and thus the same Schwarz rearrangements in $\Omega^*$). Since $\Theta\ge0$ in $\Omega^*$, one infers from~(\ref{hardy}),~(\ref{schwarz}) and Fubini's theorem that, for all $r\in(0,R]$,
$$\baa{l}
\displaystyle\int_0^r\tilde{\Theta}(s)\,\tilde{h}_L(s)\,s^{n-1}\,ds\le\frac{1}{n\alpha_n}\int_{B_r}\Theta(x)\,(h^r_L)^*(x)\,dx\le\frac{1}{n\alpha_n}\int_{B_r}\Theta(x)\,f_u^*(x)\,dx\vspace{3pt}\\
\qquad\qquad\qquad\qquad\qquad\qquad\qquad\qquad\qquad\qquad\qquad\qquad\qquad=\displaystyle\int_0^r\tilde{\Theta}(s)\,\tilde{f}_u(s)\,s^{n-1}\,ds.\eaa$$
Together with~(\ref{tildezL}) and~(\ref{tildez}), it follows that
$$-\tilde{\Theta}(r)\,r^{n-1}\,\tilde{\Lambda}(r)\,\tilde{z}_L'(r)\le-\tilde{\Theta}(r)\,r^{n-1}\,\tilde{\Lambda}(r)\,\tilde{z}'(r)\ \hbox{ for a.e. }r\in(0,R),$$
whence
$$-\tilde{z}_L'(r)\le-\tilde{z}'(r)\ \hbox{ for a.e. }r\in(0,R)$$
since $\tilde{\Lambda}\in L^{\infty}_+(0,R)$ and $\tilde{\Theta}$ is continuous and positive on $[0,R]$. Finally, since both functions~$\tilde{z}_L$ and $\tilde{z}$ are continuous on $[0,R]$, vanish at $R$ ($z_L$ and $z$ are in $H^1_0(\Omega^*)$) and are in $H^1_{loc}((0,R])$, one gets that, for all~$r\in[0,R]$,
$$\tilde{z}_L(r)=\int_r^R-\tilde{z}_L'(s)\,ds\le\int_r^R-\tilde{z}'(s)\,ds=\tilde{z}(r).$$
As a conclusion, remembering the definitions~(\ref{deftildezL}) and the radial symmetry and continuity of $z_L$ and $z$ in $\overline{\Omega^*}$, one concludes that
$$z_L(x)\le z(x)\hbox{ for all }x\in\overline{\Omega^*}.$$
The proof of Lemma~\ref{lemzLz} is thereby complete.\hfill\fin


\SE{General growth with respect to the gradient}\label{sec4}

This section is devoted to the proofs of Theorems~\ref{th2} and~\ref{th2bis}. Throughout this section, we assume~(\ref{ALambda}), (\ref{hypH}) and~(\ref{hypb}) with $\Lambda\in L^{\infty}_+(\Omega)$ and
$$1<q\le 2,$$
that is the nonlinear function $H$ is bounded from below by an at most quadratic function of~$|p|$. Furthermore, $u\in W(\Omega)$ denotes a solution of~(\ref{equ}) satisfying~(\ref{hypu}), that is $u>0$ in~$\Omega$ and~$|\nabla u|\neq 0$ on $\partial\Omega$. We recall that, even if it means redefining $u$ on a negligible subset of~$\overline{\Omega}$, one can assume without loss of generality that $u\in C^{1,\alpha}(\overline{\Omega})$ for all $\alpha\in[0,1)$.\par
Our goal is to establish the inequalities~(\ref{u*vbis}) and~(\ref{u*vepsilon}) and the quantified ones~(\ref{u*vetaubis}), (\ref{u*vetaubiseps}), (\ref{u*vetabis}) and~(\ref{u*vetabiseps}), that is to compare the Schwarz rearrangement $u^*$ of $u$ with the unique solutions~$v$ and~$v_{\epsilon}$ of~(\ref{eqv4}) and~(\ref{eqv4eps}). The strategy follows a similar scheme to that used in the previous section for the proofs of Theorems~\ref{th1} and~\ref{th1bis}. Namely, after establishing some uniform bounds on $u$ under assumption~(\ref{hypAMquadr}), we first approximate $u$ by some smooth solutions $u_j$ of some regularized equations in $\Omega$. Next, we apply the general rearrangement inequalities of Section~\ref{rearrangement} and we compare some $u^*_{j_k}$ with the solutions $v_k$ of some symmetrized equations in $\Omega^*$. Lastly, we approximate the functions $v_k$ by some solutions of equations of the type~(\ref{eqv4}) and~(\ref{eqv4eps}) in $\Omega^*$.\par
The proofs of Theorems~\ref{th2} and~\ref{th2bis} are done in Section~\ref{sec41} and the proofs of two auxiliary technical lemmas are carried out in Section~\ref{sec42}.


\subsection{Proofs of Theorems~\ref{th2} and~\ref{th2bis}}\label{sec41}

\subsubsection*{Step 1: uniform bounds on $u$ under assumption~(\ref{hypAMquadr})}

In this step, some uniform pointwise and smoothness estimates  are established under assumption~(\ref{hypAMquadr}). Actually, these quantified estimates will only be needed for Theorem~\ref{th2bis}, in which~$\Omega$ is not a ball. We recall that the sets $E_{\alpha,N,\beta}(\Omega)$ have been defined in Section~\ref{sec24}.

\begin{lem} \label{lemboundsquadr}
Under assumption~$(\ref{hypAMquadr})$ with $n\ge 2$, there are some real numbers $N>0$ and~$\beta>0$, which depend only on $\Omega$, $n$, $q$, $M$ and $r$, such that $u\in E_{1/2,N,\beta}(\Omega)$.
\end{lem}

The proof of Lemma \ref{lemboundsquadr}, which is a version of Lemma~\ref{lembounds} adapted to the case where $1<q\leq 2$, can be found in Section~\ref{sec42} below

\subsubsection*{Step 2: approximated coefficients and approximated solutions $u_j$ in $\Omega$}

This step is the same as Step~2 of the proofs of Theorems~\ref{th1} and~\ref{th1bis}. Namely, the sequences~$(H_j)_{j\in\N}$, $(A_j)_{j\in\N}=((A_{j;i,i'})_{1\le i,i'\le n})_{j\in\N}$, $(\Lambda_j)_{j\in\N}$ and $(u_j)_{j\in\N}$ satisfy~(\ref{Hkinfty}) and~(\ref{defuk}), as well as~(\ref{Aklambdak}), that is
\be\label{Aklambdakbis}\left\{\baa{l}
\displaystyle A_{j;i,i'}\mathop{\longrightarrow}_{j\to+\infty}A_{i,i'}\hbox{ uniformly in }\overline{\Omega}\hbox{ for all }1\le i,i'\le n,\ \ \displaystyle\mathop{\sup}_{j\in\N}\|A_j\|_{W^{1,\infty}(\Omega)}<+\infty,\vspace{3pt}\\
A_j\ge\Lambda_j\hbox{Id}\hbox{ in }\overline{\Omega}\ \hbox{ and }\ \|\Lambda_j^{-1}\|_{L^1(\Omega)}=\|\Lambda^{-1}\|_{L^1(\Omega)}\ \hbox{ for all }j\in\N,\vspace{3pt}\\
0<\displaystyle\mathop{{\rm{ess}}\,{\rm{inf}}}_{\Omega}\Lambda\le\mathop{\liminf}_{j\to+\infty}\Big(\mathop{\min}_{\overline{\Omega}}\Lambda_j\Big)\le\mathop{\limsup}_{j\to+\infty}\Big(\mathop{\max}_{\overline{\Omega}}\Lambda_j\Big)\le\mathop{{\rm{ess}}\,{\rm{sup}}}_{\Omega}\Lambda.\eaa\right.
\ee
Due to~(\ref{Aklambdakbis}), one can assume without loss of generality that
\be\label{minlambdaj}
\min_{\overline{\Omega}}\Lambda_j\ge\frac{{\rm{ess}}\,{\rm{inf}}_{\Omega}\Lambda}{2}=:m_{\Lambda}>0
\ee
for all $j\in\N$. Since the sequence $(u_j)_{j\in\N}$ converges to $u$ in $W^{2,p}(\Omega)$ weak for all $1\le p<+\infty$ and in $C^{1,\alpha}(\overline{\Omega})$ for all $0\le\alpha<1$, and since $u\in E_{1/2,N_u,\beta_u}(\Omega)$ for some parameters $N_u>0$ and~$\beta_u>0$ depending on $u$ (from~(\ref{hypu}) and the smoothness of $u$), one can assume without loss of generality that~(\ref{ujNbeta}) holds for all $j\in\N$, that is
\be\label{ujNbetabis}
|\nabla u_j|\neq 0\hbox{ on }\partial\Omega,\ \ u_j>0\hbox{ in }\Omega\ \hbox{ and }\ u_j\in E_{1/2,2N_u,\beta_u/2}(\Omega).
\ee\par
Furthermore, if $\Omega$ is not a ball (in which case $n\ge 2$) and if~\eqref{hypAMquadr} holds, Lemma \ref{lemboundsquadr} yields that $u\in E_{1/2,N,\beta}(\Omega)$ for some positive constants $N>0$ and~$\beta>0$ only depending on $\Omega$, $n$, $q$, $M$ and $r$. One can therefore assume, without loss of generality, that, in that case,
\be\label{ujNbetaquadr}
|\nabla u_j|\neq 0\hbox{ on }\partial\Omega,\ \ u_j>0\hbox{ in }\Omega\ \hbox{ and }\ u_j\in E_{1/2,2N,\beta/2}(\Omega)
\ee
for all $j\in \N$.

\subsubsection*{Step 3: symmetrized coefficients and the inequalities $u^*_{j_k}\le\hat{\psi}_k$, $(1+\eta_u)u^*_{j_k}\le \hat{\psi}_k$ and~$(1+\eta)u^*_{j_k}\le\hat{\psi}_k$ in~$\Omega^*$}

Let now $k\in\N$ be fixed in this step and in the next two ones. For all $j\in\N$ and $x\in\Omega$, denote
$$\baa{rcl}
B_j(x) & = & \displaystyle-\,\hbox{div}(A_j\nabla u_j)(x)-a(x,u(x),\nabla u(x))\,|\nabla u_j(x)|^q+b(x,u(x),\nabla u(x))\,u_j(x)\vspace{3pt}\\
& & \displaystyle-f(x,u(x),\nabla u(x))-2^{-k}\vspace{3pt}\\
& = & -H_j(x)-a(x,u(x),\nabla u(x))\,|\nabla u_j(x)|^q+b(x,u(x),\nabla u(x))\,u_j(x)\vspace{3pt}\\
& & \displaystyle-f(x,u(x),\nabla u(x))-2^{-k}.\eaa$$
Due to~(\ref{hypH}), (\ref{Hkinfty}) and the fact that $u_j\to u$ in (at least) $C^1(\overline{\Omega})$ as $j\to+\infty$, it follows that
$$\limsup_{j\to+\infty}\Big(\sup_{x\in\Omega}B_j(x)\Big)\le-2^{-k}<0.$$
Therefore, there is an integer $j_k\ge k$ such that $B_{j_k}(x)\le0$ for all $x\in\Omega$, that is
\be\label{Bkepsbis}\baa{rcl}
\displaystyle-\,\hbox{div}(A_{j_k}\nabla u_{j_k})(x)-a(x,u(x),\nabla u(x))\,|\nabla u_{j_k}(x)|^q & & \vspace{3pt}\\
+\,b(x,u(x),\nabla u(x))\,u_{j_k}(x)-f(x,u(x),\nabla u(x)) & \le & 2^{-k}\ (\le 1)\ \hbox{ for all }x\in\Omega.\eaa
\ee\par
As in Step~3 of the proofs of Theorems~\ref{th1} and~\ref{th1bis}, one can then apply the general results of Section~\ref{rearrangement} to the coefficients $A_{\Omega}$, $\Lambda_{\Omega}$, $\psi$, $a_{\Omega}$, $b_{\Omega}$ and $f_{\Omega}$ as in~(\ref{defcoef}), that is
$$\left\{\baa{l}
A_{\Omega}(x)=A_{j_k}(x),\ \Lambda_{\Omega}(x)=\Lambda_{j_k}(x),\ \psi(x)=u_{j_k}(x),\vspace{3pt}\\
a_{\Omega}(x)=a(x,u(x),\nabla u(x)),\ b_{\Omega}(x)=b(x,u(x),\nabla u(x)),\vspace{3pt}\\
f_{\Omega}(x)=f(x,u(x),\nabla u(x))=f_u(x),\eaa\right.$$
and to our given power~$q\in(1,2]$. Call $\rho_k:[0,\max_{\overline{\Omega}}u_{j_k}]\to[0,R]$, $E_k$, $\hat{\Lambda}_k\in L^{\infty}_+(\Omega^*)$, $\hat{\psi}_k\in W^{1,\infty}(\Omega^*)\cap H^1_0(\Omega^*)$, $\hat{a}_k\in L^{\infty}(\Omega^*)$ and $\hat{f}_k\in L^{\infty}(\Omega^*)$ the symmetrized quantities defined as in~(\ref{defrho}), (\ref{defE}), (\ref{defhatlambda}), (\ref{defhatpsi}), (\ref{defhata}) and~(\ref{defhatf}). In particular, the set $E_k$ is given as~(\ref{Ek}) and $(0<)\,a^k_1<\cdots<a^k_{m_k}=\max_{\overline{\Omega}}u_{j_k}$ denote the $m_k$ critical values of the function~$u_{j_k}$ in~$\overline{\Omega}$. All functions $\hat{\Lambda}_k$, $\hat{\psi}_k$, $\hat{a}_k$ and~$\hat{f}_k$ are radially symmetric. It follows from~(\ref{ineqlambda}),~(\ref{ineqa}),~(\ref{ineqf}) and~(\ref{Aklambdakbis}) that
\be\label{coefepsilonbis}\left\{\baa{l}
\displaystyle0<\min_{\overline{\Omega}}\Lambda_{j_k}\le\mathop{{\rm{ess}}\,{\rm{inf}}}_{\Omega^*}\hat{\Lambda}_k\le\mathop{{\rm{ess}}\,{\rm{sup}}}_{\Omega^*}\hat{\Lambda}_k\le\max_{\overline{\Omega}}\Lambda_{j_k},\vspace{3pt}\\
\|\hat{\Lambda}_k^{-1}\|_{L^1(\Omega^*)}=\|\Lambda_{j_k}^{-1}\|_{L^1(\Omega)}=\|\Lambda^{-1}\|_{L^1(\Omega)},\vspace{3pt}\\
\displaystyle\inf_{\Omega\times\R\times\R^n}a^+\le\min_{\overline{\Omega}}a(\cdot,u(\cdot),\nabla u(\cdot))^+\le\mathop{{\rm{ess}}\,{\rm{inf}}}_{\Omega^*}\hat{a}_k\le\mathop{{\rm{ess}}\,{\rm{sup}}}_{\Omega^*}\hat{a}_k\le\cdots\\
\qquad\qquad\qquad\qquad\qquad\displaystyle\cdots\le\Big(\max_{\overline{\Omega}}a(\cdot,u(\cdot),\nabla u(\cdot))^+\Big)\times\Big(\frac{\max_{\overline{\Omega}}\Lambda_{j_k}}{\min_{\overline{\Omega}}\Lambda_{j_k}}\Big)^{q-1}\le\cdots\vspace{3pt}\\
\qquad\qquad\qquad\qquad\qquad\displaystyle\cdots\le\Big(\sup_{\Omega\times\R\times\R^n}a^+\Big)\times\Big(\frac{\max_{\overline{\Omega}}\Lambda_{j_k}}{\min_{\overline{\Omega}}\Lambda_{j_k}}\Big)^{q-1},\vspace{3pt}\\
\displaystyle\inf_{\Omega\times\R\times\R^n}f\le\min_{\overline{\Omega}}f_u\le\mathop{{\rm{ess}}\,{\rm{inf}}}_{\Omega^*}\hat{f}_k\le\mathop{{\rm{ess}}\,{\rm{sup}}}_{\Omega^*}\hat{f}_k\le\max_{\overline{\Omega}}f_u\le\sup_{\Omega\times\R\times\R^n}f,\ \ \displaystyle\int_{\Omega^*}\!\hat{f}_k=\int_{\Omega}\!f_u.\eaa\right.
\ee
Proposition~\ref{key1} then implies that~(\ref{ineqkey1}) and~(\ref{schwarzuk}) hold, that is
\be\label{schwarzuk3}
0\le u_{j_k}^*(x)\le\hat{\psi}_k(x)\ \hbox{ for all }x\in\overline{\Omega^*}.
\ee\par
Furthermore, if $\Omega$ is not a ball, it follows from~(\ref{ujNbetabis}) and Proposition~\ref{key1notball} that
\be\label{schwarzukbisuquadr}
0\le(1+\eta_u)\,u_{j_k}^*(x)\le\hat{\psi}_k(x)\ \hbox{ for all }x\in\overline{\Omega^*}.
\ee
where $\eta_u>0$ only depends on $\Omega$, $n$, $N_u$ and $\beta_u$, that is on $\Omega$, $n$ and $u$. If $\Omega$ is not a ball and the assumption~(\ref{hypAMquadr}) of Theorem~\ref{th2bis} is made, it follows from~(\ref{ujNbetaquadr}) and Proposition~\ref{key1notball} that
\be\label{schwarzukbisquadr}
0\le(1+\eta)\,u_{j_k}^*(x)\le\hat{\psi}_k(x)\ \hbox{ for all }x\in\overline{\Omega^*}
\ee
where $\eta>0$ only depends on $\Omega$, $n$, $N$ and $\beta$, that is on $\Omega$, $n$ and $M$.\par
Lastly, remember that the functions $u_{j_k}$ satisfy~(\ref{ujNbetabis}) for some positive constants $N_u$ and~$\beta_u$ independent of $k$, that the functions $\Lambda_{j_k}$ satisfy~(\ref{minlambdaj}), that
$$\min_{\overline{\Omega}}b_{\Omega}\ge\inf_{\Omega\times\R\times\R^n}b=:m_b>0$$
from~(\ref{hypb}), that $\|a_{\Omega}^+\|_{L^{\infty}(\Omega)}\le\|a\|_{L^{\infty}(\Omega\times\R\times\R^n)}$ and that $\|f_{\Omega}^+\|_{L^{\infty}(\Omega)}\le\|f\|_{L^{\infty}(\Omega\times\R\times\R^n)}$. Therefore, it follows from~(\ref{Bkepsbis}) and Proposition~\ref{key2bis} with $\kappa=1$ that there exists a constant
$$\hat{\delta}>0,$$
which depends on $\Omega$, $n$, $m_b$, $m_{\Lambda}$, $\|a\|_{L^{\infty}(\Omega\times\R\times\R^n)}$, $\|f\|_{L^{\infty}(\Omega\times\R\times\R^n)}$, $N_u$ and $\beta_u$, but which does not depend on $k$, such that, for all $x\in E_k\cap\Omega^*$, there exists a point~$y\in\Sigma_{\rho_k^{-1}(|x|)}$ satisfying
$$\begin{array}{rl}
& \displaystyle -\,{\rm{div}}\big(\widehat{\Lambda}_k\nabla\hat{\psi}_k\big)(x)-\widehat{a}_k(x)\vert \nabla\hat{\psi}_k(x)\vert^q+\hat{\delta}\,\hat{\psi}_k(x)-\widehat{f}_k(x)\vspace{3pt}\\
\leq\!\! & \displaystyle-\,{\rm{div}}(A_{j_k}\nabla u_{j_k})(y)-a(y,u(y),\nabla u(y))\left\vert \nabla u_{j_k}(y)\right\vert^q+b(y,u(y),\nabla u(y))\,u_{j_k}(y)\vspace{3pt}\\
& -f(y,u(y),\nabla u(y))\vspace{3pt}\\
\leq\!\! & 2^{-k}.\end{array}$$
In other words,
\be\label{ineqpsiepsbis}
-\,{\rm{div}}\big(\widehat{\Lambda}_k\nabla\hat{\psi}_k\big)(x)-\hat{a}_k(x)\,|\nabla\hat{\psi}_k(x)|^q+\hat{\delta}\,\hat{\psi}_k(x)\le g_k(x)\ \hbox{ for all }x\in E_k\cap\Omega^*,
\ee
where $g_k$ is defined as in~(\ref{defalphak}), that is
\be\label{defalphakbis}
g_k(x)=\widehat{f}_k(x)+2^{-k}\ \hbox{ for all }x\in E_k.
\ee
The constant $\hat{\delta}>0$ will be that of the conclusion of Theorem~\ref{th2}.

\subsubsection*{Step 4: the functions $\hat{\psi}_k$ are $H^1_0(\Omega^*)\cap L^{\infty}(\Omega^*)$ weak subsolutions of~$(\ref{ineqpsiepsbis})$}

We know that, for every $k\in\N$, the function $\hat{\psi}_k$ is in $W^{1,\infty}(\Omega^*)\cap H^1_0(\Omega^*)$, and that the inequality~(\ref{ineqpsiepsbis}) holds in $E_k\cap\Omega^*$, whence almost everywhere in $\Omega^*$. We now claim that $\hat{\psi}_k$ is a weak $H^1_0(\Omega^*)\cap L^{\infty}(\Omega^*)$ subsolution of~(\ref{ineqpsiepsbis}), in the sense that
\be\label{weaksub}
\int_{\Omega^*}\hat{\Lambda}_k\nabla\hat{\psi}_k\cdot\nabla\varphi-\int_{\Omega^*}\hat{a}_k|\nabla\hat{\psi}_k|^q\,\varphi+\int_{\Omega^*}\hat{\delta}\,\hat{\psi}_k\,\varphi-\int_{\Omega^*}g_k\,\varphi\le0
\ee
for all $\varphi\in H^1_0(\Omega^*)\cap L^{\infty}(\Omega^*)$ such that $\varphi\ge0$ a.e. in $\Omega^*$.\par
The proof of this claim follows exactly the same scheme as that of Lemma~\ref{lemsubpsik} given in Section~\ref{sec32}. Indeed, one first notices that, given any test function $\varphi\in H^1_0(\Omega^*)\cap L^{\infty}(\Omega^*)$ with~$\varphi\ge0$ a.e. in $\Omega^*$, there is a sequence $(\varphi_m)_{m\in\N}$ of nonnegative $C^1_c(\Omega^*)$ functions such that~$\varphi_m\to\varphi$ as $m\to+\infty$ in $H^1_0(\Omega^*)$ (but the convergence does not hold in $L^{\infty}(\Omega^*)$ in general). Next, since $\hat{\Lambda}_k\in L^{\infty}_+(\Omega^*)$, $\hat{\psi}_k\in W^{1,\infty}(\Omega^*)$, $\hat{a}_k\in L^{\infty}(\Omega^*)$ and $g_k\in L^{\infty}(\Omega^*)$, it follows that the left-hand side of~(\ref{weaksub}) is equal to the limit as $m\to+\infty$ of the same quantities with~$\varphi_m$ instead of $\varphi$. It is therefore sufficient to show~(\ref{weaksub}) when the test function $\varphi$ belongs to $C^1_c(\Omega^*)$ and is nonnegative. Finally, one can repeat the proof of Lemma~\ref{lemsubpsik} and one just needs to replace the quantities $\hat{a}_k\,e_r\cdot\nabla\hat{\psi}_k$ and $-g_k$ by, respectively, $-\hat{a}_k|\nabla\hat{\psi}_k|^q$ ($\in L^{\infty}(\Omega^*)$) and $\hat{\delta}\,\hat{\psi}_k-g_k$ ($\in L^{\infty}(\Omega^*)$), without any other modification in the proof.\par
As a conclusion, the claim~(\ref{weaksub}) holds.

\subsubsection*{Step 5: the inequalities $u^*_{j_k}\le v_k$, $(1+\eta_u)u^*_{j_k}\le v_k$ and $(1+\eta)u^*_{j_k}\le v_k$ in $\Omega^*$}

For every $k\in\N$, let $v_k$ be the unique weak $H^1_0(\Omega^*)\cap L^{\infty}(\Omega^*)$ solution of
\be\label{defvkbis}\left\{\baa{rcll}
-\,{\rm{div}}\big(\widehat{\Lambda}_k\nabla v_k\big)-\hat{a}_k|\nabla v_k|^q+\hat{\delta}\,v_k & = & g_k & \hbox{in }\Omega^*,\vspace{3pt}\\
v_k & = & 0 & \hbox{on }\partial\Omega^*,\eaa\right.
\ee
in the sense that
$$\int_{\Omega^*}\hat{\Lambda}_k\nabla v_k\cdot\nabla\varphi-\int_{\Omega^*}\hat{a}_k|\nabla v_k|^q\,\varphi+\int_{\Omega^*}\hat{\delta}\,v_k\,\varphi=\int_{\Omega^*}g_k\,\varphi$$
for all $\varphi\in H^1_0(\Omega^*)\cap L^{\infty}(\Omega^*)$. As recalled in Section~\ref{intro} after the proof of Corollary~\ref{cor2}, the existence and uniqueness of $v_k\in H^1_0(\Omega^*)\cap L^{\infty}(\Omega^*)$ is guaranteed by Th\'eor\`eme~2.1 and the following comments of~\cite{bmp84} and by Theorem~2.1 of~\cite{bm95}, since $\hat{\Lambda}_k\in L^{\infty}_+(\Omega^*)$, $\hat{a}_k\in L^{\infty}(\Omega^*)$,~$g_k\in L^{\infty}(\Omega^*)$, $\hat{\delta}>0$ and $q\in(1,2]$. Furthermore, the function $v_k$ is radially symmetric in $\Omega^*$ by the uniqueness, since all coefficients of~(\ref{defvkbis}) are so, and, without loss of generality, $v_k$ is continuous in $\overline{\Omega^*}$ from the local continuity (Corollary~4.23 of~\cite{hl}), the radial symmetry and the fact that $v_k\in H^1_0(\Omega^*)$. Lastly, it follows from~(\ref{weaksub}) and Theorem~2.1 of~\cite{bm95} that
$$\hat{\psi}_k\le v_k\hbox{ a.e. in }\Omega^*$$
(actually, the inequality can be assumed to hold everywhere in $\overline{\Omega^*}$ since both functions $\hat{\psi}_k$ and $v_k$ can be assumed to be continuous in $\overline{\Omega^*}$ without loss of generality). One concludes from~(\ref{schwarzuk3}) that
\be\label{inequjkvkter}
0\le u_{j_k}^*\le v_k\ \hbox{ a.e. in }\Omega^*.
\ee\par
If $\Omega$ is not a ball, inequality \eqref{schwarzukbisuquadr} therefore yields
\be\label{inequjkvkteruquadr}
0\le(1+\eta_u)\,u_{j_k}^*(x)\le v_k\ \hbox{ a.e. in }\Omega^*.
\ee
Furthermore, if $\Omega$ is not a ball and the assumption~(\ref{hypAMquadr}) of Theorem~\ref{th2bis} is made, it follows from \eqref{schwarzukbisquadr} that
\be\label{inequjkvkterquadr}
0\le(1+\eta)\,u_{j_k}^*(x)\le v_k\ \hbox{ a.e. in }\Omega^*.
\ee

\subsubsection*{Step 6: the limiting inequalities $u^*\le v$, $(1+\eta_u)u^*\le v$ and $(1+\eta)u^*\le v$ in $\Omega^*$}

As in the beginning of Step~6 of the proofs of Theorems~\ref{th1} and~\ref{th1bis}, one can assume that, up to extraction of a subsequence,
\be\label{ujk*ter}
u_{j_k}^*(x)\to u^*(x)\hbox{ a.e. in }\Omega^*\hbox{ as }k\to+\infty.
\ee
On the other hand, the sequences $(\hat{\Lambda}_k)_{k\in\N}$, $(\hat{\Lambda}_k^{-1})_{k\in\N}$, $(\hat{a}_k)_{k\in\N}$ and $(g_k)_{k\in\N}$ are bounded in~$L^{\infty}(\Omega^*)$ from~(\ref{Aklambdakbis}),~(\ref{coefepsilonbis}) and~(\ref{defalphakbis}). Furthermore, $\hat{\delta}$ is positive. It follows then from Theorem~2.1 of~\cite{bmp88cras} and Theorem~1 of~\cite{bmp92} (see also Theorem~3.1 of~\cite{Porretta}) that the sequence~$(v_k)_{k\in\N}$ is bounded in $H^1_0(\Omega^*)\cap L^{\infty}(\Omega^*)$ and relatively compact in~$H^1_0(\Omega^*)$. Therefore, there exists a radially symmetric function $v\in H^1_0(\Omega^*)\cap L^{\infty}(\Omega^*)$ such that, up to extraction of a subsequence,
\be\label{vkvter}
v_k\to v\hbox{ in }H^1_0(\Omega^*)\hbox{ strong and a.e. in }\Omega^*\hbox{ as }k\to+\infty.
\ee
Hence, together with~(\ref{inequjkvkter}) and~(\ref{ujk*ter}), one infers~(\ref{u*vinfty}) with $v$ instead of $v_{\infty}$, that is
\be\label{u*vinftyter}
0\le u^*\le v\ \hbox{ a.e. in }\Omega^*.
\ee\par
If $\Omega$ is not a ball, inequality \eqref{inequjkvkteruquadr} yields
\be\label{u*vuquadr}
0\le(1+\eta_u)\,u^*(x)\le v\ \hbox{ a.e. in }\Omega^*.
\ee
Furthermore, if $\Omega$ is not a ball and the assumption~(\ref{hypAMquadr}) of Theorem~\ref{th2bis} is made, it follows from \eqref{inequjkvkterquadr} that
\be\label{u*vquadr}
0\le(1+\eta)\,u^*(x)\le v\ \hbox{ a.e. in }\Omega^*.
\ee\par
The function $v$ will be that of the conclusion of Theorem~\ref{th2}. One shall identify in the following steps the equation~(\ref{eqv4}) satisfied by $v$ and one also gets the inequality~(\ref{u*vepsilon}) involving the solution~$v_{\epsilon}$ of~(\ref{eqv4eps}).

\subsubsection*{Step 7: a limiting equation satisfied by~$v$ in $\Omega^*$}

Let us now pass to the limit in the coefficients $\hat{\Lambda}_k$, $\hat{a}_k$ and $g_k$ of~(\ref{defvkbis}). From~(\ref{Aklambdakbis}),~(\ref{coefepsilonbis}) and~(\ref{defalphakbis}), there exist some radially symmetric functions $\hat{\Lambda}\in L^{\infty}_+(\Omega^*)$, $\hat{a}\in L^{\infty}(\Omega^*)$ and~$\hat{f}\in L^{\infty}(\Omega^*)$ such that, up to extraction of some subsequence,
\be\label{chapeauxbis}
\hat{\Lambda}_k^{-1}\rightharpoonup\hat{\Lambda}^{-1},\ \ \hat{\Lambda}_k^{-q}\hat{a}_k\rightharpoonup\hat{\Lambda}^{-q}\hat{a}\ \hbox{ and }\ g_k\rightharpoonup\hat{f}\ \hbox{ in }L^{\infty}(\Omega^*)\hbox{ weak-* as }k\to+\infty,
\ee
whence
\be\label{hats2bis}
\displaystyle0<\mathop{{\rm{ess}}\,{\rm{inf}}}_{\Omega}\Lambda\le\mathop{{\rm{ess}}\,{\rm{inf}}}_{\Omega^*}\hat{\Lambda}\le\mathop{{\rm{ess}}\,{\rm{sup}}}_{\Omega^*}\hat{\Lambda}\le\mathop{{\rm{ess}}\,{\rm{sup}}}_{\Omega}\Lambda\ \hbox{ and }\ \|\hat{\Lambda}^{-1}\|_{L^1(\Omega^*)}=\|\Lambda^{-1}\|_{L^1(\Omega)}.
\ee
Furthermore,
\be\label{hats2ter}\left\{\baa{l}
\displaystyle\mathop{\inf}_{\Omega\times\R\times\R^n}f\le\mathop{{\rm{ess}}\,{\rm{inf}}}_{\Omega^*}\hat{f}\le\mathop{{\rm{ess}}\,{\rm{sup}}}_{\Omega^*}\hat{f}\le\mathop{\sup}_{\Omega\times\R\times\R^n}f,\vspace{3pt}\\
\displaystyle\int_{\Omega^*}\hat{f}=\displaystyle{\mathop{\lim}_{k\to+\infty}}\int_{\Omega^*}g_k=\displaystyle{\mathop{\lim}_{k\to+\infty}}\int_{\Omega^*}\hat{f}_k=\int_{\Omega}f_u.\eaa\right.
\ee
In~(\ref{chapeauxbis}), the function $\hat{a}=\hat{\Lambda}^q\,\hat{\Lambda}^{-q}\,\hat{a}$ is defined as $\hat{\Lambda}^q$ times the $L^{\infty}(\Omega^*)$ weak-* limit of the functions $\hat{\Lambda}_k^{-q}\hat{a}_k$.\par
The functions $\hat{\Lambda}$, $\hat{a}$ and $\hat{f}$ will be those of the conclusion of Theorem~\ref{th2}. Notice first that, from~(\ref{hats2bis}) and~(\ref{hats2ter}), the functions $\hat{\Lambda}$ and $\hat{f}$ fulfill~(\ref{hats4}). Let us now establish in this paragraph the bounds~(\ref{hats4}) for the function~$\hat{a}$. It follows from~(\ref{coefepsilonbis}) that
$$0\le\inf_{\Omega\times\R\times\R^n}a^+\le\mathop{{\rm{ess}}\,{\rm{inf}}}_{\Omega^*}\hat{a}_k\le\mathop{{\rm{ess}}\,{\rm{sup}}}_{\Omega^*}\hat{a}_k\le\Big(\sup_{\Omega\times\R\times\R^n}a^+\Big)\times\Big(\frac{\max_{\overline{\Omega}}\Lambda_{j_k}}{\min_{\overline{\Omega}}\Lambda_{j_k}}\Big)^{q-1}$$
and
$$\mathop{{\rm{ess}}\,{\rm{sup}}}_{\Omega^*}\big(\hat{\Lambda}_k^{-q+1}\hat{a}_k\big)\le\Big(\frac{1}{\min_{\overline{\Omega}}\Lambda_{j_k}}\Big)^{q-1}\times\Big(\sup_{\Omega\times\R\times\R^n}a^+\Big)\times\Big(\frac{\max_{\overline{\Omega}}\Lambda_{j_k}}{\min_{\overline{\Omega}}\Lambda_{j_k}}\Big)^{q-1}.$$
Therefore,
$$\limsup_{k\to+\infty}\Big(\mathop{{\rm{ess}}\,{\rm{sup}}}_{\Omega^*}\big(\hat{\Lambda}_k^{-q+1}\hat{a}_k\big)\Big)\le\frac{\big({\rm{ess}}\,{\rm{sup}}_{\Omega}\Lambda\big)^{q-1}}{\big({\rm{ess}}\,{\rm{inf}}_{\Omega}\Lambda\big)^{2(q-1)}}\times\sup_{\Omega\times\R\times\R^n}a^+$$
from~(\ref{Aklambdakbis}). Since, by~(\ref{chapeauxbis}),
\be\label{weaklimits}\left\{\baa{l}
\hat{\Lambda}^q\hat{\Lambda}_k^{-1}\rightharpoonup\hat{\Lambda}^{q-1}\vspace{3pt}\\
\big(\hat{\Lambda}^q\hat{\Lambda}_k^{-1}\big)\times\big(\hat{\Lambda}_k^{-q+1}\hat{a}_k\big)=\hat{\Lambda}^q\times\big(\hat{\Lambda}_k^{-q}\hat{a}_k\big)\rightharpoonup\hat{\Lambda}^q\times\big(\hat{\Lambda}^{-q}\hat{a}\big)=\hat{a}\eaa\right.\hbox{ in }L^{\infty}(\Omega^*)\hbox{ weak-*}
\ee
as $k\to+\infty$, one gets that
$$\hat{a}\le\hat{\Lambda}^{q-1}\times\frac{\big({\rm{ess}}\,{\rm{sup}}_{\Omega}\Lambda\big)^{q-1}}{\big({\rm{ess}}\,{\rm{inf}}_{\Omega}\Lambda\big)^{2(q-1)}}\times\sup_{\Omega\times\R\times\R^n}a^+\le\Big(\frac{{\rm{ess}}\,{\rm{sup}}_{\Omega}\Lambda}{{\rm{ess}}\,{\rm{inf}}_{\Omega}\Lambda}\Big)^{2(q-1)}\times\sup_{\Omega\times\R\times\R^n}a^+\ \hbox{ a.e. in }\Omega^*,$$
where the last inequality follows from~(\ref{hats2bis}). On the other hand, it follows again from~(\ref{coefepsilonbis}) that
$$\mathop{{\rm{ess}}\,{\rm{inf}}}_{\Omega^*}\big(\hat{\Lambda}_k^{-q+1}\hat{a}_k\big)\ge\Big(\frac{1}{\max_{\overline{\Omega}}\Lambda_{j_k}}\Big)^{q-1}\times\inf_{\Omega\times\R\times\R^n}a^+,$$
whence
$$\liminf_{k\to+\infty}\Big(\mathop{{\rm{ess}}\,{\rm{inf}}}_{\Omega^*}\big(\hat{\Lambda}_k^{-q+1}\hat{a}_k\big)\Big)\ge\Big(\frac{1}{{\rm{ess}}\,{\rm{sup}}_{\Omega}\Lambda}\Big)^{q-1}\times\inf_{\Omega\times\R\times\R^n}a^+$$
from~(\ref{Aklambdakbis}). Using again~(\ref{hats2bis}) and~(\ref{weaklimits}), one gets that
$$\hat{a}\ge\hat{\Lambda}^{q-1}\times\Big(\frac{1}{{\rm{ess}}\,{\rm{sup}}_{\Omega}\Lambda}\Big)^{q-1}\times\inf_{\Omega\times\R\times\R^n}a^+\ge\Big(\frac{{\rm{ess}}\,{\rm{inf}}_{\Omega}\Lambda}{{\rm{ess}}\,{\rm{sup}}_{\Omega}\Lambda}\Big)^{q-1}\times\inf_{\Omega\times\R\times\R^n}a^+\,\ge\,0\ \hbox{ a.e. in }\Omega^*.$$
Finally, the radially symmetric function $\hat{a}\in L^{\infty}(\Omega^*)$ satisfies the properties~(\ref{hats4}) of the statement of Theorem~\ref{th2}.\par
The main goal of this step is to show that $v$ is a weak $H^1_0(\Omega^*)\cap L^{\infty}(\Omega^*)$ solution of the limiting equation
\be\label{eqvlimitbis}\left\{\baa{rcll}
-\,{\rm{div}}\big(\widehat{\Lambda}\nabla v\big)-\hat{a}\,|\nabla v|^q+\hat{\delta}\,v & = & \hat{f} & \hbox{in }\Omega^*,\vspace{3pt}\\
v & = & 0 & \hbox{on }\partial\Omega^*\eaa\right.
\ee
obtained by passing formally to the limit as $k\to+\infty$ in~(\ref{defvkbis}), in the sense that
$$\int_{\Omega^*}\hat{\Lambda}\nabla v\cdot\nabla\varphi-\int_{\Omega^*}\hat{a}\,|\nabla v|^q\,\varphi+\int_{\Omega^*}\hat{\delta}\,v\,\varphi-\int_{\Omega^*}\hat{f}\,\varphi=0$$
for all $\varphi\in H^1_0(\Omega^*)\cap L^{\infty}(\Omega^*)$ (notice that all integrals in the above formula converge since $v$ and $\varphi$ belong to $H^1_0(\Omega^*)\cap L^{\infty}(\Omega^*)$, $q\in(1,2]$ and $\hat{\Lambda}$, $\hat{a}$ and $\hat{f}$ are in $L^{\infty}(\Omega^*)$).\par
The proof of~(\ref{eqvlimitbis}) will actually be a consequence of the following lemma.

\begin{lem}\label{lem41}
Let $(\delta_k)_{k\in\N}$ be a sequence of positive real numbers, let $(\lambda_k)_{k\in\N}$ be a sequence of radially symmetric $L^{\infty}_+(\Omega^*)$ functions and let $(\alpha_k)_{k\in\N}$ and $(\gamma_k)_{k\in\N}$ be two sequences of radially symmetric $L^{\infty}(\Omega^*)$ functions. Assume that the sequences $(\lambda_k)_{k\in\N}$, $(\lambda_k^{-1})_{k\in\N}$, $(\alpha_k)_{k\in\N}$ and~$(\gamma_k)_{k\in\N}$ are bounded in $L^{\infty}(\Omega^*)$ and that there are $\delta_{\infty}\in(0,+\infty)$ and some functions $\lambda_{\infty}\in L^{\infty}_+(\Omega^*)$, $\alpha_{\infty}\in L^{\infty}(\Omega^*)$ and $\gamma_{\infty}\in L^{\infty}(\Omega^*)$ such that $\delta_k\to\delta_{\infty}$ as $k\to+\infty$ and
\be\label{chapeauxter}
\lambda_k^{-1}\rightharpoonup\lambda_{\infty}^{-1},\ \ \lambda_k^{-q}\alpha_k\rightharpoonup\lambda_{\infty}^{-q}\alpha_{\infty}\ \hbox{ and }\ \gamma_k\rightharpoonup\gamma_{\infty}\ \hbox{ in }L^{\infty}(\Omega^*)\hbox{ weak-* }\hbox{ as }k\to+\infty.
\ee
Let $({V}_k)_{k\in\N}$ be the sequence of weak $H^1_0(\Omega^*)\cap L^{\infty}(\Omega^*)$ solutions of
\be\label{defvkter}\left\{\baa{rcll}
-\,{\rm{div}}\big(\lambda_k\nabla{V}_k\big)-\alpha_k|\nabla{V}_k|^q+\delta_k\,{V}_k & = & \gamma_k & \hbox{in }\Omega^*,\vspace{3pt}\\
{V}_k & = & 0 & \hbox{on }\partial\Omega^*.\eaa\right.
\ee
Then there is a radially symmetric function ${V}_{\infty}\in H^1_0(\Omega^*)\cap L^{\infty}(\Omega^*)$ such that
\be\label{weakVk}
{V}_k\to {V}_{\infty}\hbox{ in }H^1_0(\Omega^*)\hbox{ strong as }k\to+\infty,
\ee
where ${V}_{\infty}$ denotes the unique weak $H^1_0(\Omega^*)\cap L^{\infty}(\Omega^*)$ solution of
\be\label{Vinfty}\left\{\baa{rcll}
-\,{\rm{div}}\big(\lambda_{\infty}\nabla{V}_{\infty}\big)-\alpha_{\infty}\,|\nabla{V}_{\infty}|^q+\delta_{\infty}\,{V}_{\infty} & = & \gamma_{\infty} & \hbox{in }\Omega^*,\vspace{3pt}\\
{V}_{\infty} & = & 0 & \hbox{on }\partial\Omega^*,\eaa\right.
\ee
in the sense that
\be\label{defIvarphi}
I^{\varphi}:=\int_{\Omega^*}\lambda_{\infty}\nabla V_{\infty}\cdot\nabla\varphi-\int_{\Omega^*}\alpha_{\infty}\,|\nabla V_{\infty}|^q\,\varphi+\int_{\Omega^*}\delta_{\infty}\,V_{\infty}\,\varphi-\int_{\Omega^*}\gamma_{\infty}\,\varphi=0
\ee
for all $\varphi\in H^1_0(\Omega^*)\cap L^{\infty}(\Omega^*)$.
\end{lem}

In order to finish the proof of Theorem~\ref{th2}, the proof of Lemma~\ref{lem41} is postponed in Section~\ref{sec42}. Notice that, together with~(\ref{u*vinftyter}), Lemma~\ref{lem41} already provides the first part of the conclusion of Theorem~\ref{th2}, that is~\eqref{u*vbis} with $v$ satisfying~\eqref{eqv4} and~\eqref{defhatH}. It only remains to show the comparison~(\ref{u*vepsilon}) of $u^*$ with the solution $v_{\epsilon}$ of the equation~(\ref{eqv4eps}) involving a function $\hat{f}_{\epsilon}$ having the same distribution function as the function $f_u$ defined in~(\ref{deffu}). To do so, we follow the same scheme as in last steps of the proofs of Theorems~\ref{th1} and~\ref{th1bis} in Section~\ref{sec31}. Namely, having in hand that $0\le u^*\le v$ (or the inequalities $0\leq (1+\eta_u)u^*\leq v$ and $0\leq (1+\eta)u^*\leq v$) a.e. in~$\Omega^*$ where $v$ is the weak~$H^1_0(\Omega^*)\cap L^{\infty}(\Omega^*)$ solution of equation~(\ref{eqvlimitbis}), we shall approximate~$v$ by the solutions of some approximating equations where $\hat{f}$ in~(\ref{eqvlimitbis}) is replaced by some right-hand sides having the same distribution function as the function $f_u$.

\subsubsection*{Step 8: approximation of $v$ by some functions $w_k$ in $\Omega^*$}

Let $(\hat{f}_k)_{k\in\N}$ be the sequence of radially symmetric functions defined in Step~3, and remember that the sequence $(\hat{f}_k)_{k\in\N}$ is bounded in $L^{\infty}(\Omega^*)$ from~(\ref{coefepsilonbis}). Since $\hat{\Lambda}\in L^{\infty}_+(\Omega^*)$, $\hat{a}\in L^{\infty}(\Omega^*)$ and $\hat{\delta}>0$, it follows from~\cite{bm95,bmp84} that, for each $k\in\N$, there is a unique weak $H^1_0(\Omega^*)\cap L^{\infty}(\Omega^*)$ solution~$w_k$ of
\be\label{eqwkbis}\left\{\baa{rcll}
-\,{\rm{div}}\big(\widehat{\Lambda}\nabla w_k\big)-\hat{a}\,|\nabla w_k|^q+\hat{\delta}\,w_k & = & \hat{f}_k & \hbox{in }\Omega^*,\vspace{3pt}\\
w_k & = & 0 & \hbox{on }\partial\Omega^*,\eaa\right.
\ee
in the sense that
$$\int_{\Omega^*}\hat{\Lambda}\,\nabla w_k\cdot\nabla\varphi-\int_{\Omega^*}\hat{a}\,|\nabla w_k|^q\,\varphi+\int_{\Omega^*}\hat{\delta}\,w_k\,\varphi-\int_{\Omega^*}\hat{f}_k\,\varphi=0$$
for all $\varphi\in H^1_0(\Omega^*)\cap L^{\infty}(\Omega^*)$. Furthermore, the functions $w_k$ are all radially symmetric by uniqueness and since all coefficients $\hat{\Lambda}$, $\hat{a}$ and $\hat{f}_k$ are radially symmetric. Lastly, as in Step~6 above, the sequence $(w_k)_{k\in\N}$ is bounded in $H^1_0(\Omega^*)\cap L^{\infty}(\Omega^*)$ and relatively compact in~$H^1_0(\Omega^*)$, from~\cite{bmp88cras,bmp92}. There exists then a function $w_{\infty}\in H^1_0(\Omega^*)\cap L^{\infty}(\Omega^*)$ such that, up to extraction of a subsequence, $w_k\to w_{\infty}$ in $H^1_0(\Omega^*)$ strong as $k\to+\infty$. Since $\hat{f}_k=g_k-2^{-k}\rightharpoonup\hat{f}$ in $L^{\infty}(\Omega^*)$ weak-* as $k\to+\infty$ from~(\ref{defalphakbis}) and~(\ref{chapeauxbis}), it follows from Lemma~\ref{lem41} applied with $\delta_k=\delta_{\infty}=\hat{\delta}$, $\lambda_k=\lambda_{\infty}=\hat{\Lambda}$, $\alpha_k=\alpha_{\infty}=\hat{a}$, $\gamma_k=\hat{f}_k$, $\gamma_{\infty}=\hat{f}$, $w_k=V_k$, that $w_{\infty}$ is the unique weak~$H^1_0(\Omega^*)\cap L^{\infty}(\Omega^*)$ solution of the limiting equation~(\ref{eqvlimitbis}), that is $w_{\infty}=v$ by uniqueness. By uniqueness of the limit, the whole sequence $(w_k)_{k\in\N}$ converges to $v$, that is
\be\label{wkvinftybis}
w_k\to v\hbox{ in }H^1_0(\Omega^*)\hbox{ strong as }k\to+\infty.
\ee

\subsubsection*{Step 9: approximation of the function $w_K$ for $K$ large by some functions $z_l$ in $\Omega^*$}

Let $\epsilon>0$ be an arbitrary positive real number. From~(\ref{wkvinftybis}) and Sobolev embeddings, there is an integer $K\in\N$ large enough such that
\be\label{wKvbis}
\|w_K-v\|_{L^{2^*}(\Omega^*)}\le\frac{\epsilon}{2},
\ee
where $w_K$ is the unique weak $H^1_0(\Omega^*)\cap L^{\infty}(\Omega^*)$ solution of~(\ref{eqwkbis}) with $k=K$ and where the Sobolev exponant $2^*$ is defined as in Theorem~\ref{th2}.\par
Now, as in Step~9 of the proofs of Theorems~\ref{th1} and~\ref{th1bis}, there is a sequence~$(h_l)_{l\in\N}$ of~$L^{\infty}(\Omega^*)$ radially symmetric functions such that
$$h_l\rightharpoonup\hat{f}_K\hbox{ as }l\to+\infty\hbox{ in }L^{\infty}(\Omega^*)\hbox{ weak-* and }\mu_{h_l}=\mu_{f_u}\hbox{ for all }l\in\N.$$
As in the previous step of the proof of the present Theorems~\ref{th2} and~\ref{th2bis}, for every $l\in\N$, there is a unique weak $H^1_0(\Omega^*)\cap L^{\infty}(\Omega^*)$, and radially symmetric, solution $z_l$ of
\be\label{eqzlbis}\left\{\baa{rcll}
-\,{\rm{div}}\big(\widehat{\Lambda}\nabla z_l\big)-\hat{a}\,|\nabla z_l|^q+\hat{\delta}\,z_l & = & h_l & \hbox{in }\Omega^*,\vspace{3pt}\\
z_l & = & 0 & \hbox{on }\partial\Omega^*\eaa\right.
\ee
and, from Lemma~\ref{lem41} again, the functions $z_l$ converge to $w_K$ in $H^1_0(\Omega^*)$ strong as $l\to+\infty$. In particular, there is $L\in\N$ large enough such that $\|z_L-w_K\|_{L^{2^*}(\Omega^*)}\le\epsilon/2$, whence
\be\label{zLvinftybis}
\|z_L-v\|_{L^{2^*}(\Omega^*)}\le\epsilon
\ee
from~(\ref{wKvbis}).

\subsubsection*{Step 10: conclusion of the proofs of Theorems~\ref{th2} and~\ref{th2bis}}

Call
$$v_{\epsilon}=z_L\ \hbox{ and }\ \hat{f}_{\epsilon}=h_L.$$
The function $\hat{f}_{\epsilon}$ is radially symmetric and has the same distribution function as the function~$f_u$. Furthermore, $v_{\epsilon}$ solves~(\ref{eqv4eps}). Lastly, it follows from~(\ref{u*vinftyter}) and~(\ref{zLvinftybis}) that
$$\baa{rcl}
\|(u^*-v_{\epsilon})^+\|_{L^{2^*}(\Omega^*)} & \le & \|(u^*-v)^++(v-v_{\epsilon})^+\|_{L^{2^*}(\Omega^*)}\vspace{3pt}\\
& \le & \|(u^*-v)^+\|_{L^{2^*}(\Omega^*)}+\|(v-v_{\epsilon})^+\|_{L^{2^*}(\Omega^*)}\vspace{3pt}\\
& \le & 0+\epsilon=\epsilon.\eaa$$
This is the desired conclusion~(\ref{u*vepsilon}). Under the assumptions of Theorem \ref{th2bis}, one obtains similarly
$$
\|((1+\eta_u)u^*-v_{\epsilon})^+\|_{L^{2^*}(\Omega^*)}\le \epsilon
$$
when $\Omega$ is not a ball, and
$$
\|((1+\eta)u^*-v_{\epsilon})^+\|_{L^{2^*}(\Omega^*)}\le \epsilon
$$
when $\Omega$ is not a ball and \eqref{hypAMquadr} is assumed. The proofs of Theorems~\ref{th2} and \ref{th2bis} are thereby complete.\hfill\fin


\subsection{Proofs of Lemmas~\ref{lemboundsquadr} and~\ref{lem41}}\label{sec42}

This section is devoted to the proofs of the technical lemmas~\ref{lemboundsquadr} and~\ref{lem41} which were stated and used in Steps~1,~7,~8 and~9 of the proofs of Theorems~\ref{th2} and~\ref{th2bis}.\hfill\break

\noindent{\bf Proof of Lemma \ref{lemboundsquadr}.} The proof follows the same lines as the one of Lemma \ref{lembounds}.\par
{\it Step 1: a uniform bound on $\left\Vert u\right\Vert_{L^{\infty}(\Omega)}$. }ÊBy the conditions on $q$ and $r$ contained in assumption~\eqref{hypAMquadr} and Theorem \ref{thapp2} below, it follows that $u\in W(\Omega)$, and in particular $u$ is qualitatively bounded. Furthermore, $u\in H^1_0(\Omega)$ and \eqref{equ}, \eqref{hypH} and \eqref{hypAMquadr} show that
$$
-\mbox{div}(A\nabla u)+M^{-1}u\leq M\left(\left\vert \nabla u\right\vert^q+1\right)
$$
in the weak $H^1_0(\Omega)\cap L^{\infty}(\Omega)$ sense. Thus, Theorem 3.1 in \cite{Porretta} ensures that 
\begin{equation} \label{ubounded}
\left\Vert u\right\Vert_{L^{\infty}(\Omega)}\leq C_1
\end{equation}
where $C_1=M^2$.\par
{\it Step 2: a uniform bound on $\left\Vert u\right\Vert_{H^1_0(\Omega)}$.} Arguing as in the proof of Lemma \ref{lembounds}, one has
$$
M^{-1}\int_{\Omega} \left\vert \nabla u(x)\right\vert^2dx\leq \int_{\Omega} -H(x,u(x),\nabla u(x))u(x)dx.
$$
Since
$$
\left\vert H(x,u(x),\nabla u(x))\right\vert \leq M\left(1+\left\vert u(x)\right\vert^r+\left\vert \nabla u(x)\right\vert^q\right)\leq M(1+C_1^r)+M\left\vert \nabla u(x)\right\vert^q,
$$
using the fact that $q<1+2/n\le 2$ (recall that $n\geq 2$), one derives
$$M^{-1}\int_{\Omega}|\nabla u|^2\le M\,C_1\int_{\Omega}\big(1+C_1^r+|\nabla u|^q\big),$$
whence
$$
\left\Vert u\right\Vert_{H^1_0(\Omega)}\leq C_2,
$$
where $C_2>0$ only depends on $\Omega$, $q$, $M$ and $r$.\par
{\it Step $3$: a uniform estimate of $\left\Vert u\right\Vert_{W^{2,p}(\Omega)}$. }ÊSince $q<2$, the proof of Theorem \ref{thapp2} below and the quantitative assumption~\eqref{hypAMquadr} show that, for all $1\leq p<\infty$, 
\be\label{uw2p}
\left\Vert u\right\Vert_{W^{2,p}(\Omega)}\leq K_p,
\ee
where $K_p>0$ only depends on $\Omega$, $n$, $q$, $M$, $r$ and $p$. In particular,
$$
\left\Vert u\right\Vert_{C^{1,1/2}(\overline{\Omega})}\leq C_3,
$$
where $C_3>0$ only depends on $\Omega$, $n$, $q$, $M$ and $r$.\par
{\it Step 4: conclusion.} What remains to be proved is the existence of a positive constant~$\beta$, only depending on $\Omega$, $n$, $q$, $M$ and $r$ such that
$$u(x)\ge\beta\,d(x,\partial\Omega)\hbox{ for all }x\in\Omega.$$
As in the proof of Lemma \ref{lembounds}, one argues by contradiction, assuming that this conclusion does not hold. Then there are a sequence~$(A_m)_{m\in\N}$ of $W^{1,\infty}(\Omega,{\mathcal{S}}_n(\R))$ matrix fields, a sequence~$(\Lambda_m)_{m\in\N}$ of $L^{\infty}_+(\Omega)$ functions and four sequences~$(H_m)_{m\in\N}$, $(a_m)_{m\in\N}$, $(b_m)_{m\in\N}$, $(f_m)_{m\in\N}$ of continuous functions in~$\overline{\Omega}\times\R\times\R^n$, satis\-fying~(\ref{ALambda}),~(\ref{hypH}) and~(\ref{hypAMquadr}) with the same parameter $M>0$, as well as a sequence~$(u_m)_{m\in\N}$ of~$W(\Omega)$ solutions of~(\ref{equ}) satisfying~(\ref{hypu}), and a sequence~$(x_m)_{m\in\N}$ of points in $\Omega$ such that~\eqref{xm} holds. Observe first that the sequence~$(H_m(\cdot,u_m(\cdot),\nabla u_m(\cdot)))_{m\in\N}$ is bounded in $L^{\infty}(\Omega)$. Indeed, by \eqref{hypAMquadr}, for all $x\in \Omega$,
$$
\left\vert H_m(x,u_m(x),\nabla u_m(x))\right\vert \leq M\left(1+\left\vert u_m(x)\right\vert^r+\left\vert \nabla u_m(x)\right\vert^q\right),
$$
and $u_m$ and $\nabla u_m$ are uniformly bounded in $L^{\infty}(\Omega)$ by \eqref{uw2p}. As a consequence, there are a symmetric matrix field $A_{\infty}\in W^{1,\infty}(\Omega,{\mathcal{S}}_n(\R))$, a function $u_{\infty}\in W(\Omega)\cap H^1_0(\Omega)$, a point~$x_{\infty}\in\overline{\Omega}$ and two functions $H_{\infty},\ H^0\in L^{\infty}(\Omega)$ satisfying~\eqref{Ainfty},~\eqref{uinfini} and~\eqref{hopf} as in the proof of Lemma \ref{lembounds}. For every $m\in\N$ and $x\in\Omega$, one has, by \eqref{hypAMquadr},
$$
\begin{array}{lll}
\displaystyle H_m(x,u_m(x),\nabla u_m(x))& \le & \displaystyle  M\big((u_m(x))^r+|\nabla u_m(x)|^q\big)+H_m(x,0,0)\\
& \leq & \displaystyle C_4\big(u_m(x)+|\nabla u_m(x)|\big)+H_m(x,0,0),
\end{array}
$$
from~(\ref{hypAMquadr}), \eqref{uw2p} and the nonnegativity of $u_m$. Here, $C_4>0$ only depends on $\Omega$, $n$, $q$, $M$ and $r$. Finally, one concludes as for the proof of Lemma \ref{lembounds}.\hfill$\Box$\break

\noindent{\bf{Proof of Lemma~\ref{lem41}.}} We first remember that the existence and uniqueness of the weak solutions $V_k\in H^1_0(\Omega^*)\cap L^{\infty}(\Omega^*)$ of~(\ref{defvkter}) is guaranteed by~\cite{bm95,bmp84}. By uniqueness of $V_k$ and radial symmetry of all coefficients, the functions $V_k$ are all radially symmetric. Furthermore, from the boundedness assumptions made in Lemma~\ref{lem41} and from~\cite{bmp88cras,bmp92}, the sequence~$(V_k)_{k\in\N}$ is bounded in~$H^1_0(\Omega^*)$ and in $L^{\infty}(\Omega^*)$, and it is relatively compact in $H^1_0(\Omega^*)$. Therefore, there exists a radially symmetric function~$V_{\infty}\in H^1_0(\Omega^*)\cap L^{\infty}(\Omega^*)$ such that the limit~(\ref{weakVk}) holds, at least for a subsequence.\par
The goal is to show that the function~$V_{\infty}$ solves the limiting equation~(\ref{Vinfty}) in the weak~$H^1_0(\Omega^*)\cap L^{\infty}(\Omega^*)$ sense. Notice that, once this is done, then by uniqueness of the~$H^1_0(\Omega^*)\cap L^{\infty}(\Omega^*)$ solution of this limiting equation~(\ref{Vinfty}), it follows immediately that the whole sequence~$(V_k)_{k\in\N}$ converges to $V_{\infty}$ in the sense of~(\ref{weakVk}).\par
In order to show~(\ref{Vinfty}), as for the proof of Lemma~\ref{vinfty}, the strategy is to work with the one-dimensional equations satisfied by the functions $V_k$, to derive additional bounds and to pass to the limit in a certain sense. Comparing to Lemma~\ref{vinfty}, the additional difficulty will be to pass to the limit in the terms $\alpha_k\,|\nabla V_k|^q$, for which only $L^1$ bounds are available (in particular, estimates similar to~(\ref{kepsilon}) with $-\alpha_k|\nabla V_k|^q$ instead of~$\hat{a}_ke_r\cdot\nabla v_k$ may not be true). For the proof of~(\ref{Vinfty}), we consider separately the cases where the dimension~$n$ is such that~$n\ge 2$, and the case $n=1$.\par
{\it First case: $n\ge 2$.} Call
\be\label{deftilde4}\left\{\baa{ll}
\displaystyle\tilde{\lambda}_k(r)=\frac{1}{n\alpha_nr^{n-1}}\int_{\partial B_r}\lambda_k\,d\theta_r, & \displaystyle\tilde{\lambda}_{\infty}(r)=\frac{1}{n\alpha_nr^{n-1}}\int_{\partial B_r}\lambda_{\infty}\,d\theta_r,\vspace{3pt}\\
\displaystyle\tilde{\alpha}_k(r)=\frac{1}{n\alpha_nr^{n-1}}\int_{\partial B_r}\alpha_k\,d\theta_r, & \displaystyle\tilde{\alpha}_{\infty}(r)=\frac{1}{n\alpha_nr^{n-1}}\int_{\partial B_r}\alpha_{\infty}\,d\theta_r,\vspace{3pt}\\
\displaystyle\tilde{\gamma}_k(r)=\frac{1}{n\alpha_nr^{n-1}}\int_{\partial B_r}\gamma_k\,d\theta_r, & \displaystyle\tilde{\gamma}_{\infty}(r)=\frac{1}{n\alpha_nr^{n-1}}\int_{\partial B_r}\gamma_{\infty}\,d\theta_r,\vspace{3pt}\\
\displaystyle\tilde{V}_k(r)=\frac{1}{n\alpha_nr^{n-1}}\int_{\partial B_r}V_k\,d\theta_r, & \displaystyle\tilde{V}_{\infty}(r)=\frac{1}{n\alpha_nr^{n-1}}\int_{\partial B_r}V_{\infty}\,d\theta_r,\vspace{3pt}\\
\displaystyle\tilde{W}_k(r)=\frac{1}{n\alpha_nr^{n-1}}\int_{\partial B_r}\lambda_k\,\nabla V_k\cdot e_r\,d\theta_r. & \eaa\right.
\ee
From Fubini's theorem, these quantities can be defined for almost every $r\in(0,R)$. The functions $\tilde{\lambda}_k$ and $\tilde{\lambda}_{\infty}$ are in $L^{\infty}_+(0,R)$ and the functions $\tilde{\alpha}_k$, $\tilde{\alpha}_{\infty}$, $\tilde{\gamma}_k$ and $\tilde{\gamma}_{\infty}$ are in~$L^{\infty}(0,R)$. Furthermore, it follows from~(\ref{chapeauxter}) that
\be\label{chapeaux4}
\tilde{\lambda}_k^{-1}\rightharpoonup\tilde{\lambda}_{\infty}^{-1},\ \ \tilde{\lambda}_k^{-q}\tilde{\alpha}_k\rightharpoonup\tilde{\lambda}_{\infty}^{-q}\tilde{\alpha}_{\infty}\ \hbox{ and }\ \tilde{\gamma}_k\rightharpoonup\tilde{\gamma}_{\infty}\ \hbox{ in }L^{\infty}(r_0,R)\hbox{ weak-* }\hbox{ as }k\to+\infty
\ee
for every $r_0\in(0,R)$. On the other hand, the functions $\tilde{V}_k$ and $\tilde{V}_{\infty}$ belong to the space~$L^{\infty}(0,R)\cap H^1_{loc}((0,R])$, the sequence $(\tilde{V}_k)_{k\in\N}$ is bounded in $L^{\infty}(0,R)$ and in $H^1(r_0,R)$ for every $r_0\in(0,R)$, and
\be\label{Vkweak}
\tilde{V}_k\to\tilde{V}_{\infty}\hbox{ in }H^1(r_0,R)\hbox{ strong as }k\to+\infty
\ee
for every $r_0\in(0,R)$, from~(\ref{weakVk}). Lastly, the sequence $(\tilde{W}_k)_{k\in\N}$ is bounded in $L^2(r_0,R)$ for every $r_0\in(0,R)$ and it is straightforward to check that
\be\label{tildeWk}
\tilde{W}_k=\tilde{\lambda}_k\,\tilde{V}_k'\ \hbox{ a.e. in }(0,R)
\ee
for every $k\in\N$.\par
By testing~(\ref{defvkter}) against radially symmetric functions $\varphi\in C^1_c(\Omega^*)$ which vanish in a neighborhood of $0$, it follows that, for every $k\in\N$, the function $\tilde{W}_k$ is in $W^{1,1}(r_0,R)$ for every $r_0\in(0,R)$, and that
$$-\tilde{W}_k'(r)=\frac{n-1}{r}\,\tilde{\lambda}_k(r)\,\tilde{V}_k'(r)+\tilde{\alpha}_k(r)\,|\tilde{V}_k'(r)|^q-\delta_k\,\tilde{V}_k(r)+\tilde{\gamma}_k(r)\ \hbox{ a.e. in }(0,R).$$
Therefore, the sequence $(\tilde{W}_k')_{k\in\N}$ is bounded in $L^1(r_0,R)$ for every $r_0\in(0,R)$, whence $(\tilde{W}_k)_{k\in\N}$ is bounded in $W^{1,1}(r_0,R)$ for every $r_0\in(0,R)$. As a consequence, there is a function $\tilde{W}_{\infty}\in L^{\infty}_{loc}((0,R])$ such that, up to extraction of a subsequence,
\be\label{tildeWinfty}\left\{\baa{l}
\displaystyle\tilde{W}_k\displaystyle{\mathop{\rightharpoonup}_{k\to+\infty}}\tilde{W}_{\infty}\hbox{ in }L^{\infty}(r_0,R)\hbox{ weak-*},\vspace{3pt}\\
\tilde{W}_k\displaystyle{\mathop{\to}_{k\to+\infty}}\tilde{W}_{\infty}\hbox{ in }L^p(r_0,R)\hbox{ strong for every }1\le p<+\infty,\eaa\right.\hbox{ for every }r_0\in(0,R).
\ee
Together with~(\ref{chapeaux4}) and~(\ref{tildeWk}), one gets that $\tilde{V}_k'=\tilde{\lambda}_k^{-1}\tilde{W}_k\rightharpoonup\tilde{\lambda}_{\infty}^{-1}\tilde{W}_{\infty}$ as $k\to+\infty$ in, say,~$L^2(r_0,R)$ weak, for every $r_0\in(0,R)$. Remembering~(\ref{Vkweak}), one concludes that
\be\label{tildeWinftybis}
\tilde{W}_{\infty}=\tilde{\lambda}_{\infty}\,\tilde{V}_{\infty}'\ \hbox{ a.e. in }(0,R).
\ee\par
We shall now show that $V_{\infty}$ is a weak $H^1_0(\Omega^*)\cap L^{\infty}(\Omega^*)$ solution of~(\ref{Vinfty}), that is $I^{\varphi}=0$ for every $\varphi\in H^1_0(\Omega^*)\cap L^{\infty}(\Omega^*)$, where $I^{\varphi}$ is defined in~(\ref{defIvarphi}). Let us fix such a function $\varphi$. Since $n\ge 2$, there is a sequence~$(\phi_m)_{m\in\N}$ of~$C^1_c(\Omega^*)$ functions such that, for every $m\in\N$,
$$0\le\phi_m\le 1\hbox{ in }\Omega^*,\ \phi_m=1\hbox{ in }B_{2^{-m-2}R},\ \phi_m=0\hbox{ in }\Omega^*\backslash B_{2^{-m-1}R}$$
and the sequence $(\phi_m)_{m\in\N}$ is bounded in $H^1_0(\Omega^*)$. Call $\varphi_m=\varphi\,(1-\phi_m)$. The functions $\varphi_m$ are all in $H^1_0(\Omega^*)\cap L^{\infty}(\Omega^*)$ and they are such that
$$\varphi_m=0\hbox{ a.e. in }B_{2^{-m-2}R}\hbox{ and }\varphi_m-\varphi=-\varphi\,\phi_m=0\hbox{ a.e. in }\Omega^*\backslash B_{2^{-m-1}R}.$$
One has
\be\label{varphim}\baa{rcl}
I^{\varphi_m}-I^{\varphi} & = & \displaystyle-\int_{B_{2^{-m-1}R}}\lambda_{\infty}\,\nabla V_{\infty}\cdot\nabla(\varphi\,\phi_m)+\int_{B_{2^{-m-1}R}}\alpha_{\infty}\,|\nabla V_{\infty}|^q\,\varphi\,\phi_m\vspace{3pt}\\
& & \displaystyle-\int_{B_{2^{-m-1}R}}\delta_{\infty}\,V_{\infty}\,\varphi\,\phi_m+\int_{B_{2^{-m-1}R}}\gamma_{\infty}\,\varphi\,\phi_m\eaa
\ee
for every $m\in\N$. The last three integrals converge to $0$ as $m\to+\infty$ from Lebesgue's dominated convergence theorem, since $\alpha_{\infty}\in L^{\infty}(\Omega^*)$, $V_{\infty}\in H^1_0(\Omega^*)\cap L^{\infty}(\Omega^*)$ and $\varphi\in L^{\infty}(\Omega^*)\,(\cap H^1_0(\Omega^*))$ and since the sequence $(\phi_m)_{m\in\N}$ is bounded in $L^{\infty}(\Omega^*)$. The first integral of~(\ref{varphim}) also converges to $0$ as $m\to+\infty$ by Cauchy-Schwarz inequality and Lebesgue's dominated convergence theorem, since $\lambda_{\infty}|\nabla V_{\infty}|$ is in $L^2(\Omega^*)$, $\varphi$ is in $H^1_0(\Omega^*)\cap L^{\infty}(\Omega^*)$ and the sequence~$(\phi_m)_{m\in\N}$ is bounded in $H^1_0(\Omega^*)$. Finally, $I^{\varphi_m}-I^{\varphi}\to0$ as $m\to+\infty$. So, in order to show that $I^{\varphi}=0$ for all $\varphi\in H^1_0(\Omega^*)\cap L^{\infty}(\Omega^*)$, it is sufficient to show it for all~$\varphi\in H^1_0(\Omega^*)\cap L^{\infty}(\Omega^*)$ such that $\varphi=0$ almost everywhere in a ball with positive radius centered at the origin.\par
Let then $\varphi$ be a fixed function in $H^1_0(\Omega^*)\cap L^{\infty}(\Omega^*)$ such that $\varphi=0$ almost everywhere in~$B_{r_0}$ for some $r_0\in(0,R)$, and let us show that $I^{\varphi}=0$. It follows from Fubini's theorem and the definitions~(\ref{deftilde4}) that
\be\label{int0}
I^{\varphi}=\int_{r_0}^R\Big(\tilde{\lambda}_{\infty}(r)\,\tilde{V}_{\infty}'(r)\,\Phi(r)-\tilde{\alpha}_{\infty}(r)\,|\tilde{V}_{\infty}'(r)|^q\,\phi(r)+\delta_{\infty}\,\tilde{V}_{\infty}(r)\,\phi(r)-\tilde{\gamma}_{\infty}(r)\,\phi(r)\Big)\,dr,
\ee
where the functions $\phi$ and $\Phi$ are defined as in~(\ref{defphiPhi}) and are in $L^2(r_0,R)$. Furthermore,~$\phi\in L^{\infty}(r_0,R)$. Since $\tilde{\lambda}_k\,\tilde{V}_k'=\tilde{W}_k\rightharpoonup\tilde{W}_{\infty}=\tilde{\lambda}_{\infty}\,\tilde{V}_{\infty}'$ as $k\to+\infty$ in $L^{\infty}(r_0,R)$ weak-* from~(\ref{tildeWk}),~(\ref{tildeWinfty}) and~(\ref{tildeWinftybis}), one infers that
\be\label{int1}
\int_{r_0}^R\tilde{\lambda}_k(r)\,\tilde{V}_k'(r)\,\Phi(r)\,dr\to\int_{r_0}^R\tilde{\lambda}_{\infty}(r)\,\tilde{V}_{\infty}'(r)\,\Phi(r)\,dr\ \hbox{ as }k\to+\infty.
\ee
Furthermore,
\be\label{alphaVk}
\tilde{\alpha}_k\,|\tilde{V}_k'|^q=\big(\tilde{\lambda}_k^{-q}\,\tilde{\alpha}_k\big)\times\big(\tilde{\lambda}_k^q\,|\tilde{V}_k'|^q\big)
\ee
and
\be\label{alphaVkbis}
\tilde{\lambda}_k^{-q}\,\tilde{\alpha}_k\mathop{\rightharpoonup}_{k\to+\infty}\tilde{\lambda}_{\infty}^{-q}\,\tilde{\alpha}_{\infty}\hbox{ in }L^{\infty}(r_0,R)\hbox{ weak-*}
\ee
from~(\ref{chapeaux4}). On the other hand,
\be\label{alphaVk3}
\tilde{\lambda}_k^q\,|\tilde{V}_k'|^q=|\tilde{W}_k|^q\mathop{\to}_{k\to+\infty}|\tilde{W}_{\infty}|^q=\tilde{\lambda}_{\infty}^q\,|\tilde{V}_{\infty}'|^q\hbox{ in, say, }L^1(r_0,R)\hbox{ strong},
\ee
from~(\ref{tildeWinfty}) and~(\ref{tildeWinftybis}): indeed, since the sequence $(W_k)_{k\in\N}$ is bounded in $L^{\infty}(r_0,R)$ from~(\ref{tildeWinfty}), it follows that, for every $\epsilon>0$, there is a constant $C_{\epsilon}>0$ such that
$$\big||\tilde{W}_k(r)|^q-|\tilde{W}_{\infty}(r)|^q\big|\le C_{\epsilon}\big|\tilde{W}_k(r)-\tilde{W}_{\infty}(r)\big|+\epsilon\hbox{ for every }k\in\N\hbox{ and a.e. }r\in(r_0,R),$$
and the $L^1(r_0,R)$ convergence of $|\tilde{W}_k|^q$ to $|\tilde{W}_{\infty}|^q$ follows then from the $L^1(r_0,R)$ convergence of $\tilde{W}_k$ to $\tilde{W}_{\infty}$ by~(\ref{tildeWinfty}). Putting together~(\ref{alphaVk}),~(\ref{alphaVkbis}) and~(\ref{alphaVk3}) leads to
$$\tilde{\alpha}_k\,|\tilde{V}_k'|^q\mathop{\rightharpoonup}_{k\to+\infty}\big(\tilde{\lambda}_{\infty}^{-q}\,\tilde{\alpha}_{\infty}\big)\times\big(\tilde{\lambda}_{\infty}^q\,|\tilde{V}_{\infty}'|^q\big)=\tilde{\alpha}_{\infty}\,|\tilde{V}_{\infty}'|^q\hbox{ in }L^1(r_0,R)\hbox{ weak},$$
whence
\be\label{int2}
\int_{r_0}^R\tilde{\alpha}_k(r)\,|\tilde{V}_k'(r)|^q\,\phi(r)\,dr\to\int_{r_0}^R\tilde{\alpha}_{\infty}(r)\,|\tilde{V}_{\infty}'(r)|^q\,\phi(r)\,dr\ \hbox{ as }k\to+\infty
\ee
since $\phi\in L^{\infty}(r_0,R)$. Similarly, it follows from~(\ref{chapeaux4}),~(\ref{Vkweak}) and the convergence of $\delta_k$ to $\delta_{\infty}$, that
\be\label{int3}
\int_{r_0}^R\Big(\delta_k\,\tilde{V}_k(r)\,\phi(r)-\tilde{\gamma}_k(r)\,\phi(r)\Big)\,dr\mathop{\to}_{k\to+\infty}\int_{r_0}^R\Big(\delta_{\infty}\,\tilde{V}_{\infty}(r)\,\phi(r)-\tilde{\gamma}_{\infty}(r)\,\phi(r)\Big)\,dr.
\ee
Finally, it follows from~(\ref{int0}),~(\ref{int1}),~(\ref{int2}) and~(\ref{int3}) that
$$\underbrace{\int_{r_0}^R\Big(\tilde{\lambda}_k(r)\,\tilde{V}_k'(r)\,\Phi(r)-\tilde{\alpha}_k(r)\,|\tilde{V}_k'(r)|^q\,\phi(r)+\delta_k\,\tilde{V}_k(r)\,\phi(r)-\tilde{\gamma}_k(r)\,\phi(r)\Big)\,dr}_{=:I^{\varphi}_k}\to I^{\varphi}$$
as $k\to+\infty$. But it follows from Fubini's theorem, the definitions~(\ref{deftilde4}) and the fact that $\varphi=0$ in $B_{r_0}$, that
$$I^{\varphi}_k=\displaystyle\int_{\Omega^*}\lambda_k\nabla V_k\cdot\nabla\varphi-\int_{\Omega^*}\alpha_k\,|\nabla V_k|^q\,\varphi+\int_{\Omega^*}\delta_k\,V_k\,\varphi-\int_{\Omega^*}\gamma_k\,\varphi$$
for every $k\in\N$. Hence, $I^{\varphi}_k=0$ for every $k\in\N$, owing to the definition of $V_k$ in~(\ref{defvkter}). As a conclusion, $I^{\varphi}=0$ for every $\varphi\in H^1_0(\Omega^*)\cap L^{\infty}(\Omega^*)$ such that $\varphi=0$ almost everywhere in a neighborhood of $0$. As already emphasized, this implies that $I^{\varphi}=0$ for every $\varphi\in H^1_0(\Omega^*)\cap L^{\infty}(\Omega^*)$ in this case $n\ge 2$.\par
{\it Second case: $n=1$.} In this case, we work directly with the functions $\lambda_k$, $\lambda_{\infty}$, $\alpha_k$, $\alpha_{\infty}$, $\gamma_k$, $\gamma_{\infty}$, $V_k$ and $V_{\infty}$ defined in $\Omega^*=(-R,R)$ and satisfying~(\ref{chapeauxter}) and~(\ref{weakVk}). Call
$$W_k=\lambda_k\,V_k'.$$
The sequence $(W_k)_{k\in\N}$ is bounded in $L^2(-R,R)$. Furthermore, it follows from~(\ref{defvkter}) that each function $W_k$ is in $W^{1,1}(-R,R)$ and that
$$-W_k'(r)=\alpha_k(r)\,|V_k'(r)|^q-\delta_k\,V_k(r)+\gamma_k(r)\ \hbox{ a.e. in }(-R,R).$$
Hence, the sequence $(W_k)_{k\in\N}$ is bounded in $W^{1,1}(-R,R)$ and there exists a function $W_{\infty}\in L^{\infty}(-R,R)$ such that
$$\left\{\baa{l}
W_k\rightharpoonup W_{\infty}\hbox{ as }k\to+\infty\hbox{ in }L^{\infty}(-R,R)\hbox{ weak-*},\vspace{3pt}\\
W_k\to W_{\infty}\hbox{ as }k\to+\infty\hbox{ in }L^p(-R,R)\hbox{ strong for every }1\le p<+\infty.\eaa\right.$$
One then infers as in the case $n\ge2$ that
$$W_{\infty}=\lambda_{\infty}\,V_{\infty}'\ \hbox{ a.e. in }(-R,R).$$
The same arguments as the ones used in the last part of the proof in the case $n\ge 2$ then imply that, for every $\varphi\in H^1_0(-R,R)\,(\subset L^{\infty}(-R,R))$,
$$\baa{l}
0=\displaystyle\int_{-R}^R\lambda_k(r)\,V'_k(r)\,\varphi'(r)\,dr-\int_{-R}^R\alpha_k(r)\,|V'_k(r)|^q\,\varphi(r)\,dr+\int_{-R}^R\delta_k\,V_k(r)\,\varphi(r)\,dr\vspace{3pt}\\
\qquad\qquad\qquad\qquad\displaystyle-\int_{-R}^R\gamma_k(r)\,\varphi(r)\,dr\vspace{3pt}\\
\qquad\displaystyle{\mathop{\to}_{k\to+\infty}}\displaystyle\int_{-R}^R\lambda_{\infty}(r)\,V'_{\infty}(r)\,\varphi'(r)\,dr-\int_{-R}^R\alpha_{\infty}(r)\,|V'_{\infty}(r)|^q\,\varphi(r)\,dr+\int_{-R}^R\delta_{\infty}\,V_{\infty}(r)\,\varphi(r)\,dr\vspace{3pt}\\
\qquad\qquad\qquad\qquad\displaystyle-\int_{-R}^R\gamma_{\infty}(r)\,\varphi(r)\,dr.\eaa$$
As a consequence, the limiting integral is equal to $0$ for every $\varphi\in H^1_0(-R,R)\,(\subset L^{\infty}(-R,R))$, which means that $V_{\infty}$ is the weak solution of~(\ref{Vinfty}).\par
The proof of Lemma~\ref{lem41} is thereby complete in both cases $n\ge 2$ and $n=1$.\hfill\fin


\SE{Appendix}\label{secapp}

We first give in Section~\ref{sec51} some sufficient conditions on $H$ for the existence and uniqueness of a solution $u$ to problem~(\ref{equ}) in a domain $\Omega$ when $H$ grows at most linearly in $|p|$ . We refer to the comments after the statements of Theorem~\ref{th1} and Corollary~\ref{cor1} for the presentation of these conditions. In Section~\ref{sec52}, when~$H$ grows at most quadratically in~$|p|$, we give some sufficient conditions on~$H$ for any weak solution of~(\ref{equ}) to be in~$W(\Omega)$.


\subsection{The case of at most linear growth with respect to the gradient}\label{sec51}

\begin{theo}\label{thapp}
Consider problem~$(\ref{equ})$ in a bounded domain~$\Omega\subset\R^n$ of class $C^2$ and assume that~$(\ref{ALambda})$ holds with $\Lambda\in L^{\infty}_+(\Omega)$.
\begin{itemize}
\item[{\rm{(i)}}] If $H$ satisfies~$(\ref{hypH2})$, then problem~$(\ref{equ})$ has at most one solution in $H^1_0(\Omega)$.
\item[{\rm{(ii)}}] If $H$ satisfies~$(\ref{hypH2})$ and~$(\ref{hypH3})$, then problem~$(\ref{equ})$ has exactly one solution $u$ in $H^1_0(\Omega)$.\par
\item[{\rm{(iii)}}] If $H$ satisfies~$(\ref{hypH2})$, $(\ref{hypH3})$ and~$(\ref{hypH4})$, then $u\in W(\Omega)$.\par
\item[{\rm{(iv)}}] If $H$ satisfies~$(\ref{hypH2})$,~$(\ref{hypH3})$,~$(\ref{hypH4})$ and~$(\ref{hypH5})$, then $u>0$ in $\Omega$ and $\partial_{\nu}u<0$ on $\partial\Omega$, where $\nu$ denotes the outward unit normal to $\partial\Omega$.
\end{itemize}
\end{theo}

The proof of the ``uniqueness'' part relies on the following form of the weak maximum principle:

\begin{pro} \label{weakmax}
Let~$\Omega\subset\R^n$ be a bounded domain of class $C^2$ and assume that~$(\ref{ALambda})$ holds with $\Lambda\in L^{\infty}_+(\Omega)$. Let $u\in H^1_0(\Omega)$ and assume that there exists a function $g\in L^2(\Omega)$ such that~$u(x)$ and~$g(x)$ have the same sign for almost every $x\in \Omega$, and a function $B\in L^n(\Omega)$ such that
\begin{equation} \label{inequ}
-{\rm{div}}(A\nabla u)+g\leq B\left\vert \nabla u\right\vert
\end{equation}
in the weak $H^1_0(\Omega)$ sense. Then $u\leq 0$ in $\Omega$.
\end{pro}

\noindent{\bf{Proof.}} The proof is inspired by the one of Proposition 2.1 in \cite{Porretta}. Let $k>0$. Using $v_k:=(u-k)^+\in H^1_0(\Omega)$ as a test function in \eqref{inequ}, one obtains
$$
\int_{\Omega} A(x)\nabla u(x)\cdot \nabla v_k(x)\,dx+ \int_{\Omega} g(x)v_k(x)\,dx \leq \int_{\Omega}  B(x)\left\vert \nabla u(x)\right\vert v_k(x)\,dx.
$$
Notice that
$$
\int_{\Omega} g(x)v_k(x)dx\geq 0.
$$
Indeed, if $u(x)\leq k$, $v_k(x)=0$, while, if $u(x)>k$, then $u(x)>0$, so that $g(x)\geq 0$ by the assumption on $g$. Therefore,
$$
\int_{\Omega} A(x)\nabla u(x)\cdot \nabla v_k(x)\,dx\leq \int_{E_k}  B(x)\left\vert \nabla u(x)\right\vert v_k(x)\,dx,
$$
where
$$
E_k:=\big\{x\in \Omega;\ u(x)>k\mbox{ and } \left\vert \nabla u(x)\right\vert>0\big\}.
$$
One concludes as in the proof of Proposition 2.1 in \cite{Porretta}, using in particular the fact that~$B\in L^n(\Omega)$.\hfill\fin\break

The proof of the ``existence'' part in Theorem \ref{thapp} relies on a fixed point argument, summarized in the following proposition:

\begin{pro} \label{fixed}
Let~$\Omega\subset\R^n$ be a bounded domain of class $C^2$ and assume that~$(\ref{ALambda})$ holds with $\Lambda\in L^{\infty}_+(\Omega)$ and that $H$ satisfies~\eqref{hypH2} and~\eqref{hypH3}. Then:\par
\begin{itemize}
\item[{\rm{(i)}}] For all $v\in H^1_0(\Omega)$, there exists a unique function $u$ in $H^1_0(\Omega)$ solving
\begin{equation} \label{uv}
-{\rm{div}}(A\nabla u)+H(x,v,\nabla v)=0.
\end{equation}
This allows to define $u:=Tv$, and $T$ is continuous from $H^1_0(\Omega)$ into itself.
\item[{\rm{(ii)}}] The map $T$ is compact from $H^1_0(\Omega)$ into itself $($in the sense that $T$ maps bounded sets into precompact sets$)$.
\item[{\rm{(iii)}}] There exists $M\geq 0$ such that, for all $u\in H^1_0(\Omega)$ and all $\sigma\in [0,1]$ such that $u=\sigma Tu$, one has
$$\left\Vert u\right\Vert_{H^1_0(\Omega)}\leq M.$$
\end{itemize}
\end{pro}

\noindent{\bf{Proof.}} The proof is inspired by the one of Proposition 2.2 in \cite{Porretta}.\par
{\it Proof of}~(i). Observe first that there exists a constant $C>0$ such that, for all $v\in H^1_0(\Omega)$, $H(\cdot,v,\nabla v)\in H^{-1}(\Omega)$ and
\begin{equation} \label{h-1}
\left\Vert H(\cdot,v,\nabla v)\right\Vert_{H^{-1}(\Omega)}\leq C\big(\left\Vert v\right\Vert_{H^1_0(\Omega)}+1\big).
\end{equation}
Indeed, by \eqref{hypH3}, $\left\vert H(x,v(x),\nabla v(x))\right\vert\leq \beta(x)\left(1+\left\vert v(x)\right\vert+\left\vert \nabla v(x)\right\vert\right)$ in $\Omega$. If $n\geq 2$, since $1+\left\vert v\right\vert+\left\vert \nabla v\right\vert$ belongs to $L^2(\Omega)$ and $\beta\in L^n(\Omega)$, $H(\cdot,v,\nabla v)\in L^{\frac{2n}{n+2}}(\Omega)\subset H^{-1}(\Omega)$. If $n=1$, since $1+\left\vert v\right\vert+\left\vert v^{\prime}\right\vert$ still belongs to $L^2(\Omega)$ and $\beta\in L^2(\Omega)$, $H(\cdot,v,v^{\prime})\in L^1(\Omega)\subset H^{-1}(\Omega)$. The Lax-Milgram theorem then shows that there exists a unique $u\in H^1_0(\Omega)$ solving \eqref{uv}. Moreover, there exists $C>0$ such that, for all $v\in H^1_0(\Omega)$,
\begin{equation} \label{Tbounded}
\left\Vert Tv\right\Vert_{H^1_0(\Omega)}\leq C\big(\left\Vert v\right\Vert_{H^1_0(\Omega)}+1\big).
\end{equation}\par
Let us now prove that $T$ is continuous. Let $(v_l)_{l\in\N}$ and $v$ in $H^1_0(\Omega)$ be such that $v_l\rightarrow v$ in $H^1_0(\Omega)$ as $l\to+\infty$. We first claim that $H(\cdot,v_l,\nabla v_l)\rightarrow H(\cdot,v,\nabla v)$ in $H^{-1}(\Omega)$ as $l\to+\infty$. Indeed, by~\eqref{hypH2},
$$
\begin{array}{lll}
\displaystyle \left\vert H(\cdot,v_l,\nabla v_l)-H(\cdot,v,\nabla v)\right\vert & \leq & \displaystyle \left\vert H(\cdot,v_l,\nabla v_l)-H(\cdot,v,\nabla v_l)\right\vert + \left\vert H(\cdot,v,\nabla v_l)-H(\cdot,v,\nabla v)\right\vert\vspace{3pt}\\
& \leq  & \displaystyle\omega\left\vert v_l-v\right\vert + \alpha\,\big(1+\left\vert v\right\vert^{2/n}\big)\left\vert \nabla (v_l-v)\right\vert.
\end{array}
$$
On the one hand, $v_l-v\to0$ in $L^2(\Omega)\subset H^{-1}(\Omega)$ as $l\to+\infty$. On the other hand, when~$n\geq 3$, since $\alpha\in L^{n^2/2}(\Omega)$, $1+\left\vert v\right\vert^{2/n}\in L^{n^2/(n-2)}(\Omega)$ by Sobolev embeddings and $\nabla (v_l-v)\to0$ in $L^2(\Omega)$ as $l\to+\infty$, the H\"older inequality entails that $\alpha\,(1+\left\vert v\right\vert^{2/n})\left\vert \nabla (v_l-v)\right\vert\to0$ in $L^{2n/(n+2)}(\Omega)\subset H^{-1}(\Omega)$ as $l\to+\infty$. If $n=2$, since $\alpha\in L^r(\Omega)$ for some $r>2$, $1+\left\vert v\right\vert\in L^s(\Omega)$ for all $s\in (2,+\infty)$ and $\nabla (v_l-v)\rightarrow 0$ in $L^2(\Omega)$ as $l\to+\infty$, it follows that~$\alpha\,\left(1+\left\vert v\right\vert\right)\left\vert \nabla (v_l-v)\right\vert\to0$ as $l\to+\infty$ in $L^{\chi}(\Omega)$ for all $1<\chi<2r/(r+2)$ and~$L^{\chi}(\Omega)\subset H^{-1}(\Omega)$. Finally, when $n=1$, since $\alpha\in L^2(\Omega)$, $1+\left\vert v\right\vert^2\in L^{\infty}(\Omega)$ and~$(v_l-v)^{\prime}\rightarrow 0$ in $L^2(\Omega)$ as $l\to+\infty$, one gets that $\alpha\left(1+\left\vert v\right\vert\right)\left\vert (v_l-v)^{\prime}\right\vert\to0$ in $L^1(\Omega)\subset H^{-1}(\Omega)$ as $l\to+\infty$. To sum up, $H(\cdot,v_l,\nabla v_l)\rightarrow H(\cdot,v,\nabla v)$ in $H^{-1}(\Omega)$ as $l\to+\infty$, in any dimension $n\ge 1$. Therefore, if $u_l=Tv_l$ and $u=Tv$, one has
$$
-\div(A\nabla (u-u_l))=H(\cdot,v_l,\nabla v_l)-H(\cdot,v,\nabla v)
$$
and since the right-hand side goes to $0$ in $H^{-1}(\Omega)$ as $l\to+\infty$, one obtains that $Tv_l\rightarrow Tv$ in~$H^1_0(\Omega)$ as $l\to+\infty$, which shows that $T$ is continuous.\hfill\break\par
{\it Proof of}~(ii). We want to show that $T$ is compact. Let $(v_l)_{l\in\N}$ be a bounded sequence in~$H^1_0(\Omega)$ and denote, for $l\in\N$,
$$u_l=Tv_l.$$
Since \eqref{Tbounded} holds, the sequence $(u_l)_{l\in\N}$ is bounded in $H^1_0(\Omega)$, so that, up to a subsequence, there exists $u\in H^1_0(\Omega)$ such that $u_l\rightarrow u$ weakly in $H^1(\Omega)$, strongly in $L^2(\Omega)$ and almost everywhere in $\Omega$, as $l\to+\infty$. For all $l\in\N$,
$$
\int_{\Omega} A(x)\nabla u_l(x)\cdot \nabla (u_l-u)(x)\,dx=-\int_{\Omega} H(x,v_l(x),\nabla v_l(x))\,(u_l-u)(x)\,dx,
$$
so that
\be\label{inequlvl}
\begin{array}{lll}
\displaystyle\mathop{{\rm{ess}}\,{\rm{inf}}}_{\Omega}\Lambda\,\times\int_{\Omega}Ê\left\vert \nabla (u_l-u)(x)\right\vert^2 dx  & \leq & \displaystyle \int_{\Omega} A(x)\nabla (u_l-u)(x)\cdot \nabla (u_l-u)(x)\,dx\vspace{3pt}\\
& = & \displaystyle -\int_{\Omega} H(x,v_l(x),\nabla v_l(x))\,(u_l-u)(x)\,dx\vspace{3pt}\\
& & \displaystyle-\int_{\Omega} A(x)\nabla u(x)\cdot \nabla (u_l-u)(x)\,dx\vspace{3pt}\\
& = & \displaystyle -\int_{\Omega} H(x,v_l(x),\nabla v_l(x))\,(u_l-u)(x)\,dx\,+\,\varepsilon_l,
\end{array}
\ee
where $\varepsilon_l\to0$ as $l\to+\infty$ by the weak convergence of $u_l$ to $u$ in $H^1_0(\Omega)$.\par
Let us now estimate the first term of the last right-hand side of the previous formula. For all $l\in\N$ and $m>0$, one has
\be\label{inequlvl2}
\begin{array}{l}
\displaystyle \left\vert \int_{\Omega} H(x,v_l(x),\nabla v_l(x))\,(u_l-u)(x)\,dx \right\vert\vspace{3pt}\\
\qquad\qquad\qquad\qquad\leq\displaystyle\underbrace{\int_{\left\vert u_l-u\right\vert\leq m} \beta(x)\left(1+\left\vert v_l(x)\right\vert+\left\vert \nabla v_l(x)\right\vert\right)\left\vert (u_l-u)(x)\right\vert dx}_{=:A_{m,l}}\vspace{3pt}\\
\qquad\qquad\qquad\qquad\quad+\ \displaystyle\underbrace{\int_{\left\vert u_l-u\right\vert> m} \beta(x)\left(1+\left\vert v_l(x)\right\vert+\left\vert \nabla v_l(x)\right\vert\right)\left\vert (u_l-u)(x)\right\vert dx}_{=:B_{m,l}}.
\end{array}
\ee
Let us first show that $A_{m,l}\to0$ as $l\to+\infty$ for all $m>0$. Indeed, on the one hand, if~$n\geq 3$, since $\beta\in L^n(\Omega)$ and the sequence $(1+\left\vert v_l\right\vert+\left\vert \nabla v_l\right\vert)_{l\in\N}$ is bounded in $L^2(\Omega)$, the H\"older inequality yields
$$
\begin{array}{l}
\displaystyle \int_{\left\vert u_l-u\right\vert\leq m} \beta(x)\left(1+\left\vert v_l(x)\right\vert+\left\vert \nabla v_l(x)\right\vert\right)\left\vert (u_l-u)(x)\right\vert dx\vspace{3pt}\\
\qquad\qquad\qquad\leq\displaystyle \left\Vert \beta\right\Vert_{L^n(\Omega)}\times \left\Vert 1+\left\vert v_l\right\vert + \left\vert \nabla v_l\right\vert\right\Vert_{L^2(\Omega)}\times\left(\int_{\left\vert u_l-u\right\vert\leq m} \left\vert (u_l-u)(x)\right\vert^{\frac{2n}{n-2}} dx\right)^{\frac{n-2}{2n}}\vspace{3pt}\\
\qquad\qquad\qquad\leq\displaystyle C_1\times\left(\int_{\left\vert u_l-u\right\vert\leq m} \left\vert (u_l-u)(x)\right\vert^{\frac{2n}{n-2}} dx\right)^{\frac{n-2}{2n}},
\end{array}
$$
where $C_1>0$ is independent of $l$ and $m>0$, and the dominated convergence theorem therefore entails that $\lim_{l\rightarrow +\infty} A_{m,l}=0$ for all $m>0$. On the other hand, if $n=1$ or $n=2$, the Cauchy-Schwarz inequality yields
$$
\begin{array}{l}
\displaystyle \int_{\left\vert u_l-u\right\vert\leq m} \beta(x)\left(1+\left\vert v_l(x)\right\vert+\left\vert \nabla v_l(x)\right\vert\right)\left\vert (u_l-u)(x)\right\vert dx\vspace{3pt}\\
\qquad\qquad\leq\displaystyle \left(\int_{\Omega} \left(1+\left\vert v_l(x)\right\vert+\left\vert \nabla v_l(x)\right\vert\right)^2dx\right)^{\frac 12}\times\left(\int_{\left\vert u_l-u\right\vert\leq m} \beta^2(x)\left\vert (u_l-u)(x)\right\vert^2dx\right)^{\frac 12},
\end{array}
$$
and since $\beta\in L^2(\Omega)$, the dominated convergence theorem gives again that $\lim_{l\rightarrow +\infty} A_{m,l}=0$ for all $m>0$.\par
Let us now bound $B_{m,l}$ from above. Since the sequence $(u_l-u)_{l\in\N}$ is bounded in $H^1_0(\Omega)$, the Sobolev embeddings imply that~$(u_l-u)_{l\in\N}$ is bounded in $L^{2^{\ast}}(\Omega)$, where
$$
2^{\ast}=\left\{
\begin{array}{ll}
\displaystyle \frac{2n}{n-2} & \mbox{ if }n\geq 3,\vspace{3pt}\\
\displaystyle \mbox{ any }s\in (2,+\infty) & \mbox{ if }n=2,\vspace{3pt}\\
\infty & \mbox{ if }n=1.
\end{array}
\right.
$$
Since the sequence $(1+\left\vert v_l\right\vert+\left\vert \nabla v_l\right\vert)_{l\in\N}$ is bounded in $L^2(\Omega)$, one infers from the H\"older inequality that, for all $l\in\N$ and $m>0$,
\be\label{Bml}
B_{m,l}\leq C_2\left(\int_{\left\vert u_l-u\right\vert> m}\beta(x)^tdx\right)^{1/t},
\ee
where $t$ is given in~\eqref{hypH3} and $C_2>0$ is independent of $l\in\N$ and $m>0$.\par
Finally, it follows from~\eqref{inequlvl},~\eqref{inequlvl2},~\eqref{Bml} and $\lim_{l\to+\infty}\epsilon_l=\lim_{l\to+\infty}A_{m,l}=0$ that, for all~$m>0$,
$$
\displaystyle\mathop{{\rm{ess}}\,{\rm{inf}}}_{\Omega}\Lambda\,\times\limsup_{l\to+\infty}\int_{\Omega}Ê\left\vert \nabla(u_l-u)\right\vert^2\leq C_2\sup_{l\in\N}\left(\int_{\left\vert u_l-u\right\vert> m}\beta(x)^tdx\right)^{1/t}.
$$
Since the measure $E\mapsto \int_E\beta(x)^tdx$ is absolutely continuous with respect to the Lebesgue measure, for all $\varepsilon>0$, there exists $\delta>0$ such that $\int_E\beta(x)^tdx\leq \varepsilon$ whenever $\left\vert E\right\vert<\delta$. As a consequence, since the sequence $(u_l-u)_{l\in\N}$ is bounded in $L^2(\Omega)$ and
$$
\Big\vert\big\{x\in \Omega;\ \left\vert (u_l-u)(x)\right\vert> m\big\}\Big\vert\leq\frac{\Vert u_l-u\Vert_{L^2(\Omega)}^2}{m^2}
$$
for all $l\in\N$ and $m>0$, there holds
$$
\lim_{m\rightarrow +\infty} \sup_{l\in\N}\left(\int_{\left\vert u_l-u\right\vert> m}\beta(x)^tdx\right)^{1/t}=0.
$$
We finally conclude that
$$
\lim_{l\rightarrow +\infty} \int_{\Omega}Ê\left\vert \nabla(u_l-u)\right\vert^2=0,
$$
which proves that $T$ is compact.\hfill\break\par
{\it Proof of~}(iii). Assume by contradiction that (iii) is false. For all $l\in\N$, pick up a function~$u_l\in H^1_0(\Omega)$ and a number $\sigma_l\in [0,1]$ such that
$$u_l=\sigma_lT(u_l)\ \hbox{ and }\ \left\Vert u_l\right\Vert_{H^1_0(\Omega)}>l.$$
Define $v_l=u_l/\left\Vert u_l\right\Vert_{H^1_0(\Omega)}$, so that $\left\Vert v_l\right\Vert_{H^1_0(\Omega)}=1$, for all $l\in\N$. Up to a subsequence, there exists $v\in H^1_0(\Omega)$ such that $v_l\rightarrow v$ weakly in $H^1(\Omega)$, strongly in $L^2(\Omega)$ and almost everywhere in $\Omega$, as $l\to+\infty$. For each $l\in\N$, the function $v_l$ satisfies
\be\label{eqvl}
-\div(A\nabla v_l)=-\sigma_l\frac{H(x,u_l,\nabla u_l)}{\|u_l\|_{H^1_0(\Omega)}}
\ee
in the weak $H^1_0(\Omega)$ sense, and
$$\left|-\sigma_l\frac{H(x,u_l,\nabla u_l)}{\|u_l\|_{H^1_0(\Omega)}}\right|\le\beta(x)\left(\frac{1}{\|u_l\|_{H^1_0(\Omega)}}+|v_l(x)|+|\nabla v_l(x)|\right)$$
by~\eqref{hypH3}. By multiplying~\eqref{eqvl} by $v_l-v$ and arguing as in point~(ii) above, one infers that~$v_l\to v$ strongly in $H^1_0(\Omega)$ as $l\to+\infty$. In particular, $\left\Vert v\right\Vert_{H^1_0(\Omega)}=1$.\par
On the other hand, it follows from~\eqref{hypH2} and~\eqref{hypH3} that, for all $l\in\N$,
$$\begin{array}{lll}
\displaystyle -\div(A\nabla u_l)+\omega^{-1}\sigma_lu_l & \leq & \displaystyle -\div(A\nabla u_l)+\sigma_l\left(H(x,u_l,\nabla u_l)-H(x,0,\nabla u_l)\right)\vspace{3pt}\\
& = & \displaystyle -\sigma_lH(x,0,\nabla u_l)\vspace{3pt}\\
& \leq & \displaystyle \sigma_l \beta(x) \left(1+\left\vert \nabla u_l\right\vert\right)
\end{array}
$$
in the weak $H^1_0(\Omega)$ sense, whence
$$-\div(A\nabla v_l)+\omega^{-1}\sigma_l v_l\leq \sigma_l \beta(x) \left(\frac 1{\left\Vert u_l\right\Vert_{H^1_0(\Omega)}}+\left\vert \nabla v_l\right\vert\right)
$$
in the weak $H^1_0(\Omega)$ sense. Up to a subsequence, one can assume that $\sigma_l\to\sigma\in[0,1]$ as $l\to+\infty$. Therefore, letting $l$ go to $+\infty$ in the previous inequality, one obtains
$$
-\div(A\nabla v)+\omega^{-1}\sigma v\leq \sigma\,\beta(x)\left\vert \nabla v(x)\right\vert
$$
in the weak $H^1_0(\Omega)$ sense. Since $\sigma\,\beta\in L^n(\Omega)$ for any $n\ge 1$, Proposition~\ref{weakmax} above or Proposition 2.1 in \cite{Porretta} show that $v\leq 0$ a.e. in $\Omega$. Similarly, there holds
$$
\begin{array}{lll}
\displaystyle -\div(A\nabla u_l)+\omega\,\sigma_lu_l & \geq & \displaystyle -\div(A\nabla u_l)+\sigma_l\left(H(x,u_l,\nabla u_l)-H(x,0,\nabla u_l)\right)\vspace{3pt}\\
& = & \displaystyle -\sigma_lH(x,0,\nabla u_l)\vspace{3pt}\\
& \geq & \displaystyle -\sigma_l \beta(x) \left(1+\left\vert \nabla u_l\right\vert\right)
\end{array}
$$
and
$$
-\div(A\nabla v_l)+\omega\,\sigma_l v_l\geq -\sigma_l \beta \left(\frac 1{\left\Vert u_l\right\Vert_{H^1_0(\Omega)}}+\left\vert \nabla v_l\right\vert\right),
$$
whence
$$
-\div(A\nabla v)+\omega\sigma v\geq -\sigma \beta \left\vert \nabla v\right\vert
$$
in the weak $H^1_0(\Omega)$ sense. One gets as above that $-v\leq 0$ a.e. in~$\Omega$, so that $v=0$ a.e. in~$\Omega$, which is impossible since $\left\Vert v\right\Vert_{H^1_0(\Omega)}=1$. The proof of Proposition~\ref{fixed} is thereby complete.\hfill\fin\break

The proof of Theorem \ref{thapp} also requires the following regularity result:

\begin{lem} \label{sobolev}
Let~$\Omega\subset\R^n$ be a bounded domain of class $C^2$ and assume that~$(\ref{ALambda})$ holds with~$\Lambda\in L^{\infty}_+(\Omega)$, that $H$ satisfies~\eqref{hypH3} and~\eqref{hypH4}, and that $n\ge 2$. If $u\in H^1_0(\Omega)\cap W^{2,p_0}(\Omega)$ is a weak solution of~\eqref{equ} with $1\le p_0<n$, then $u\in W^{2,p_1}(\Omega)$ with $p_1>n$ given by
$$
\frac 1{p_1}=\frac 1{p_0}-\frac 1n.
$$
\end{lem}

\noindent{\bf{Proof.}} By Sobolev embeddings, $u\in W^{1,p_1}(\Omega)$. Since
$$|H(x,u(x),\nabla u(x))|\le\|\beta\|_{L^{\infty}(\Omega)}\times\big(1+|u(x)|+|\nabla u(x)|\big),$$
the function $x\mapsto H(x,u(x),\nabla u(x))$ is in $L^{p_1}(\Omega)$. Therefore, $u\in W^{2,p_1}(\Omega)$ by elliptic regula\-rity.\hfill\fin\break

\noindent{\bf{Proof of Theorem \ref{thapp}.}} The proof if inspired by the one of Theorem 2.1 in \cite{Porretta}. We first establish~(i). Let $u_1$ and $u_2$ be two solutions of~\eqref{equ} in $H^1_0(\Omega)$ and define $u:=u_1-u_2$. Then
$$
\begin{array}{lll}
\displaystyle -\div (A\nabla u)(x) +H(x,u_1,\nabla u_1)-H(x,u_2,\nabla u_1) & = & \displaystyle H(x,u_2,\nabla u_2)-H(x,u_2,\nabla u_1)\vspace{3pt}\\
& \leq & \displaystyle \alpha(x)\left(1+\left\vert u_2(x)\right\vert^{2/n}\right)\left\vert \nabla u(x)\right\vert
\end{array}
$$
in the weak $H^1_0(\Omega)$ sense. If $g(x):=H(x,u_1(x),\nabla u_1(x))-H(x,u_2(x),\nabla u_1(x))$ for all $x\in \Omega$, then, by~\eqref{hypH2}, the $L^2(\Omega)$ function $g$ has the same sign as~$u$ almost everywhere in $\Omega$. Moreover, if $n\geq 3$, $u_2\in L^{2n/(n-2)}(\Omega)$, which entails that the function
$$
B(x):=\alpha(x)\left(1+\left\vert u_2(x)\right\vert^{2/n}\right)
$$
belongs to $L^n(\Omega)$, since $\alpha\in L^{n^2/2}(\Omega)$. When $n=2$ (resp. $n=1$), argue similarly, using the facts that $\alpha\in L^r(\Omega)$ with $2<r<+\infty$ (resp. $\alpha\in L^2(\Omega)\subset L^1(\Omega)$) and $u_2$ belongs to $L^s(\Omega)$ for all $s<+\infty$ (resp. to $L^{\infty}(\Omega)$), to conclude that $B\in L^n(\Omega)$. In all cases, Proposition~\ref{weakmax} yields that $u\leq 0$ in $\Omega$. Intertwining the roles of $u_1$ and $u_2$ yields $u=0$, so that $u_1=u_2$.\par
Assertion (ii) follows at once from Proposition \ref{fixed} and the Leray-Schauder fixed point theorem.\par
For (iii), since $\beta\in L^{\infty}(\Omega)$, the function $x\mapsto H(x,u,\nabla u)$ belongs to $L^2(\Omega)$, where, by~(i) and~(ii), $u$ denotes the solution of~\eqref{equ} in~$H^1_0(\Omega)$. That means that $-\div(A\nabla u)\in L^2(\Omega)$, and the standard elliptic regularity yields $u\in W^{2,2}(\Omega)$. If $n\ge 2$, using Lemma~\ref{sobolev} repeatedly (a finite number of times only depending on $n$), one obtains that $u\in W^{2,p}(\Omega)$ for some~$p>n$, so that $u\in W^{1,\infty}(\Omega)$. It follows that $x\mapsto H(x,u(x),\nabla u(x))\in L^{\infty}(\Omega)$ and therefore~$u\in W^{2,\tilde{p}}(\Omega)$ for all $1\le\tilde{p}<\infty$ by elliptic regularity. If $n=1$, then $u\in W^{2,2}(\Omega)\subset W^{1,\infty}(\Omega)$ and the same conclusion that $u\in W(\Omega)$ follows.\par
For (iv), one has, by~\eqref{hypH2} and~\eqref{hypH5},
$$
-\mbox{div}(A\nabla u)+H(x,u,\nabla u)-H(x,0,\nabla u)+H(x,0,\nabla u)=0
$$
and
$$-\mbox{div}(A\nabla u)+\omega\,u+\varpi\,|\nabla u|\ge0$$
in the weak $H^1_0(\Omega)$ sense (and a.e. in $\Omega)$, with $\omega>0$ and~$\varpi\ge0$ given in~(\ref{hypH2}) and~(\ref{hypH5}). It follows that~$u\geq 0$ in $\Omega$ by Proposition~\ref{weakmax} above or Proposition~2.1 of~\cite{Porretta}. Furthermore, $u$ does not identically vanish in $\Omega$ since $H(\cdot,0,0)$ is not identically zero. As for the solution $u_{\infty}$ of~\eqref{uinfty}, it follows then from Theorem~9.6 in~\cite{gt} and the discussion there about the strong maximum principle and the Hopf lemma for strong solutions, that $u>0$ in $\Omega$ and $\partial_{\nu}u<0$ on $\partial\Omega$. The proof of Theorem~\ref{thapp} is thereby complete.\hfill\fin


\subsection{The case of at most quadratic growth with respect to the gradient}\label{sec52}

For problem~(\ref{equ}) with an at most quadratic growth in $|\nabla u|$, we give some sufficient conditions on $H$ for any weak solution of~(\ref{equ}) to actually belong to the space~$W(\Omega)$.

\begin{theo}\label{thapp2}
Consider problem~$(\ref{equ})$ in a bounded domain~$\Omega\subset\R^n$ of class $C^2$ and assume that~$(\ref{ALambda})$ holds with $\Lambda\in L^{\infty}_+(\Omega)$. If $H$ satisfies~$(\ref{hypHrq})$ with $1\le q\le 2$, $q<1+2/n$, $M\ge 0$, $r\ge0$ and $r(n-2)<n+2$, then any weak solution $u\in H^1_0(\Omega)$ of~$(\ref{equ})$ belongs to~$W(\Omega)$.
\end{theo}

Similarly to the case of at most linear growth with respect to the gradient, the proof of Theorem \ref{thapp2} goes through the following regularity result:

\begin{lem} \label{sobolevbis}
Under the assumptions of Theorem~$\ref{thapp2}$ with $n\geq 3$, let
$$p_0\in \left(\frac{2n}{n+2},\frac n2\right)$$
and $u\in H^1_0(\Omega)\cap W^{2,p_0}(\Omega)$ be a solution of \eqref{equ}. Write
$$\frac 1{p_0}=\frac{n+2}{2n}-\eta$$
with $\eta>0$. Then there exists $\varepsilon(n,\eta)>0$ only depending on the dimension $n$ and on $\eta$ such that $u\in W^{2,p_1}(\Omega)$ where $p_1>1$ meets 
\begin{equation} \label{condp2}
\frac 1{p_1}<\frac 1{p_0}-\varepsilon(n,\eta).
\end{equation}
\end{lem}

\noindent{\bf{Proof.}} By Sobolev embeddings, $u\in W^{1,p_2}(\Omega)$ with $1/p_2=1/p_0-1/n$ and $u\in L^{p_3}(\Omega)$ with~$1/p_3=1/p_0-2/n$. It follows that
$\left\vert u\right\vert^r\in L^{p_3/r}(\Omega)$, and 
$$\baa{rcl}
\displaystyle\frac r{p_3}<\frac{n+2}{n-2}\left(\frac 1{p_0}-\frac 2n\right) & = & \displaystyle\frac 1{p_0}+\frac 4{(n-2)p_0}-\frac{2n+4}{n(n-2)}\vspace{3pt}\\
& = & \displaystyle\frac 1{p_0}+\frac 4{n-2}\left(\frac{n+2}{2n}-\eta\right)-\frac{2n+4}{n(n-2)}\ =\ \frac 1{p_0}-\frac{4\eta}{n-2}.\eaa
$$
Moreover, $\left\vert \nabla u\right\vert^q\in L^{p_2/q}(\Omega)$, and 
$$
\frac q{p_2}<\frac{n+2}{n}\left(\frac 1{p_0}-\frac 1n\right)=\frac 1{p_0}+\frac{2}{np_0}-\frac{n+2}{n^2}=\frac 1{p_0}+\frac 2n\left(\frac{n+2}{2n}-\eta\right)-\frac{n+2}{n^2}=\frac 1{p_0}-\frac{2\eta}n.
$$
Thus, $x\mapsto H(x,u(x),\nabla u(x))\in L^{p_1}(\Omega)$ where $p_1>1$ meets \eqref{condp2}. Elliptic regularity gives the conclusion.\hfill\fin\break

\noindent{\bf{Proof of Theorem \ref{thapp2}.}} Assume first that $n\geq 3$. Since $u\in H^1_0(\Omega)$, $u\in L^{2n/(n-2)}(\Omega)$. This shows that $\left\vert u\right\vert^r\in L^{2n/(r(n-2))}(\Omega)$ and $\left\vert \nabla u\right\vert^q\in L^{2/q}(\Omega)$. This entails that the function $x\mapsto H(x,u(x),\nabla u(x))\in L^{p_0}(\Omega)$ with
$$p_0=\min\left(\frac{2n}{r(n-2)},\frac 2q\right)>\frac {2n}{n+2}.$$
Applying standard elliptic estimates and Lemma \ref{sobolevbis} a finite number of times yields the desired conclusion, as for the proof of item (iii) in Theorem~\ref{thapp}. The case where $n\leq 2$ is similar and easier. \hfill\fin


\subsection{Recovering Talenti's seminal results} \label{backtalenti}

In this section, we explain how to derive Talenti's results in~\cite{Talenti76} and~\cite{talenti} from Theorem \ref{th1} of the present work, under the assumption that the matrix field $A$ is continuous in~$\overline{\Omega}$. Let us concentrate on the result~(\ref{ineqtalenti}) under conditions~(\ref{hyptalenti}) and~(\ref{eqtalenti}) obtained in~\cite{talenti}, which clearly encompasses the one obtained in~\cite{Talenti76}. In the proof, $u$ denotes the unique $H^1_0(\Omega)$ solution of
\be\label{equtalenti}\left\{\baa{rcll}
-\div(A\nabla u)+\alpha\cdot\nabla u+b\,u & = & f & \hbox{in }\Omega,\vspace{3pt}\\
u & = & 0 & \hbox{on }\partial\Omega.\eaa\right.
\ee\par
First of all, if $\|f\|_{L^{\infty}(\Omega)}=0$, then, by uniqueness (Theorem~8.3 in~\cite{gt}), $u=0$ a.e. in $\Omega$ and the solution $v$ of~(\ref{eqtalenti}) is equal to $0$ a.e. in $\Omega^*$, whence the conclusion~(\ref{ineqtalenti}) is immediate.\par
One can then assume in the sequel that~$\|f\|_{L^{\infty}(\Omega)}>0$. Let $\widetilde{u}\in H^1_0(\Omega)$ be the unique solution of 
\be\label{tildeu}\left\{\baa{rcll}
-\div(A\nabla \widetilde{u})+\alpha\cdot \nabla \widetilde{u} & = & \left\vert f\right\vert & \hbox{in }\Omega,\vspace{3pt}\\
\widetilde{u} & = & 0 & \hbox{on }\partial\Omega.\eaa\right.
\ee
The weak maximum principle (Theorem~8.1 in~\cite{gt}) implies that $\tilde{u}\ge0$ a.e. in $\Omega$. Since $b\ge0$ a.e. in~$\Omega$, the function $\tilde{u}$ is then a weak supersolution of the equation~(\ref{equtalenti}) satisfied by $u$ as well as of the one satisfied by $-u$ (which is obtained by replacing $f$ by $-f$), whence  $u\leq \widetilde{u}$ and $-u\le\tilde{u}$ a.e. in $\Omega$. Therefore, $|u|\le\tilde{u}$ a.e. in $\Omega$ and
\be\label{uutildestar}
|u|^{\ast}\leq \widetilde{u}^{\ast}\hbox{ a.e. in }\Omega^{\ast}.
\ee\par
Pick up a sequence of matrix fields $(A_k)_{k\in\N}$ in $W^{1,\infty}(\Omega,{\mathcal S}_n(\R))$ converging to the matrix field $A$ in~$L^{\infty}(\Omega,{\mathcal S}_n(\R))$. There exists an increasing function $\varphi:\N\rightarrow \N$ such that, for all~$k\in\N$, $A_{\varphi(k)}-A\geq-2^{-k}\hbox{Id}$ in $\overline{\Omega}$. If $B_{k}:=A_{\varphi(k)}+2^{-k}\hbox{Id}$, then $B_{k}\rightarrow A$ as $k\to+\infty$ in~$L^{\infty}(\Omega,{\mathcal S}_n(\R))$ and $B_{k}\geq \hbox{Id}$ in $\overline{\Omega}$ for all $k\in\N$. Therefore, one can assume without loss of generality that $A_k\geq \hbox{Id}$ for all~$k\in\N$. Let also $(f_k)_{k\in\N}$ be a sequence of continuous functions in~$\overline{\Omega}$, such that $0\le f_k\le\|f\|_{L^{\infty}(\Omega)}$ in $\overline{\Omega}$ and $f_k\to|f|$ in $L^p(\Omega)$ as $k\to+\infty$ for all~$1\le p<+\infty$. Without loss of generality, one can assume that $\|f_k\|_{L^{\infty}(\Omega)}>0$ for all $k\in\N$.\par
On the one hand, for each $k\in\N$, let now $\widetilde{u}_k$ be the unique $H^1_0(\Omega)$ solution of
\be\label{tildeuk}\left\{\baa{rcll}
-\div(A_k\nabla\widetilde{u}_k)+\alpha\cdot\nabla\widetilde{u}_k & = & f_k & \hbox{in }\Omega,\vspace{3pt}\\
\widetilde{u}_k & = & 0 & \hbox{on }\partial\Omega,\eaa\right.
\ee
whose existence and uniqueness follow from Theorem~8.3 in~\cite{gt}. Since $\hbox{Id}\le A_k\le C\,\hbox{Id}$ in~$\overline{\Omega}$ for all $k\in\N$, where $C$ is a constant independent of $k$, and since $\alpha\in L^{\infty}(\Omega,\R^n)$ and the sequence $(f_k)_{k\in\N}$ is bounded in $L^{\infty}(\Omega)$, Corollary~8.7 in~\cite{gt} implies that the sequence~$(\widetilde{u}_k)_{k\in\N}$ is bounded in $H^1_0(\Omega)$. Up to extraction of a subsequence, there is then $\tilde{u}_{\infty}\in H^1_0(\Omega)$ such that~$\tilde{u}_k\to\tilde{u}_{\infty}$ as $k\to+\infty$ in $H^1(\Omega)$ weakly and in $L^2(\Omega)$ strongly. By testing~(\ref{tildeuk}) against an arbitrary test function $\varphi\in H^1_0(\Omega)$ and by passing to the limit as $k\to+\infty$, it easily follows that $\tilde{u}_{\infty}$ is a weak $H^1_0(\Omega)$ solution of~(\ref{tildeu}), whence $\tilde{u}_{\infty}=\tilde{u}$ by uniqueness and, furthermore, the whole sequence $(\tilde{u}_k)_{k\in\N}$ converges to $\tilde{u}$ in (at least) $L^1(\Omega)$. As a consequence, up to extraction of a subsequence, it follows from~\cite{ct} that
\be\label{tildeuuk}
\tilde{u}_k^*\to\tilde{u}^*\hbox{ as }k\to+\infty\hbox{ in }L^1(\Omega^*)\hbox{ and a.e. in }\Omega^*.
\ee\par
On the other hand, for each $k\in\N$, let $\overline{u}_k$ be the unique $H^1_0(\Omega)$ solution of
\be\label{overlineuk}\left\{\baa{rcll}
-\div(A_k\nabla\overline{u}_k)-\tilde{\alpha}\,|\nabla\overline{u}_k| & = & f_k & \hbox{in }\Omega,\vspace{3pt}\\
\overline{u}_k & = & 0 & \hbox{on }\partial\Omega,\eaa\right.
\ee
where $\tilde{\alpha}=\|\,|\alpha|\,\|_{L^{\infty}(\Omega)}$. The existence and uniqueness of $\overline{u}_k$ are guaranteed by Theorem~2.1 in~\cite{Porretta}. Furthermore, since $\alpha\cdot\nabla\overline{u}_k\ge-\tilde{\alpha}\,|\nabla\overline{u}_k|$ a.e. in $\Omega$, the function $\overline{u}_k$ is a weak supersolution of the equation~(\ref{tildeuk}), whence $\tilde{u}_k\le\overline{u}_k$ a.e. in $\Omega$ by the maximum principle (Theorem~8.1 in~\cite{gt}), and
\be\label{tildeoverlineuk}
\tilde{u}_k^*\le\overline{u}_k^*\hbox{ a.e. in }\Omega^{\ast}.
\ee\par
Furthermore, one can write the term $-\tilde{\alpha}\,|\nabla\overline{u}_k|$ in~(\ref{overlineuk}) as $-\tilde{\alpha}\,|\nabla\overline{u}_k(x)|=\overline{\alpha}_k(x)\cdot\nabla\overline{u}_k(x)$ in~$\Omega$, where the vector field $\overline{\alpha}_k\in L^{\infty}(\Omega,\R^n)$ can be defined as $\overline{\alpha}_k(x)=-\tilde{\alpha}\,\nabla\overline{u}_k(x)/|\nabla\overline{u}_k(x)|$ if~$|\nabla\overline{u}_k(x)|>0$ and~$\overline{\alpha}_k(x)=0$ if $|\nabla\overline{u}_k(x)|=0$. From elliptic regularity theory, one then has~$\overline{u}_k\in W(\Omega)$, while the Hopf lemma and the strong maximum principle show that $\overline{u}_k>0$ in $\Omega$ and $\partial_{\nu}\overline{u}_k<0$ on~$\partial\Omega$ (from Theorem~9.6 in~\cite{gt} and the discussion there, as in the end of the proof of Lemma~\ref{lembounds} above). For each $k\in\N$, the assumptions of Theorem~\ref{th1} are therefore satisfied by the function~$\overline{u}_k$, which solves an equation of the type~(\ref{equ}) with the matrix field $A_k\in W^{1,\infty}(\Omega,\mathcal{S}_n(\R))$ satisfying~(\ref{ALambda}) with $\Lambda=1$, and with the continuous function $H(x,s,p)=-\tilde{\alpha}\,|p|-f_k(x)$ satisfying~(\ref{hypH}) with $q=1$, $a(x,s,p)=\tilde{\alpha}$, $b(x,s,p)=0$ and~$f(x,s,p)=f_k(x)$. It follows then from the proof of Theorem~\ref{th1}, namely the inequa\-lity~(\ref{u*z}) in Remark~\ref{remu*z}, that
\be\label{overlineu*kzk}
\overline{u}^*_k\le z_k\hbox{ a.e. in }\Omega^*,
\ee
where $z_k$ is the unique $H^1_0(\Omega^*)$ solution of the equation corresponding to~(\ref{eqz}), that is here
$$\left\{\baa{rcll}
-\Delta z_k+\tilde{\alpha}\,e_r\cdot\nabla z_k & = & f_k^* & \hbox{in }\Omega^*,\vspace{3pt}\\
z_k & = & 0 & \hbox{on }\partial\Omega^*.\eaa\right.$$\par
Lastly, by~\cite{ct}, one can assume without loss of generality that $f^*_k\to|f|^*$ as $k\to+\infty$ a.e. in $\Omega^*$ (and notice that $\|f^*_k\|_{L^{\infty}(\Omega^*)}=\|f_k\|_{L^{\infty}(\Omega)}\le\|f\|_{L^{\infty}(\Omega)}$ for all $k\in\N$). Therefore, as for the solutions $\tilde{u}_k$ of~(\ref{tildeuk}), it follows easily that
\be\label{zkv}
z_k\to v\hbox{ as }k\to+\infty\hbox{ in }H^1(\Omega^*)\hbox{ weakly, }L^2(\Omega^*)\hbox{ strongly and a.e. in }\Omega^*,
\ee
where $v$ is the unique $H^1_0(\Omega^*)$ solution  of
$$\left\{\baa{rcll}
-\Delta v+\tilde{\alpha}\,e_r\cdot\nabla v & = & |f|^* & \hbox{in }\Omega^*,\vspace{3pt}\\
v & = & 0 & \hbox{on }\partial\Omega^*,\eaa\right.$$
that is~(\ref{eqtalenti}). As a conclusion, gathering~(\ref{uutildestar}),~(\ref{tildeuuk}),~(\ref{tildeoverlineuk}),~(\ref{overlineu*kzk}) and~(\ref{zkv}) yields $|u|^{\ast}\leq v$ a.e. in $\Omega^*$, which is the desired result~(\ref{ineqtalenti}).


\end{document}